\documentclass[12pt, leqno]{amsart}

\usepackage{youngtab}
\usepackage{graphicx}
\usepackage{amssymb,amsmath,amsthm,amscd}
\usepackage{mathrsfs,mathdots}
\usepackage{enumerate,stmaryrd}
\usepackage{upgreek}
\usepackage[usenames,dvipsnames]{color}
\usepackage[colorlinks=true, pdfstartview=FitV,
 linkcolor=blue,citecolor=blue,urlcolor=blue]{hyperref}

\usepackage[all]{xy}

\usepackage{verbatim}
\usepackage{oldgerm}
\usepackage{yhmath}
\usepackage[normalem]{ulem}  
\usepackage{xspace}
\usepackage{xcolor,cancel}
\usepackage{tikz}
\tikzset{dynkdot/.style={circle,draw,scale=.38}}

\usepackage{marginnote}

\numberwithin{equation}{section} \allowdisplaybreaks[4]

\usepackage[colorinlistoftodos]{todonotes}
\newcommand{\nc}{\newcommand}

\newcounter{myc}

\renewcommand{\le}{\leqslant}

\renewcommand{\ge}{\geqslant}
\renewcommand{\preceq}{\preccurlyeq}

\nc{\op}{\operatorname}

\theoremstyle{plain}
\newtheorem{lemma}{Lemma}[section]
\newtheorem{proposition}[lemma]{Proposition}
\newtheorem{theorem}[lemma]{Theorem}
\newtheorem{maintheorem}[lemma]{Main Theorem}

\newtheorem*{problem}{Problem}

\newtheorem{sublemma}{Sublemma}
\newtheorem{corollary}[lemma]{Corollary} 

\theoremstyle{definition}
\newtheorem{remark}[lemma]{Remark}

\newtheorem{definition}[lemma]{Definition}
\newtheorem*{convention}{Convention}

\nc{\Prop}{\begin{proposition}} \nc{\enprop}{\end{proposition}}
\nc{\Th}{\begin{theorem}} \nc{\enth}{\end{theorem}}
\nc{\Lemma}{\begin{lemma}} \nc{\enlemma}{\end{lemma}}
\nc{\Cor}{\begin{corollary}} \nc{\encor}{\end{corollary}}
\nc{\Rem}{\begin{remark}} \nc{\enrem}{\end{remark}}
\nc{\Def}{\begin{definition}} \nc{\edf}{\end{definition}}
\nc{\Sub}{\begin{sublemma}} \nc{\ensub}{\end{sublemma}}
\nc{\Prob}{\begin{problem}} \nc{\enprob}{\end{problem}}
\nc{\MT}{\begin{maintheorem}}
  \nc{\enMT}{\end{maintheorem}}

\nc{\qtext}{\quad\text}

\nc{\shc}{\mathcal{C}}

\newcommand{\Q}{\mathbb{Q}}

\newcommand{\frakC}{\mathfrak{C}}

\newcommand{\tx}{{\widetilde{x}}}
\newcommand{\ty}{{\widetilde{y}}}

\newcommand{\Seq}{\Sigma}
\newcommand{\T}{\mathcal{T}}
\newcommand{\dT}{\mathrm{T}}

\newcommand{\A}{\mathcal{A}}
\newcommand{\im}{i}
\newcommand{\jm}{j}

\newcommand{\Ca}{\mathscr{C}}
\nc{\F}{\mathcal{F}}

\newcommand{\D}{\mathscr{D}}
\newcommand{\Dd}{\mathcal{D}}

\nc{\HOM}{\on{H\textsc{om}}}
\newcommand{\M}{\mathcal{M}}
\newcommand{\W}{\mathsf{W}}

\newcommand{\LL}{\mathcal{L}}
\newcommand{\RR}{\mathcal{R}}

\newcommand{\bS}{\ms{1mu}{\mathsf{S}}}

\newcommand{\Z}{\mathbb{\ms{1mu}Z\ms{1mu}}}
\newcommand{\seteq}{\mathbin{:=}}
\newcommand{\hd}{{\operatorname{hd}}}
\newcommand{\soc}{{\operatorname{soc}}}
\nc{\ov}[1]{\overline{#1}} \nc{\wit}[1]{\widetilde{#1}}

\newcommand{\hI}{\widehat{I}}

\nc{\dis}[1]{\vert\ms{1.8mu}{#1}\ms{2mu}\vert}
\nc{\disp}[1]{\dis{#1}_\phi}

\nc{\Wlmj}[3]{\W_{#2,#3}^{(#1)}} \nc{\Mkl}[2]{\M_\ttww(#1,#2)}
\nc{\rmat}[1]{{\mathbf
r}_{\mspace{-2mu}\raisebox{-.5ex}{${\scriptstyle{#1}}$}}}
\nc{\po}[1]{p_{\mspace{-2mu}\raisebox{-.5ex}{${\scriptstyle{#1}}$}}}
\nc{\ddfrac}[2]{  \dfrac{#2}{#1} }
\newcommand{\on}{\operatorname}
\nc{\de}{\on{\textfrak{d}}} \nc{\tL}{\widetilde{\Lambda}}
\nc{\tl}{\widetilde{\lambda}} \nc{\mqs}{(-q^2)}
\nc{\Cquiver}{\upsigma}

\nc{\mut}[1]{{\mu}_{\mspace{-2mu}\raisebox{-.5ex}{${\scriptstyle{#1}}$}}}

\newcommand{\g}{{\mathfrak{g}}}
\newcommand{\fing}{{\mathsf{g}}}
\newcommand{\finn}{{\mathsf{n}}}

\newcommand{\Anw}{A_q(\finn(w))}

\newcommand{\n}{{\mathfrak{n}}}
\newcommand{\isoto}[1][]{\mathop{\xrightarrow%
[{\raisebox{.3ex}[0ex][.3ex]{$\scriptstyle{#1}$}}]%
{{\raisebox{-.6ex}[0ex][-.6ex]{$\mspace{2mu}\sim\mspace{2mu}$}}}}}

\newcommand{\soplus}{\mathop{\mbox{\normalsize$\bigoplus$}}\limits}

\newcommand{\ww}{ \textbf{\textit{w}}}
\newcommand{\ttww}{{\widetilde{\ww}} }

\nc{\nconv}{\mathop{\mbox{\large $\odot$}}}
\nc{\nnconv}{\mathop{\mbox{\large $\star$}}}

\nc{\lb}{\llbracket} \nc{\rb}{\rrbracket}

\nc{\la}{\lambda} \nc{\La}{\Lambda}

\nc{\finite}{{\mathrm{fin}}} \nc{\gf}{{\g_\finite}}
\nc{\tLa}{\widetilde{\Lambda}} \nc{\ve}{\varepsilon}
\nc{\ep}{\epsilon} \nc{\vp}{\varphi} \nc{\lan}{\langle}
\nc{\ran}{\rangle} \nc{\Uqg}{U_q(\g)} \nc{\Aqg}{A_q(\g)}
\nc{\Aqn}{A_q(\n)} \nc{\ual}{\upalpha\ms{1mu}}
\nc{\sPi}{\mathsf{\Pi}}
\nc{\Pif}{\Pi_\finite} 
\nc{\sLa}{\mathsf{\Lambda}\ms{1mu}} \nc{\al}{\alpha} \nc{\be}{\beta}
\nc{\ga}{\gamma} \nc{\wt}{\operatorname{wt}}
\nc{\ch}{\operatorname{ch}}

\nc{\norm}{{\mathrm{norm}}} \nc{\aff}{{\mathrm{aff}}}
\nc{\Maf}{M_\aff} \nc{\ev}{{\mathrm{even}}} \nc{\od}{{\mathrm{odd}}}
\nc{\Sev}{\Seq^{\ev}} \nc{\Sod}{\Seq^{\od}} \nc{\Spl}{\Seq^{+}}
\nc{\Smi}{\Seq^{-}} \nc{\low}{{\mathrm{low}}}
\nc{\upper}{{\mathrm{up}}} \nc{\one}{{\bf{1}}}
\nc{\To}[1][{\hspace{2ex}}]{\xrightarrow{\,#1\,}}
\nc{\te}{\tilde{e}} \nc{\tw}{{\underline{w}}}
\nc{\hw}{{\widehat{w}_0}} \nc{\hhw}{\widehat{\tw}_0}
\nc{\tww}{\ww} \nc{\tuu}{{\mathsf{u}}} \nc{\tel}{\tilde{e}^\low}
\nc{\teu}{\tilde{e}^\upper} \nc{\tf}{\tilde{f}}
\nc{\tfl}{\tilde{f}^\low} \nc{\tfu}{\tilde{f}^\upper}
\nc{\tE}{\widetilde{E}} \nc{\tF}{\widetilde{F}}
\nc{\tFF}{\widetilde{\F}} \nc{\tB}{\widetilde{B}}
\nc{\tz}{\tilde{z}}
\nc{\tQ}{\hspace{-.2ex}\textbf{\textit{Q}}}
\nc{\Ft}{\F^\dT}
\nc{\Seed}{\mathcal{S}}
\newenvironment{rouge}
{\color{red}} {}

\nc{\ber}{\begin{rouge}} \nc{\er}{\end{rouge}}

\newcommand{\berm}{\begin{rouge}{}\marginnote{\fbox{\scshape\lowercase{M}}}
}
\newcommand{\berMH}{\begin{rouge}{}\marginnote{\fbox{\scshape\lowercase{MH}}}
}

\newcommand{\bero}{\begin{rouge}{}\marginnote{\fbox{\scshape\lowercase{O}}}
{}}

\newenvironment{blue}
{\color{Dandelion}} {}

\nc{\beb}{\begin{rouge}\marginnote{\fbox{\scshape\lowercase{M}}}\end{rouge}
\begin{blue}} \nc{\eb}{\end{blue}}

\nc{\ble}[1]{\underline{#1}}

\newenvironment{bluk}
{\relax\color{blue}} {\hspace*{.5ex}\relax}

\newcommand{\bek}{\begin{bluk}}
\newcommand{\ek}{\end{bluk}}

\newenvironment{yellow}
{\relax\color{Dandelion}}
{\hspace*{.5ex}\relax}

\newcommand{\bey}{\begin{yellow}}
\newcommand{\ey}{\end{yellow}}

\nc{\cor}{\mathbf{k}} \nc{\tens}{\mathop\otimes}
\nc{\gmod}{\mbox{-$\mathrm{gmod}$}}
\nc{\gMod}{\mbox{-$\mathrm{gMod}$}} \nc{\Md}{\mbox{-$\mathrm{Mod}$}}
\nc{\md}{\mbox{-$\mathrm{mod}$}} \nc{\uqm}{\mathscr{C}_\g}

\nc{\proj}{\mbox{-$\mathrm{proj}$}}
\nc{\gproj}{\mbox{-$\mathrm{gproj}$}}
\nc{\smod}{\mbox{-$\mathrm{mod}$}}
\nc{\nmod}{\mbox{-$\mathrm{nilmod}$}}

\nc{\seed}{\mathscr{S}}

\newcommand{\cmf}{\mathsf{C}} 

\nc{\Rnorm}{\mathrm{R}^{\mathrm{norm}}}
\nc{\Runiv}{\mathrm{R}^{\ms{1mu}\mathrm{univ}}}
\nc{\Rren}{\mathrm{R}^{\ms{1mu}\mathrm{ren}}} 
\nc{\col}{\colon} \nc{\epiTo}[1][]{\xymatrix{\ar@{->>}[r]^-{{#1}}&}}
\nc{\epito}{\twoheadrightarrow}
\nc{\monoTo}[1][]{\xymatrix{\ar@{>->}[r]^-{{#1}}&}}
\nc{\monogets}[1][]{\xymatrix{&\ar@{_{(}->}[l]^-{{#1}}}}
\nc{\sym}{\mathfrak{S}} \nc{\rl}{\mathsf{Q}} \nc{\prl}{\rl_+}
\nc{\crl}{\mathsf{Q}^\vee} \nc{\pcrl}{\crl_+} \nc{\Qq}{{\Q(q)}}
\nc{\wl}{\mathsf{P}}   
\nc{\wlf}{\mathsf{P}_\finite}   
\nc{\Oint}{\mathcal{O}_{{\mathrm{int}}}}
\newcommand{\scbul}{{\,\raise1pt\hbox{$\scriptscriptstyle\bullet$}\,}}

\nc{\conv}{\mathop{\mathbin{\mbox{\large $\circ$}}}}
\newcommand{\hconv}{\mathbin{\scalebox{.9}{$\nabla$}}}

\nc{\scrA}{\mathscr{A}}

\nc{\pv}{  \to\updownarrow\gets }
\nc{\nv}{  \longleftrightarrow {\raise -1pt\hbox{$\hspace{-2ex}\begin{matrix}\downarrow \\[-1ex] \uparrow\end{matrix}$}} }

\newcommand{\Hom}{\operatorname{Hom}}

\renewcommand{\Im}{\op{Im}}

\newcommand{\ex}{{\mathrm{ex}}}
\nc{\K}{\mathsf{K}} \nc{\Kex}{{\K}^{\mathrm{ex}}}
\nc{\Uex}{\Uppsi_{\mathrm{ex}}}
\nc{\Kfr}{{\K}^{\mathrm{f\mspace{.01mu}r}}}
\nc{\ben}{\begin{enumerate}}
\nc{\ee}{\end{enumerate}} 
\nc{\bnum}{\begin{enumerate}[{\rm(i)}]}
\nc{\bna}{\begin{enumerate}[{\rm(a)}]} 
\nc{\bc}{\begin{cases}}
\nc{\ec}{\end{cases}}

\newenvironment{myequation}
{\relax\setlength{\arraycolsep}{1pt}\begin{eqnarray}}
{\end{eqnarray}}
\newenvironment{myequationn}
{\relax\setlength{\arraycolsep}{1pt}\begin{eqnarray*}}
{\end{eqnarray*}}

\nc{\eq}{\begin{myequation}} \nc{\eneq}{\end{myequation}}
\nc{\eqn}{\begin{myequationn}} \nc{\eneqn}{\end{myequationn}}
\nc{\hs}{\hspace*} \nc{\hl}{\hspace{-.5ex}}

\newenvironment{myarray}[1]{\relax\setlength{\arraycolsep}{.5pt}

\begin{array}{#1}}{\end{array}\relax}

\newcommand{\ba}{\begin{myarray}}
\newcommand{\ea}{\end{myarray}}

\nc{\noi}{\noindent} \nc{\ang}[1]{\langle{#1}\rangle}
\nc{\fr}{{\mathrm{fr}}} \nc{\qt}[1]{\quad\text{#1}}
\nc{\ol}{\overline} \nc{\true}{\delta} \nc{\ms}{\mspace}

\nc{\vs}{\vspace*} \nc{\bl}{\bigl(} \nc{\br}{\bigr)}
\nc{\bep}{\ol{\ep}} \nc{\bal}{\,\ol{\al}}
\nc{\qtq}[1][{and}]{\quad\text{#1}\quad}
\nc{\set}[2]{\{{#1}\bigm|{#2}\}}
\nc{\lset}[2]{\left\{{#1}\bigm|{#2}\right\}}
\nc{\rmo}{{\rm(}}
\nc{\rmf}{{\rm)}\xspace} \nc{\Proof}{\begin{proof}}
\nc{\QED}{\end{proof}}
\nc{\monoto}[1][]{\xymatrix@C=2ex{\ar@{>->}[r]^-{{#1}}&}\ms{-8mu}}
\nc{\etens}{\boxtimes} \nc{\height}[1]{\vert{#1}\vert}
\nc{\Lrev}{L^{\bek{\rev}\ek}}
\nc{\rev}{\mathrm{rev}} \nc{\fw}{\Lambda} \nc{\uqpg}{U_q'(\g)}
\nc{\tp}{\ms{1.5mu}{\widetilde{p}}\ms{2mu}}
\nc{\Deg}{\mathrm{\ms{1mu}Deg\ms{1mu}}} \nc{\Bg}{\mathcal{G}}

\nc{\wb}[1]{\mbox{$\rule[-1.1ex]{0ex}{2ex}#1$}}

\nc{\bwr}{\mbox{\large$\wr$}} \nc{\vphi}{\varphi}
\nc{\G}{\mathcal{G}} \nc{\tD}{\widetilde{\mathrm{De}}\mathrm{g}}
\nc{\Li}{\La^\infty} \nc{\Di}{\Deg^\infty}
\nc{\zero}{\ms{2mu}\mathrm{zero}\ms{2mu}} \nc{\cwl}{\wl^\vee}
\nc{\rc}{renormalizing coefficient\xspace}
\nc{\cz}{{\cor[z^{\pm1}]}}
\nc{\ake}[1][2ex]{\rule[-.5ex]{0ex}{#1}}
\nc{\akete}[1][0ex]{\rule[{#1}]{0ex}{1ex}}
\nc{\akew}[1][2ex]{\rule[-1ex]{#1}{0ex}} \nc{\rd}{{}^*\ms{-3mu}}
\nc{\st}[1]{\{{#1}\}}
\nc{\seq}[1]{ \boldsymbol{(} {#1} \boldsymbol{)}   }
\nc{\corh}{\widehat{\cor}}
\nc{\czt}{\cz^\times} \nc{\eps}{\varepsilon} \nc{\rr}{rationally
renormalizable\xspace} \nc{\QHA}{\mathrm{QHA}} \nc{\Ker}{\on{Ker}}
\nc{\usq}{U'_q(\g)} \nc{\Wf}{\W_\finite} \nc{\If}{I_\finite}
\nc{\ord}{\mathrm{ord}}
\nc{\Qd}{\mathscr{Q}} \nc{\e}{\mathrm{e}} \nc{\snoi}{\smallskip\noi}
\nc{\mnoi}{\medskip\noi}
\nc{\ci}{{\ms{1mu}\mathfrak{c}\ms{1mu}}}
\newcommand{\tc}{{\widetilde{\ci}}}
\nc{\Iff}{I_\fing} \nc{\CM}[1][s]{M[{#1},0\}}
\nc{\CMp}[1][s]{M[{#1},0\}'} \nc{\cM}[1][s]{\mathsf{M}_{#1}}
\nc{\fM}{\mathsf{M}}
\nc{\cMp}[1][s]{\mathsf{M}'_{#1}}
\nc{\Vi}{\mathrm{Vi}} 
\nc{\Vo}{\mathrm{Vo}}
\nc{\GLS}{\mathrm{Q}_{\mathrm{GLS}}}
\nc{\HL}{\mathrm{Q}_{\mathrm{HL}}} \nc{\vpi}{\varpi}
\nc{\bpi}{\ol{\pi}}
\nc{\dC}{\ms{1mu}\mathsf{h}^\vee} 
\nc{\ca}{completely $\Uplambda$-admissible\xspace}
\nc{\KR}{Kirillov-Reshetikhin\xspace}
\nc{\qA}{quasi-admissible\xspace}
\nc{\Ad}{affine determinantial\xspace}
\nc{\AD}{Affine determinantial\xspace}
\nc{\SW}{Schur-Weyl\xspace}
\nc{\xiz}{\xi^{\mathrm{zz}}} \nc{\Kin}{\K_{\mathrm{in}}}
\nc{\Kb}{\K_{\mathrm{b}}} \nc{\Kout}{\K_{\mathrm{out}}} \nc{\rE}{
\mathsf{E} } \nc{\rW}{ \mathcal{W} } \nc{\catCO}{\Ca^0_\g}
\nc{\sig}{\upsigma(\g)} \nc{\sigZ}{\upsigma_0(\g)}
\nc{\cato}{\catCO} \nc{\cat}{\Ca} \nc{\catg}{\Ca_\g}
\nc{\Set}{{\mathrm{Set}}} \nc{\BHL}{\tB_{\mathrm{HL}}}
\nc{\R}{\mathbb{R}} \nc{\BGLS}{\tB_{\mathrm{GLS}}}
\nc{\hIg}{{\hI_\g}} \nc{\KK}{;\K,\Kex}
\nc{\KKC}{;\K(\frakC),\Kex(\frakC)} \nc{\Fd}{\F_\Dd}
\nc{\KKCp}{;\K(\frakC'),\Kex(\frakC')}
\nc{\hF}{\widehat{\F}} \nc{\Rc}{R_\cmf} \nc{\pbw}{PBW-pair\xspace}
\nc{\pbws}{PBW-pairs\xspace} \nc{\Refl}{\mathscr{S}}
\nc{\Reflinv}{{\Refl}^{-1}} \nc{\Rt}{\mathsf{L}}
\nc{\qa}{$\Uplambda$-quasi-admissible\xspace}
\nc{\Lad}{$\Uplambda$-admissible\xspace}
\nc{\htens}{\hconv}
\nc{\Isf}{I_\fing}
\nc{\scb}{\scalebox}
\nc{\vrho}{\varrho}
\nc{\dsum}{\displaystyle\sum}
\nc{\dtens}{\mathop{\mbox{\scb{1.1}{\akete[-.6ex]\normalsize$\bigotimes$}}}\limits}
\nc{\cq}{\Check{q}}
\nc{\ud}[1]{\underline{#1}}

\nc{\sfc}{\mathsf{c}}
\nc{\sfa}{\mathsf{a}}
\nc{\bg}{\mathbf g}
\nc{\weyl}{\mathsf{W}}
\nc{\sL}{\mathsf L}
\nc{\cartan}{\mathsf{C}}
\nc{\afr}{affreal\xspace}

\nc{\infseq}{\bold i}

\renewcommand{\emptyset}{\varnothing}

\nc{\efe}{effective end\xspace}
\nc{\efes}{effective ends\xspace}
\nc{\ibox}{$\mathbf{i}$-box\xspace}
\nc{\iboxes}{$\mathbf{i}$-boxes\xspace}
\nc{\mcf}[1][{$[a,b]$}]{maximal commuting family of \iboxes in #1\xspace}
\nc{\mcfs}[1][{$[a,b]$}]{maximal commuting families of \iboxes in #1\xspace}
\nc{\adc}[1][{$[a,b]$}]{admissible chain of \iboxes in #1\xspace}
\nc{\ac}{admissible chain\xspace}
\nc{\Cw}{\cat_w}
\nc{\cl}{\colon}
\nc{\Fex}{\F_\ex}
\nc{\bM}{\mathsf{M}}
\nc{\hB}{\widehat{B}}
\nc{\ttP}{\mathsf{P}}
\nc{\bout}[1]{\ber\sout{#1}\er}

\newenvironment{magem}{\relax\color{magenta}}{\relax}

\newcommand{\bema}{\begin{magem}}
\newcommand{\ema}{\end{magem}}

\newcommand{\RD}{{{\mathscr{D}}}}

\newlength{\mylength}
\title[Exchange matrices  of I-boxes]{Exchange matrices  of I-boxes}

\author[M. Kashiwara]{Masaki Kashiwara}
\thanks{The research of M.\ Kashiwara was supported by Grant-in-Aid for Scientific Research (B) 23K20206, Japan Society for the Promotion of Science.}
\address[M. Kashiwara]{
 Kyoto University Institute for Advanced Study, Research Institute for Mathematical Sciences, Kyoto University, Kyoto 606-8502, Japan}
\email[M. Kashiwara]{masaki@kurims.kyoto-u.ac.jp}

\author[M. Kim]{Myungho Kim}
\address[M. Kim]{Department of Mathematics, Kyung Hee University, Seoul 02447, Korea}
\email[M. Kim]{mkim@khu.ac.kr}
\thanks{The research of M.\ Kim was supported by the National Research Foundation of
Korea (NRF) Grant funded by the Korea government(MSIT)
(NRF-2022R1F1A1076214 and NRF-2020R1A5A1016126).}



 \date{November 26, 2025}

\begin{document}

\begin{abstract} 
Admissible chains of 
$\bold i$-boxes are important combinatorial tools in the monoidal categorification of cluster algebras,
as they provide seeds of the cluster algebra. 
In this paper, we explore the properties of maximal commuting families of  $\bold i$-boxes in a more general setting,  and define a certain matrix associated with such a family, which we call the exchange matrix.
It turns out that, when considering the cluster algebra structure on the Grothendieck rings, this matrix is indeed the exchange matrix of 
the seed  associated with the family, both in certain categories of modules over quantum affine algebras and over quiver Hecke algebras.
We prove this by constructing  explicit short exact sequences that represent the mutation relations.
\end{abstract}

\maketitle
\tableofcontents

\section{Introduction}

A cluster algebra is a subalgebra of a field of rational functions, generated by a special set of elements called  \emph{cluster variables}.
Since  the introduction by Fomin and Zelevinsky in \cite{FZ02},
cluster algebras have been studied extensively in connection with many areas of mathematics.
The cluster variables are grouped into overlapping subsets known as \emph{clusters}. 
A \emph{cluster monomial} is a monomial of cluster variables in a cluster. 
Note that  a \emph{seed},  the pair consisting of a cluster and its \emph{exchange matrix},  can be obtained from an initial seed through a sequence of  inductive procedures known as \emph{mutations}  (\cite{FZ02}).
In \cite{HL10}, Hernandez and Leclerc introduced the notion of \emph{monoidal categorification of cluster algebras} in their study of finite-dimensional modules over a quantum affine algebra $\uqpg$.
Let  $\cat_\g$ be the category of finite-dimensional integrable modules over $\uqpg$. 
It is revealed  in \cite{KKOP22} that there are many monoidal subcategories $\cat^{\mathcal D, \widehat {\underline w}_0,[a,b]}_{\g}$  of $\cat_\g$ which serve as  monoidal categorifications of
 the cluster algebra structure on their Grothendieck rings.
For an interval $[a,b]$ in $\Z$, the category  $\cat^{\mathcal D, \widehat {\underline w}_0,[a,b]}_{\g}$ is defined as the smallest full subcategory of $\cat_\g$ that is stable by taking subquotients,  extensions and tensor products,  and contains the \emph{affine cuspidal modules} $S^{\mathcal D,\widehat {\underline w}_0 }_k$ $(a\le k\le b)$. 
Recall that  by choosing  a family ${\mathcal D}$ of simple modules in $\cat_\g$, called a \emph{strong duality datum associated with $\bg$},  we obtain a  functor $\F_{\mathcal D}$ from $ R^\bg\gmod$ to $\cat_\g$, where $\bg$ is a simply-laced finite type  Lie algebra,  
and $R^\bg\gmod$ is the category of finite-dimensional graded modules over the quiver Hecke algebra $R^\bg$ of type $\bg$ (see \cite{KKOP20C}).
For each reduced expression $\underline w_0 = s_{i_1}\cdots s_{i_r}$ of the longest element $w_0$ of the Weyl group of $\bg$, there exists a distinguished family $\{S^{\bg, \underline w_0}_k\}_{1\le k\le r}$ of simple modules in  $R^\bg\gmod$, called the \emph{cuspidal modules}. 
Then for each $k\in \Z$ we obtain the affine cuspidal modules by setting 
$S^{\mathcal D,\widehat {\underline w}_0 }_k\seteq \F_{\mathcal D}(S^{\bg, \underline w_0}_k)$ for  $1\le k\le r$, and $S^{\mathcal D,\widehat {\underline w}_0 }_{k+r}\seteq  \RD(S^{\mathcal D,\widehat {\underline w}_0 }_k)$ for $k\in\Z$, where $\D(X)$  denotes the right dual of a module $X$ and $\widehat {\underline w}_0 $ denotes a specific infinite sequence indexed by $\Z$ extending $\underline w_0$.

The category $\cat^{\mathcal D, \widehat {\underline w}_0,[a,b]}_{\g}$ provides a monoidal categorification of a cluster algebra (\cite{KKOP22}),  that is,  the Grothendieck ring $K(\cat^{\mathcal D, \widehat {\underline w}_0,[a,b]}_{\g})$ has a cluster algebra structure, and the cluster monomials correspond to simple modules in the category.  
Note that the mutation of cluster variable is realized as a short exact sequence within the category.
Recall that  two simple modules $M$ and $N$ are said to \emph{strongly commute} if  the tensor product $M\tens N$ is simple. Every cluster in $K(\cat^{\mathcal D, \widehat {\underline w}_0,[a,b]}_{\g})$  corresponds to a family of simple modules in  the category that  strongly commute with each other. 
Although  any cluster can be obtained from the initial cluster after a sequence of mutations,
it is difficult to describe all of them explicitly.

\medskip
However, there are notable families of strongly commuting simple modules that can be described very explicitly: the \emph{affine determinantial modules associated with an admissible chain of \iboxes with extent $[a,b]$} (\cite[Section 4, 5]{KKOP22}). 
Let us recall these notions more precisely.
Let $\bold i = (i_k)_{k\in\Z} $ be a certain infinite sequence  $ \widehat {\underline w}_0$
which is an extension of $\underline w_0$.
We call $i_k$ the \emph{color of $k$} and 
we  call an interval $[x,y]$ in $\Z$ such that  $i_x=i_y$  an \emph{$\bold i$-box}.  
Then the affine determinantial module $M^{\mathcal D,\widehat {\underline w}_0 }[x,y]$ associated with the $\bold i$-box $[x, y]$ is defined as the head of the decreasingly ordered tensor product of all affine cuspidal modules within the interval $[x,y]$ that share the same color as $x$.
Specifically, we have  
$M^{\mathcal D,\widehat {\underline w}_0 }[x,y]= \hd(S^{\mathcal D,\widehat {\underline w}_0 }_y \tens S^{\mathcal D,\widehat {\underline w}_0 }_{y_-}  \tens \cdots \tens S^{\mathcal D,\widehat {\underline w}_0 }_{x_+} \tens  S^{\mathcal D,\widehat {\underline w}_0 }_{x} )$, where $z_{\pm}$ denotes the integers adjacent  to $z$ that have the same color as $z$ (see \eqref{eq: nota +,-}). 
We say that two \iboxes $[x,y]$ and $[x',y']$ \emph{commute} if either 
$x_-<x'\le y'<y_+$ or $ x'_-<x\le y<y'_+$ (see Definition~\ref{def:com}). 
In other words, two \iboxes are said to commute if the extension of
one of the \iboxes by one step to the left and  one step to the right 
properly contains the other $\bold i$-box. 
Remarkably, this simple combinatorial condition on a pair of  \iboxes ensures that the corresponding pair of the affine determinantial modules strongly commute with each other. 

An admissible chain of \iboxes with extent $[a,b]$ can be described by inductively  constructing  an increasing sequence of intervals $\{ \tc_k \}_{1 \le k \le b-a+1}$, which are referred to as the \emph{envelopes} of the \iboxes.
We start by choosing an $\bold i$-box $\ci_1=[x_1,x_1]$ for some $x_1\in [a,b]$.
 We set the envelope $\tc_1$ of $\ci_1$ to be the interval $\ci_1$ itself.
To obtain the next envelope $\tc_2$, we have two choices:  either enlarge $\tc_1$ to the left by $1$ unit or to the right by $1$ unit. 
In each case, we determine the corresponding $\bold i$-box $\ci_2=[x_2,y_2]$ by finding the largest $\bold i$-box in the envelope $\tc_2$ containing the newly created integer in $\tc_2 \setminus \tc_1$.
By repeating this procedure until we reach $\tc_{b-a+1}=[a,b]$, we obtain a family of \iboxes  $\frakC=\{\ci_k=[x_k,y_k]\}_{1\le k \le b-a+1}$ along with  a  sequence of envelopes
$\{\tc_k\}_{1\le k \le b-a+1}$.
It turns out that any pair of \iboxes in $\frakC$ commute with each other.
Moreover  the family $\frakC$ is a maximal commuting family of \iboxes in $[a,b]$.
Hence we obtain a commuting family $\{M^{\mathcal D,\widehat {\underline w}_0 }[x_k,y_k]\}_{1\le k\le b-a+1}$ of simple modules in the category $\cat^{\mathcal D, \widehat {\underline w}_0,[a,b]}_{\g}$.
If we assume further that the duality datum $\mathcal D$ arises from a $\rm Q$-datum, that is $\mathcal D=\mathcal D_{\mathcal Q}$ for a $\rm Q$-datum $\mathcal Q$ (see \cite[Section 6]{KKOP22}), then this family  forms a cluster in the cluster algebra $K(\cat^{\mathcal D_Q, \widehat {\underline w}_0,[a,b]}_{\g})$ (\cite[Theorem 8.1]{KKOP22}).

There is a useful procedure called the \emph{box move}, which produces a new admissible chain of \iboxes from a given one.
Consider the set $\{ \tc_k\}_{1\le k\le  b-a+1}$ of envelopes associated with an admissible chain $\frakC$. 
This set consists of intervals within $[a,b]$ that are totally ordered by inclusion.
Conversely, any collection of $b-a+1$ intervals in $[a,b]$ that are totally ordered by inclusion can serve as the set of envelopes of some admissible chain of \iboxes in $[a,b]$.
When it is possible to \emph{move} an envelope $\tc_k$ to the right or to the left, we define $B_k(\frakC)$ at $k$ to be the admissible chain of \iboxes
obtained by such a move of $\tc_{k}$. 
The operation $B_k$ is referred to as the {\em box move}  at $k$
(see the paragraph below Proposition \ref{prop:adch}). 
A box move either  permutes  the \iboxes in $\frakC$ or replaces an $\bold i$-box in $\frakC$ with a new $\bold i$-box that was not in $\frakC$. 
It turns out that the latter case corresponds to the mutation of a cluster variable in the cluster algebra $K(\cat^{\mathcal D, \widehat {\underline w}_0,[a,b]}_{\g})$.
The short exact sequence in the category representing this mutation is referred to as the \emph{$T$-system}.

\medskip
The concept of an admissible chain of \iboxes was crucial in demonstrating that the category $\cat^{\mathcal D_{\mathcal Q}, \widehat {\underline w}_0,[a,b]}_{\g}$ serves as a monoidal categorification of its Grothendieck ring.
 However, an interesting problem remains unsolved in \cite{KKOP22}: while the cluster variables, i.e., the affine determinantial modules associated with  the \iboxes in the family $\frakC$, are described explicitly as mentioned above,  the exchange matrix associated with this family has not been  explicitly determined there.
Note that the exchange matrix is important because it encodes how to obtain the new cluster variables from the cluster through mutations.
Although the exchange matrix can be obtained through a sequence of mutations associated with box moves, 
this does not mean that the matrix is explicitly described. 
In this paper, we present an explicit form of the exchange matrix 
 associated with the admissible chain of \iboxes.

\bigskip

We shall now explain the results  of this paper.
We explore the concepts of \iboxes and admissible chains of \iboxes for \emph{arbitrary sequences} $\bold i$,
not necessarily arising from a reduced expression of a Weyl group element.
For example, 
it is true in general that an admissible chain of \iboxes with extent $[a,b]$ forms a maximal commuting family of \iboxes in $[a,b]$.

Our main %
focus is on
 maximal commuting families of \iboxes with extent $[a,b]$.
We show that, for every maximal commuting family $\F$of \iboxes in an interval $[a,b]$, there exists an admissible chain of \iboxes
$\frakC=\{\ci_k\}_{1\le k \le b-a+1}$ such that $\F=\st{\ci_k\mid 1\le k \le b-a+1}$.
This allows us to study the maximal commuting families of \iboxes through the lens of admissible chains.
Additionally, we develop and study the notion of the \emph{\efe} of an $\bold i$-box $[x,y]$ in a maximal commuting family $\F$ of \iboxes.
Note that there can be several admissible chains of \iboxes that result in the same family $\F$.
Nevertheless, there exists a unique $z\in \{x,y\}$ such that for any admissible chain $\frakC$ that  results in $\F$, if $[x,y] = \ci_k$ for some $k$, then $\{z\} = \tc_k \setminus \tc_{k-1}$.
We call $z$ the {\em \efe } of $[x,y]$ in $\F$. The \efe enables us to study and describe many properties of the maximal commuting family $\F$ without invoking  admissible  chains.
\medskip

We define the \emph{exchange matrix} $\tB(\F)$ for a maximal commuting family $\F$ of \iboxes, which is one of the main contributions of this paper.
To this end, we consider the case where $\bold i$ is a sequence within the index set $I$  of 
a symmetrizable Cartan matrix $\cartan=(\sfc_{i,j})_{i,j\in I}$ (not necessarily of  simply-laced finite type).
For a maximal commuting family $\F$ of \iboxes,
we define the exchange matrix of $\F$ as a skew-symmetrizable
matrix $\hB(\F)=(b_{[x,y],[x',y']})_{[x,y],[x',y']\in \F}$. 
See \eqref{eq:ex_matrix} for its precise description. 
We want to emphasize that the entry $b_{[x,y],[x',y']}$ is determined by some local information surrounding $[x,y]$ and $[x',y']$ within $\F$, making its calculation almost immediate once the pair $[x,y]$ and $[x',y']$  is given.
Moreover, every $b_{[x,y],[x',y']}$ belongs to
  $\st{1,-1,\sfc_{i_x,i_{x'}},-\sfc_{i_x,i_{x'}},0}$. 
  
One may associate a quiver to the exchange matrix,  which consists of two types of arrows: the \emph{horizontal arrows}  connecting \iboxes with the same colors, and the  \emph{vertical arrows} connecting the \iboxes with different colors. 
Note that, up to this point, the matrix $\hB(\F)$ has only been defined and is not yet related with a  cluster algebra.

We further analyze the exchange matrix in detail.
The horizontal arrows are readily understood by the definition.
Note that for each $i$ in $I$, there exists a largest $\bold i$-box in $\F$ of color $i$. These \iboxes in $\F$ are called \emph{frozen}, and while the remaining \iboxes are called \emph{exchangeable}.
We denote by $\F_\fr$ the set of frozen \iboxes and by
$\F_\ex$ the set of exchangeable \iboxes. 

For each exchangeable $\bold i$-box $[x,y]$,  we characterize the sets of  \iboxes $[x',y']$ that have incoming vertical arrows to 
$[x,y]$  and those that have outgoing vertical arrows from $[x,y]$.  
This characterization depends on the various configurations of the horizontal arrows adjacent to $[x,y]$, requiring a case-by-case study.
It is important to characterize these subsets because, if $\F$ is associated with a cluster in a cluster algebra, then the products of variables within these subsets correspond to the two monomials in the mutation relation of  the cluster variable corresponding to $[x,y]$.

\medskip 
Next, we relate monoidal categories to the maximal commuting families of \iboxes.
To this end, in addition to  the category $\cat^{\mathcal D_{\mathcal Q}, \widehat {\underline w}_0,[a,b]}_{\g}$ mentioned earlier, we consider another  class of  monoidal categories,  denoted by $\cat_w$.
 This category $\cat_w$  is a full subcategory of $R^\bg \gmod$ associated with $w$.  Here  $\bg$ is a symmetrizable Kac-Moody algebra,  $R^\bg\gmod$ is the category of finite-dimensional graded modules over the quiver Hecke algebra of type $\bg$, and $w$ is an element of the Weyl group of $\bg$.  
By selecting a reduced expression of $w$, we obtain the set of \emph{cuspidal modules}, which generates the category $\cat_w$ within $R^\bg \gmod$.
Recall that the Grothendieck ring of $K(\cat_w)$ is known to be isomorphic to the \emph{quantum unipotent coordinate ring} $A_q(\n(w))$ (\cite{KL09, R08}),  and it has a quantum cluster algebra structure (\cite{GY17}).
Let the category $\cat$ be either $\cat_w$ or $\cat^{\mathcal D_{\mathcal Q}, \widehat {\underline w}_0,[a,b]}_{\g}$. The remaining discussions  are valid in both cases. 
For the case  $\cat=\cat_w$,    we take $\bold i$
to be a reduced expression of $w$, and for the case  $\cat=\cat^{\mathcal D_{\mathcal Q}, \widehat {\underline w}_0,[a,b]}_{\g}$,  we take $\bold i$ to be the sequence $\widehat {\underline w}_0$ mentioned earlier. 
Let $\F$ be a maximal commuting family of \iboxes. 
In the case $\cat=\cat_w$, we can also  associate $\F$ with a family of \emph{determinantial modules}, which strongly commute with each other in $\cat_w$ (\cite{KKOP24}).
The main result of this paper is Theorem \ref{thm:main}.  It establishes that the pair $\bl[\F],\tB(\F)\br$ is a seed of the cluster algebra $K(\cat)$ (Theorem \ref{thm:main} \eqref{main_ii}), where $\tB(\F)$ is the restriction of $\hB(\F)$
to $\F\times\F_\ex$. 
  It is proved by constructing short exact sequences in $\cat$
  which represent mutations associated with $\tB(\F)$.
Remark that, although we know that the Grothendieck ring $K(\cat)$
has a cluster algebra structure,
we only know that $\cat$ is a  monoidal categorification of the cluster algebra $K(\cat)$ in the case $\cat=\cat^{\mathcal D_{\mathcal Q}, \widehat {\underline w}_0,[a,b]}_{\g}$ and the case $\cat=\cat_w$ when $\g$ is symmetric.
(\cite{KKOP22, KKKO18}). 
In the case $\cat=\cat_w$ when $\bg$  is non-symmetric, this has not yet been established.

\medskip 
Lastly, we briefly mention a result \cite{Contu24} by Contu, which was informed to the authors while writing up this paper.
In \cite{Contu24} it is shown  that, for the case $\cat=\cat^{\mathcal D_{\mathcal Q}, \widehat {\underline w}_0,[a,b]}_{\g}$, the exchange matrix $B(\frakC)$ of the cluster associated with an admissible chain 
$\frakC$ of \iboxes is given by 
$B(\frakC) = P(\frakC)^{-1} B^{[a,b]}(\underline w_0) (P(\frakC)^t)^{-1}$, where $P(\frakC)$ and $B^{[a,b]}(\underline w_0)$ are specific matrices provided explicitly.
Hence  the matrix $\tB(\F)$ should coincide with 
(a truncation of) $B(\frakC)$
where $\frakC$ is an admissible chain of \iboxes such that $\F=\{\ci_k\}_{1\le k\le b-a+1}$.  
It is interesting to note that the main ingredient  in \cite{Contu24} is so called  the \emph{additive categorification of cluster algebras}, which is quite different from the methods used in this paper.

\bigskip
This paper is organized as follows. 
Section 2 is devoted to the combinatorics of \iboxes.  
In Section 3, we provide the definition of the matrix $\tB(\F)$.
In Section 4, we analyze vertical arrows explicitly
according to the configuration of horizontal arrows.
In Section 5, we recall  the monoidal categories $\cat_w$ and $\cat^{\mathcal D_{\mathcal Q}, \widehat {\underline w}_0,[a,b]}_{\g}$, and  prove Theorem \ref{thm:main}.

\mnoi
{\bf Acknowledgments.}\quad 
We thank Se-jin Oh and Euiyong Park for many fruitful discussions and Se-jin Oh for informing us of the paper \cite{Contu24}. The results of this paper were obtained during the second author's visit to the Research Institute for Mathematical Sciences (RIMS), Kyoto University from November 2022 to July 2023. The second author gratefully acknowledges the hospitality of the people at RIMS during his visit.
Finally, the authors would like to thank the anonymous referee for carefully reading the paper and providing many valuable comments that helped improve it.

\begin{convention}  Throughout this paper, we use the following convention.
\ben
\item For a statement $\ttP$, we set $\delta(\ttP)$ to be $1$ or $0$ depending on whether $\ttP$ is true or not. In particular, we set $\delta_{i,j}=\delta(i=j)$.
\item For an object $X$ of finite length in an abelian category, we denote by
$\hd(X)$ the head of $X$, the largest semisimple quotient of $X$, and by $\soc(X)$ the socle of $X$, the largest semisimple subobject of $X$. 
\item For $a,b \in\Z$ with $a\le b$, we set 
\begin{align*}
  & [a,b] =\{  k \in \Z \ | \ a \le k \le b\},
\end{align*}
and call it an \emph{interval}. 
\ee
\end{convention}

\section{ Combinatorics of \iboxes}
\subsection{Admissible chain of \iboxes}

Let $I$ be an index set and let $\bold i=(\im_k)_{k\in S}$ be a sequence in $I$ over  $S$, where $S$ is an interval in $\Z$.

Define
\eq \label{eq: nota +,-}
\ba{l}
s_+\seteq\min\bl\{t\in S  \mid s<t ,\; \im_t=\im_s \} \cup \{\infty\}\br >s, \\
s_-\seteq\max\bl\{t\in S
\mid t<s, \; \im_t=\im_s \} \cup \{-\infty\}\br<
s,\\
s(\jm)^+\seteq\min\bl\{t\in S \mid s \le  t,\; \im_t=\jm \} \cup \{\infty\}\br
\ge s, \\
s(\jm)^-\seteq\max\bl\{t\in S \mid t \le s,\; \im_t=\jm \} \cup \{-\infty\}\br\le s
\quad\qt{for $s\in S$ and $\jm\in I$.}\ea
\eneq

\bigskip

An interval $\ci=[x,y]$ in $S$ is called an \emph{$\bold i$-box} if $\im_x=\im_y$.
We refer  to  $\im_x$ as the \emph{color of  $[x,y]$},  denoted by $\im_\ci$.

For an $\bold i$-box $[x,y]$,  we set
\eqn
[x,y]_\phi \seteq \set{s\in [x,y]}{i_s=i_x}.
\eneqn

For an interval $[x,y]$ in $S$, we define the $\bold i$-boxes
\begin{align}\label{def:int}
[x,y\}\seteq[x,y(\im_x)^-] \quad \text{ and } \quad \{x,y]\seteq[x(\im_y)^+,y].
\end{align}

{\em In the sequel we only consider intervals and \iboxes
  in an interval $S$ \emph{of finite length} unless otherwise mentioned.}

\begin{definition}\label{def:com}
We say that the $\bold i$-boxes $[x_1,y_1]$ and $[x_2,y_2]$ \emph{commute} if
we have either
$$(x_1)_- < x_2 \le y_2 <  (y_1)_+ \qtq[or] (x_2)_- < x_1 \le y_1 <  (y_2)_+.$$
\end{definition}

\Lemma\label{lem:same}
Let $\ci$ and $\ci'$ be \iboxes with the same color.
If they commute, then we have either
$\ci\subset\ci'$ or $\ci'\subset\ci$.
\enlemma
\Proof
Since $\ci=[x,y]$ and $\ci'=[x',y']$ have the same color,
$x_-<x'\le y'<y_+$ (resp.\ $(x')_-<x\le y<(y')_+$) is equivalent to
$x\le x'\le y'\le y$ (resp.\ $x'\le x\le y\le y'$).
\QED

The following lemma follows immediately from the definition.
\Lemma
Assume  that the \iboxes $[x,y]$ and $[x',y']$ commute.
\bnum
\item If $x\le x'_-$, then $y'<y_+$,
\item If $y_+\le y'$, then $x'_-<x$.
\ee
\enlemma

\begin{definition}  \label{def: ad ch of i box}\hfill
\ben
\item A sequence $\frakC$ of \iboxes
$$ \frakC = (\ci_k = [x_k,y_k] )_{1 \le k \le l} \qtext{for $l \in \ \Z_{\ge 1}$ }$$
is called \emph{an admissible chain of \iboxes} if
$$\text{$\tc_k = [\tx_k,\ty_k] \seteq \bigcup_{1 \le j \le k} [x_j,y_j]$ is an interval with $|\tc_k| = k$ for $k\in [1,l]$}$$
and
$$ \text{either   $[x_k,y_k] = [\tx_k,\ty_k\}$ or $[x_k,y_k] = \{\tx_k,\ty_k]$  for  each  $k\in [1,l]$.}
$$

\item The interval $\tc_k$  is called
the {\em envelope} of $\ci_k$, and $\tc_l$ is called the {\em extent {\rm(}or range{\rm)}}
of $\frakC$. We understand $\tc_0=\emptyset$.
\end{enumerate}
\end{definition}
The set of envelopes of $\frakC$ is totally ordered by inclusion. 

Note that to give an admissible chain of \iboxes in $[a,b]$ is equivalent to give
an increasing sequence $(\tc_k)_{1\le k\le b-a+1}$ of intervals in $[a,b]$ such that $|\tc_k|=k$.
Indeed, the $\bold i$-boxes in an admissible chain $\frakC$ are uniquely determined by its envelopes:
\begin{align}\label{eq: T_k}
\ci_k=[x_k,y_k] =T_{k-1}[\tx_{k},\ty_{k}] \seteq
\bc  \ [\tx_{k},\ty_{k}\} &\text{ {\rm (i)} if $\tx_k=\tx_{k-1}-1$,}\\
\ \{\tx_{k},\ty_{k}]&\text{ {\rm (ii)} if $\ty_k=\ty_{k-1}+1$} \ec
\end{align}
for $1<k\le l$. In case  {\rm (i)} in~\eqref{eq: T_k}, we write
$T_{k-1}=\LL$, while $T_{k-1}=\RR$ in case (ii).

Thus, to an admissible chain of $\bold i$-boxes of length $l$, we can
associate a pair $(x, \mathfrak{T})$ consisting of an integer $x$
and a sequence $\mathfrak{T} = ( T_1,T_2,\ldots, T_{l-1})$ such that
$T_i \in \{ \LL,\RR\}$ $(1 \le i < l)$,
$$x_1=y_1=\tx_1=\ty_1 ,  \quad  [\tx_k,\ty_k] =
\begin{cases}
[\tx_{k-1}-1,\ty_{k-1}] & \text{ if } T_{k-1}=\LL, \\
[\tx_{k-1},\ty_{k-1}+1]  & \text{ if } T_{k-1}=\RR.
\end{cases}$$
Note that this association is bijective.

\Prop [{\cite[Lemma 5.2]{KKOP22}}]\label{prop:adch}
Let $\frakC= \{\ci_k\}_{1\le k\le l}$ be an admissible chain of $\bold i$-boxes.
Then we have
\eqn
&&(x_k)_-<\tx_k\le x_k\le y_k\le\ty_k<(y_k)_+,\\
&&\text{$(x_k)_-<x_j\le y_j<(y_k)_+ $ for $1\le j\le k\le l$.}
\eneqn
In particular, 
any pair of \iboxes in $\frakC$ commute with each other. 
\enprop

Let $\frakC= \{\ci_k\}_{1\le k\le l}$ be an admissible chain of $\bold i$-boxes with the associated pair $(x, \mathfrak{T})$. 
For $1\le k <  \ell$ we say that $\ci_k$ is \emph{movable} if $k=1$ or $T_{k-1}\neq T_k$ $(k\ge 2)$.

For a movable $\ci_k$ in $\frakC$, the {\em box move} of $\frakC$ at
$k$ is the  admissible chain
$B_k(\frakC)$ whose associated pair $(x',\mathfrak{T}')$ is given by
\begin{align*}
& {\rm (i)} \ \begin{cases}
x' = x \pm  1  &\text{if $k=1$ and $T_1 = \RR$ (resp.\ $\LL$),}  \\
x'=x & \text{if $k>1$,}
\end{cases} \allowdisplaybreaks \\
& {\rm (ii)} \ T'_{s} = T_s \ \  \text{ for } s \not\in \{ k-1,k\}
\quad \text{ and } \quad  \ T'_{s}\ne T_s\ \ \text{ for } s \in \{
k-1,k\}.
\end{align*}
That is, $B_k(\frakC)$ is the admissible chain obtained from
$\frakC$ by shifting $\tc_k$ by $1$ to the right or to the left inside
$\tc_{k+1}$. 

\begin{proposition} [{\cite[Proposition 5.6, Proposition 5.7]{KKOP22}}]
  \label{prop:boxmove}
Let $\frakC=\seq{\ci_k}_{1 \le k \le l}$ be an admissible chain of
$\bold i$-boxes and let $k_0$ be a movable $\bold i$-box $(1\le k_0<l)$. Set
$B_{k_0}(\frakC) = \seq{\ci_k'}_{1 \le k \le l}$. 
\bnum
\item
Assume that
$\im_{\tilde{x}_{k_0+1}} \ne \im_{\tilde{y}_{k_0+1}}$, i.e., $\tc_{k_0+1}$ is
not an $\bold i$-box. Then we have
$$\ci_k' = \ci_{\mathfrak{s}_{k_0}(k)} \qtext{for}\quad 1\le k \le l, $$
where $\mathfrak{s}_{k_0} \in \sym_l$ is the transposition of $k_0$
and $k_0+1$.
\item 
 Assume that $\im_{\tilde{a}_{k_0+1}} = \im_{\tilde{b}_{k_0+1}}$, i.e.,
$\tc_{k_0+1}$ is an $\bold i$-box. Set $\ci_{k_0+1}=[x,y]$ with
$x=\tx_{k_0+1}$ and $y=\ty_{k_0+1}$.
Then we have \bna
\item $\ci_{k_0}=[x_+,y]$ and $\ci'_{k_0}=[x,y_-]$ if $T_{k_0-1}=\RR$,\label{it:1}
\item $\ci_{k_0}=[x,y_-]$ and $\ci'_{k_0}=[x_+,y]$ if $T_{k_0-1}=\LL$.\label{it:2}
\end{enumerate}
\ee
\end{proposition}

Note that any two admissible  chains 
with the same range are related by a sequence of box moves.

\subsection{Maximal commuting family of $\bold i$-boxes}

In \cite[Proposition 5.3]{KKOP22},
it is proved that for any admissible chain $\frakC= \{\ci_k\}_{1\le k\le l}$  of $\bold i$-boxes with extent $[a,b]$,
$\F=\st{\ci_k\mid 1\le k\le l}$ is a {\em maximal commuting family} of $\bold i$-boxes in $[a,b]$, i.e., maximal among the
commuting families of $\bold i$-boxes in $[a,b]$.

\begin{proposition}[{\cite[Proposition 5.3]{KKOP22}}]\label{prop:Prop5.3}
Let $\frakC=\seq{\ci_k}_{1\le k \le l}$ be an admissible chain of
\iboxes and let $\ci$ be an \ibox such that $\ci\subset
\tc_l$ and $\ci$ commutes with all $\ci_k$ $(1 \le k \le l)$. Then
there exists $s\in[1,l]$ such that $\ci = \ci_s$.
\end{proposition}

We shall prove its converse.

\begin{lemma} \label{lem:fourboxes}
  Let $\F$ be a commuting family of $\bold i$-boxes.
  Let $s\le t$.
  Then we have
  either $|\st{y\mid [s,y]\in\F, y\le t}|\le 1$
  or $|\st{x\mid [x,t]\in\F, s\le x }|\le 1$.
\enlemma
\begin{proof}
  Assuming that  $|\st{y\mid [s,y]\in\F, y\le t}|\ge 2$
  and $|\st{x\mid [x,t]\in\F, s\le x}|\ge 2$, let us derive a contradiction.
  Assume that $[s,y_1],[s,y_2]\in\F$ with $s\le y_1<y_2\le t$ and
  $[x_1,t],[x_2,t]\in\F$ with $s\le x_1<x_2\le t$.
  
   Since $[s,y_1]$ and $[x_2,t]$ 
  commute, we have either
  $t<(y_1)_+$ or $(x_2)_-<s$.
  On the other hand, we have
  $(y_1)_+\le y_2\le t$ and $s\le x_1\le (x_2)_-$, which is a contradiction.
  \QED

\begin{lemma}\label{lem:fc}
Let $\F$ be a commuting family of $\bold i$-boxes in $[a,b]$. Then there exists an admissible chain $\frakC$ of $\bold i$ -boxes with extent $[a,b]$ such that any member of $\F$ appears in $\frakC$.
\end{lemma}
\begin{proof}
  We argue by induction on $b-a= l-1 $.
  If $a=b$, then the assertion is trivial.
  Assume that $a<b$.
  Since $[a,b\}$ and $\{a,b]$ commute with all the \iboxes in $[a,b]$,
  $\F'\seteq \F\cup\st{[a,b\},\{a,b]}$ is a commuting family of \iboxes in $[a,b]$.
  Hence, replacing $\F$ with $\F'$, we may assume from the beginning that
  $[a,b\},\;\{a,b]\in\F$.

  \smallskip
By Lemma~\ref{lem:fourboxes}, we have either
  \bna
\item
  $[a,b\}$ is a unique \ibox in $\F$ with $a$ as its end, or
  \item
  $\{a,b]$ is a unique \ibox in $\F$ with $b$ as its end.
  \ee
  In case (a), $\F''=\F\setminus\st{[a,b\}}$
  is a commuting family of \iboxes in $[a+1,b]$.
  Hence by the induction hypothesis,
  there exists an admissible chain $\frakC' = \{\ci_k\}_{1\le k\le l-1}$
  of \iboxes with extent $[a+1,b]$,
  which contains $\F'$.
  Then $\frakC= \{\ci_k\}_{1\le k\le l}$ with $\ci_l=[a,b\}$ satisfies the desired condition.

  The case (b) can be treated similarly.
\end{proof}

Lemma~\ref{lem:fc} says that the converse of Proposition~\ref{prop:Prop5.3}
is true.

\Cor  \label{cor:chain for max family}
Let $\F$ be a maximal commuting family of $\bold i$-boxes in $[a,b]$. Then there exists an admissible chain $\frakC= \{\ci_k\}_{1\le k\le l}$  of $\bold i$-boxes with extent $[a,b]$ such that $\F=\set{\ci_k}{1\le k \le l }$.
\encor
Note that the chain $\frakC$ associated to $\F$ in the above corollary is not unique in general,
as seen by Proposition~\ref{prop:boxmove} (i).
The following corollary easily follows from the existence of
an admissible chain.
\Cor Let $\F$ be a \mcf. Then $|\F|=b-a+1$.
\encor

\Prop \label{prop:efe}
Let $\F$ be a maximal commuting family of $\bold i$-boxes in an interval
$[a,b]$. Then for any $\ci=[x,y]\in \F$, there exists a unique $z\in\{x,y\}$ such that for any admissible chain $\frakC= \{\ci_k\}_{1\le k\le l}$  of $\bold i$-boxes  with extent $[a,b]$ consisting of $\bold i$-boxes in $\F$, if $\ci=\ci_k$ then $\st{z}=\tc_k\setminus\tilde{\ci}_{k-1}$.
\enprop
\Proof
Let $\frakC$ be an admissible chain $\frakC= \{\ci_k\}_{1\le k\le l}$  of $\bold i$-boxes  with extent $[a,b]$ consisting of $\bold i$-boxes in $\F$, and let
us take $k$ such that 
$\ci_k=[x,y]$ and let $\tc_k=[\tx,\ty]$ be its envelope.

If $x=y$, then the assertion is obvious. Hence we may  assume that $x<y$. Then we have $k>1$.

We divide into two cases:
\bna
\item
  there exists $[x,y'] \in \F$ such that $y'<y$,
  \item
  there is no $[x,y']\in\F$ such that $y'<y$.
\ee

\mnoi
(i) Assume (a).
We shall show that $\st{y}=\tc_k\setminus\tc_{k-1}$.

Take $1\le j\le l$ such that $\ci_j=[x,y']$.
Let $\tilde \ci_j = [\tilde x',\tilde y']$ be its envelope.
Since $\tilde y\ge y\ge (y')_+>\tilde y'$, we have $j <k$.
Hence, we have $x \in \ci_j \subset \tilde \ci_{k-1}$, which implies that
$y \notin \tilde \ci_{k-1}$.

\mnoi
(ii) Assume (b).
We shall show that $\st{x}=\tc_k\setminus\tc_{k-1}$.

Assuming that $x<y$ and $\st{y}=\tc_k\setminus\tc_{k-1}$, let us derive a
contradiction. Set $\tc_k=[\tilde x,y]$ with $\tilde x\le x$.
Then $\tc_{k-1}=[\tilde x,y-1]$.
Hence $x\in\tc_{k-1}$. Take the smallest $j\ge1$ such that
$x\in\tc_j$. Then  $1\le j <k$ and
$\tilde \ci_j = \{x\} \sqcup \tilde \ci_{j-1}$.
Hence $\ci_j=[x,y']$ for some $y'$.
Since $\ci_j\subset\tc_{k-1}$, we have $y'<y$, which contradicts (b). 
\QED 

\Def
We call $z$ in Proposition \ref{prop:efe}  the  \emph{\efe of $[x,y]$}.
\end{definition}

\Cor \label{cor:efeiff}
Let $\F$ be a maximal commuting family of $\bold i$-boxes in $[a,b]$ and let $[x,y]\in \F$.
Then
\begin{enumerate}[\rm (i)]
\item if there exists $[x,y'] \in \F$ such that $y'<y$, then $y$ is the \efe,
\item if there exists no  $[x,y'] \in \F$ such that $y'<y$, then $x$ is the \efe,
\item if there exists $[x',y] \in \F$ such that $x<x'$, then $x$ is the \efe,
\item if there exists no $[x',y] \in \F$ such that $x<x'$, then $y$ is the \efe.
  \setcounter{myc}{\value{enumi}}
\end{enumerate}
Summing up, we have

\bnum
  \setcounter{enumi}{\value{myc}}
\item 
the following three conditions are equivalent:
\bna
   \item $y$ is the \efe of $[x,y]$,
  \item
    $x=y$ or there exists $[x,y'] \in \F$ such that $y'<y$,
  \item
    there exists no $[x',y] \in \F$ such that $x<x'$.
 
    \ee
  \item   the following three conditions are equivalent:
    \bna
    \item $x$ is the \efe of $[x,y]$,    
  \item
$x=y$ or there exists $[x',y] \in \F$ such that $x<x'$,
  \item
    there exists no $[x,y'] \in \F$ such that $y'<y$.
\ee 
\end{enumerate}
\end{corollary}
\Proof
(i) and (ii) are shown in the proof of Proposition \ref{prop:efe}, and (iii) and (iv) follow by symmetry. 
The other assertions are immediate.
\QED

By the definition, the following lemma holds.
\Lemma \label{lem:efebijab}
Let $\F$ be a maximal commuting family of $\bold i$-boxes in $[a,b]$. 
The map from $\F$ to $[a,b]$ given by
\eqn
\F \  \ni\  \ci \longmapsto\; \text{the \efe  of $\ci$}\;\in [a,b]
\eneqn
is a bijection from $\F$ to $[a,b]$.
\enlemma
\Proof
Let $\frakC=\{\ci_k\}_{1\le k\le l}$ be an admissible chain of $\bold i$-boxes with extent $[a,b]$ such that $\F = \set{\ci_k}{1\le k\le l}$.
Then the above map sends $\ci_k$
to $\tc_k\setminus\tc_{k-1}$.
It is evidently bijective.
\QED

\Lemma \label{lem:middle}
Let $\F$ be a maximal commuting family of $\bold i$-boxes in $[a,b]$. 
\begin{enumerate}[\rm (i)]
\item  Let $[x,y]$, $[x',y]$ be $\bold i$-boxes in $\F$ such that $[x,y]\subset [x',y]$. 
Then for any $\bold i$-box $[x'',y]$ such that $[x,y]\subset [x'',y] \subset [x',y]$, we have $[x'',y]\in \F$.
\item  Let $[x,y]$, $[x,y']$ be $\bold i$-boxes in $\F$ such that $[x,y]\subset [x,y']$. 
Then for any $\bold i$-box $[x,y'']$ such that $[x,y]\subset [x,y''] \subset [x,y']$, we have $[x,y'']\in \F$.
\end{enumerate}
\enlemma
\Proof
Since the proof of (ii) is similar, we prove only (i). 
We may assume that $x<x''<x'$. 
By Lemma~\ref{lem:efebijab}, there exists an \ibox $\ci\in\F$ such that $x''$ is an \efe of $\ci$.
Then Lemma~\ref{lem:same} implies that
$[x,y]\subset\ci\subset[x',y]$.
Hence we obtain $\ci=[x'',y]$.
\QED

The above lemma, along with Corollary~\ref{cor:efeiff}, implies the following result.
\Cor\label{cor:ibpm}
Let $\F$ be a maximal commuting family of $\bold i$-boxes in $[a,b]$, and let $[x,y]\in \F$.
Then we have
\begin{enumerate}[\rm (i)]
\item $x$ is the \efe of $[x,y]$ if and only if  $x=y$ or $[x_+,y] \in \F$,
\item $y$ is the \efe of $[x,y]$ if and only if  $x=y$ or $[x,y_-] \in \F$.
\end{enumerate}
\end{corollary}

\Def
Let $\F$ be a maximal commuting family of $\bold i$-boxes in $[a,b]$. We set
\eq
&&\F_{\fr} = \set{[x,y]\in \F}{x_- < a \ \text{and} \ b < y_+}, \\
&&\F_{\ex} = \F \setminus  \F_\fr. \nonumber
\eneq
\end{definition}
Note that 
\eq \F_\fr = \set{[a(j)^+,b(j)^-]}{j \in \{i_a,\ldots,i_b\}},\eneq
since $[a(j)^+,b(j)^-]$ commutes with all the \iboxes in $[a,b]$.

\Lemma \label{lem:equicolored}
Let $\F$ be a maximal commuting family of \iboxes in $[a,b]$, and let $[x,y] \in\F_\ex$.
Then either $[x_-,y] \in \F$ or $[x,y_+]\in \F$.
 Note that $[x_-,y]$ and $[x,y_+]$ cannot both belong to the commuting family $\F$, since they do not commute.
\enlemma
\Proof
 Since $[x,y] \in\F_\ex$, we have either  $a\le x_-$ or $y_+\le b$.

(i) Assume that $a\le x_-$.
By Lemma~\ref{lem:efebijab}, there exists $\ci\in\F$ such that
$x_-$ is the \efe of $\ci$.
By Lemma~\ref{lem:same}, we have
$[x,y]\subset \ci$, and hence $\ci$ has the form $[x_-,z]$ for some $z\ge y$.
If $z=y$, then we have done.
Hence we may assume that $z\ge y_+$.
Then Corollary~\ref{cor:ibpm} implies that $[x,z]\in\F$.
Finally Lemma~\ref{lem:middle} implies that $[x,y_+]\in\F$.

\snoi
(ii) The case $y_+\le b$ can be treated similarly.
\QED

\Lemma\label{lem:oppefe}
Let $\F$ be a maximal commuting family of $\bold i$-boxes in $[a,b]$.
Assume that $x$ is the \efe of $[x,y]\in\F$
and $y'$ is the \efe of $[x',y']\in\F$.
\bnum
\item If $x'<x$, then $y<y'$,\label{item:yy'}
  \item if $y'<y$, then $x<x'$,\label{item:xx'}
\ee
\enlemma
\Proof
$\frakC=\{\ci_k\}_{1\le k\le l}$ be an admissible chain of $\bold i$-boxes with extent $[a,b]$ such that $\F = \set{\ci_k}{1\le k\le l}$.
Let us take $j$ and $k$ such that
$\ci_j=[x,y]$ and $\ci_k=[x',y']$.
Then their envelopes are given by
$\tc_j=[x,\ty]$ and $\tc_k=[\tx',y']$.

\snoi
(i)\ 
Since $\tx'\le x'<x$, we have $\tc_j\subset\tc_k$ and hence
$y\le \ty\le y'$. 
If $y=y'$, then $[x,y']\in\F$ and hence
Corollary~\ref{cor:efeiff} implies that $x'$ is the \efe of $[x',y']$, which is a contradiction.
Hence $y<y'$.

\snoi
(ii) is proved similarly.
\QED

\begin{lemma}\label{lem:ydyp}
Let $\F$ be a maximal commuting family of $\bold i$-boxes in $[a,b]$. 
\bnum
\item If $x$ and  $x'$ are the \efes of
  $[x,y], [x',y']\in \F$ respectively and $x\le x'$,  then $y'< y_+$.
\item If $y$ and  $y'$ are the \efes of
  $[x,y], [x',y']\in \F$ respectively  and $y'\le y$,  then $x_-< x'$.
\ee
\end{lemma}
\begin{proof}
  Since the proof is similar we prove  only (i).
  If $x=x'$, then $[x,y]=[x',y']$ by Lemma~\ref{lem:efebijab},
  and hence we have $y'<y_+$.
  Hence we may assume from the beginning that $x<x'$.
  
Let $\frakC=\{\ci_k\}_{1\le k\le l}$ be an admissible chain of $\bold i$-boxes with extent $[a,b]$ such that $\F = \set{\ci_k}{1\le k\le l}$. 
Let $\tc_k$ and $\tc_j$ be the envelopes of $[x,y]$ and $[x',y']$,  respectively. Then they are written as $\tc_j=[x',\ty']$ and $\tc_k=[x,\ty]$
for some $\ty'\ge x'$ and $\ty\ge x$.
Since $x<x'$, we have $j\le k$ and
$\tc_j \subset \tc_k$.
Hence we have
\eqn
y' \le \ty' \le \ty<y_+.
\eneqn
\end{proof}

\begin{lemma} \label{lem:xj+zefe}

  Let $\F$ be a maximal commuting family of $\bold i$-boxes in $[a,b]$, and
$[x,y]\in\F$  such that $x\neq y$. 
\bnum
\item
  If $x$ is the \efe of $[x,y]$ and 
  $x'$ satisfies $a\le x'\le y$ and $x'_-  \le x $, then there exists $y'\ge x'$ such that
  $[x', y']\in \F$  with \efe $x'$.
\item
  If $y$ be the \efe of $[x,y]$
  and $y'$ satisfies $x\le y'\le b$ and $ y\le   y'_+ $, then
there exists $x'\le y'$ such that $[x',y']\in \F$ with  \efe $y'$.
\ee
\end{lemma} 
\Proof
Because the proof is similar, we will prove (i) only.

By Lemma~\ref{lem:efebijab}, there exists a unique \ibox $\ci\in\F$
with $x'$ as its \efe.

If $\ci=[x', y']$ for some $y'$, then we have
done.

Assume that there is no such $y'$.
Then there exists $u$ such that $\ci=[u,x']$
with $u< x'$.  Then we have $u\le x'_-\le x$.

Corollary~\ref{cor:ibpm}, along with $u<x'$,
implies that $[u,x'_-]\in\F$.
The same corollary implies that $[x_+,y]\in\F$.
Then the commutativity of $[u,x'_-]$ and $[x_+,y]$
implies either $y<(x'_-)_+$ or $(x_+)_-<u$, which
contradicts $x'\le y$ and $u\le x$.
\QED

\begin{lemma} \label{lem:xj+yj-}
Let $\F$ be a maximal commuting family of $\bold i$-boxes in $[a,b]$. 
Assume that
$j\in I$
and  $[x,y] \in \F$ with  \efe $x$.
If $[x(j)^+,y']\in \F$, 
then,  for any $y''$ such that $i_{y''}=j$ and $y'<y''<y$, the \ibox $[x(j)^+, y'']$ belongs to $\F$,  with \efe $y''$

\end{lemma}
\Proof
Let $\ci$ be the \ibox with $y''$ as its \efe.
Then Lemma~\ref{lem:same} implies that $[x(j)^+,y']\subset\ci$ and hence
$\ci=[z,y'']$ for some $z\le x(j)^+$.

On the other hand, since
$x$ is the \efe of $[x,y]$, $y''$ is the \efe of $[z,y'']$
and $y''<y$,
Lemma~\ref{lem:oppefe}\;\eqref{item:xx'} implies that $x<z$,
which implies $z=x(j)^+$ since $z\le x(j)^+$.
\QED
\

\subsection{Structure of $\F_j$}
For  a maximal commuting family  $\F$ of $\bold i$-boxes in $[a,b]$ and  $j\in \{i_a,\ldots,i_b\}$, we set 
\eqn 
\F_j=\set{[x,y]\in \F}{i_x=j}. 
\eneqn
Then we have
\eqn\F=\bigsqcup_{j\in \{i_a,\ldots,i_b\}} \F_j.\eneqn

\Lemma \label{lem:Fj}
Let $\F$ be a maximal commuting family of $\bold i$-boxes in $[a,b]$, and let $j\in \{i_a,\ldots,i_b\}$.
Then there exists a unique
increasing sequence of \iboxes
$\{[x_k,y_k]\}_{1\le k\le m} $ in $\F_j$ such that
\bna
\item $\F_j=\st{[x_k,y_k]\mid 1\le k\le m}$,
\item $|[x_k,y_k]_\phi|=k$  for  $1\le k\le m$,  where
  $m=|\st{k\in[a,b]\mid i_k=j}|$,
  \item
    $[x_{k}, y_k]=[(x_{k+1})_+, y_{k+1}]$ or
    $[x_{k+1},(y_{k+1})_-]$ for $1 \le k\le m-1$,
\item $\F_\fr \cap \F_j =\st{[x_m,y_m]=\{[a(j)^+,b(j)^{-}]}$.
\ee
\enlemma
\Proof
Starting from $[x_m,y_m]=[a(j)^+,b(j)^{-}]$, we can define $[x_k,y_k]$
($1\le k\le m$)
inductively by Corollary~\ref{cor:ibpm}.
Since 
$\F_j$ is totally ordered by Lemma~\ref{lem:same},
$\F_j=\st{[x_k,y_k]\mid 1\le k\le m}$.
\QED

\Def
Let $\F$ be a maximal commuting family of $\bold i$-boxes in $[a,b]$
and $j\in I$.
 An $\bold i$-box $[x,y]$ in  $\F_j$ is said to be \emph{in the right corner} (of $\F_j$) if $[x,y_-]$ and $[x_-,y]$ belong to $\F_j$. 
 An $\bold i$-box $[x,y]$ in $\F_j$ is said to be \emph{in the left corner}  (of $\F_j$) if $[x_+,y]$ and $[x,y_+]$ belong to $\F_j$.
 \end{definition}
Note that if $[x,y]$ is in the right corner, then $y$ is the \efe of $[x,y]$ and $x_-$ is the \efe of $[x_-,y]$ by Corollary~\ref{cor:efeiff}. 
If $[x,y]$ is in the left corner, then $x$ is the \efe of $[x,y]$ and $y_+$ is the \efe of $[x,y_+]$.

\Lemma \label{lem:rightleft} 
Let $\F$ be a maximal commuting family of $\bold i$-boxes in $[a,b]$.
\begin{enumerate}[\rm (i)]
\item If  $[x,y] \in \F$ with \efe $y$,   then there exists a unique $y' \le y$ such that 
  $[x,y']\in \F$ 
  and either $[x,y']$ is in the left corner  or $x=y'$.

\item  If $[x,y]\in \F$ with \efe $x$, then there exists a unique $x'\ge x$ such that $[x', y]\in \F$ and
 either $[x',y]$ is in the right corner  or $x'=y$. 

\end{enumerate}

\enlemma
\begin{proof}
 Since the proofs are similar, we prove only (i). 

 \snoi
(i) We may assume that $x<y$.
Hence by Corollary \ref{cor:efeiff}, there exists $[x,y']\in \F$ such that $y'<y$. 
Take the smallest $y'$ among them.  
Since $y'<y'_+\le y$, we have $[x,y'_+] \in\F$ by Lemma \ref{lem:middle}.
Assume that $y'\neq x$. 
Then either $[x,y'_-] \in \F$ or $[x_+,y'] \in \F$ by Corollary~\ref{cor:ibpm},
and hence $[x_+,y'] \in \F$ by the choice of $y'$.
 It follows that $[x,y']$ is in the left corner.
The uniqueness follows from the fact that  $y''$ is the \efe of $[x,y'']$ for any $[x,y'']\in \F$ with $y'<y''$ by Corollary \ref{cor:efeiff}.
\end{proof}

We have the following proposition by Corollary \ref{cor:efeiff}.
\Prop 
Let $\F$ be a maximal commuting family of $\bold i$-boxes in $[a,b]$, and let $j\in \{i_a,\ldots,i_b\}$.
Let  $\{[x_k,y_k]\}_{1\le k\le m} $ be the enumeration of $\F_j$  as in \
Lemma~\ref{lem:Fj}.
Assume that $1\le p<q \le m$.
\begin{enumerate}[\rm (i)]
\item  If $p=1$ or $[x_p,y_p]$ is in the left corner, $[x_q,y_q]$ is in the right corner, and  $[x_k,y_k]$ is neither in the left corner nor in the right corner for  $p< k< q$,
then $y_k$ is the \efe of $[x_k,y_k]$ and $x_k=x_p$ for $p<k\le q$. 

\item If $p=1$ or  $[x_p,y_p]$ is in the right corner, $[x_q,y_q]$ is in the left corner, and  $[x_k,y_k]$ is neither in the left corner nor in the right corner  for $p< k< q$,
then $x_k$ is the \efe of $[x_k,y_k]$ and $y_k=y_p$ for $p<k\le q$.  

\item If $[x_q,y_q] \in \F_\fr$ or $[x_q,y_q]$ is in the left corner,  $[x_p,y_p]$ is in the right corner, and $[x_k,y_k]$ is neither in the left corner nor in the right corner for  $p< k< q$, then $x_k$ is the \efe of $[x_k,y_k]$ and $y_k=y_p$ for $p<k\le q$.

\item If $[x_q,y_q] \in \F_\fr$ or $[x_q,y_q]$ is in the right corner,  $[x_p,y_p]$ is in the left corner, and  $[x_k,y_k]$ is neither in the left corner nor in the right corner for  $p< k< q$, then $y_k$ is the \efe of $[x_k,y_k]$ and $x_k=x_p$ for $p<k\le q$. 
\end{enumerate}
\enprop

\section{Exchange matrices}

\subsection{Skew-symmetrizable exchange matrices} 
Let $\K =\Kex \sqcup \Kfr$ be a finite index set.
We call $\Kex$ the set of exchangeable indices, and $\Kfr$
the set of frozen indices. 
A matrix $\tB=(b_{s,t})_{s\in\K, t\in\Kex}$ is called a \emph{skew-symmetrizable exchange matrix} if $b_{s,t}\in\Z$, $|\set{s\in \K}{b_{s,t}\neq 0}| < \infty$ for all $t \in \Kex$, and the principal part  $B\seteq(b_{s,t})_{s\in\Kex, t\in \Kex}$  is skew-symmetrizable; i.e., there exists a tuple  $(d_s)_{s\in \Kex}\in \Z_{>0}^{\Kex}$ such that $d_s b_{s,t} = -d_t b_{t,s}$ for $s,t\in \Kex$.
We call $(d_s)_{s\in \Kex}$ a \emph{skew-symmetrizer}. 

We extend $\tB$ to $(b_{s,t})_{(s,t)\in (\K\times \K)\setminus(\Kfr\times\Kfr)}$ by
$d_s b_{s,t} = -d_t b_{t,s}$, 
if $(d_t)_{t\in \Kfr} \in\Z^{\Kfr}_{>0}$ is given.

For $k\in \Kex$, the \emph{mutation of $\tB$ in direction $k$} is the matrix $\mu_k(\tB)=(b'_{s,t})$ where
\eqn
b'_{s,t} = \begin{cases}
-b_{s,t} & \text{if $s=k$ or $t=k$,} \\
b_{s,t}+(-1)^{\delta(b_{s,k}<0)} [b_{s,k}b_{k,t}]_+ & \text{otherwise,}
\end{cases}
\eneqn
where $[a]_+\seteq\max(a,0)$.

If $\tB$ is a skew-symmetrizable exchange matrix, then so is $\mu_k(\tB)$ with the same skew-symmetrizer $(d_s)_{s\in \Kex}$ for any $k\in\Kex$.
\medskip

One may associate
 a quiver 
to a skew-symmetrizable exchange matrix $\tB$ (together with a skew-symmetrizer $\st{d_s}_{s\in\K}$)  by drawing an arrow
\eqn
 s\To[\ d_sb_{s,t}\ ] t  \qtext{ whenever   $s,t\in \K$ and $b_{s,t}>0$.}
\eneqn

\subsection{Exchange matrix for a maximal commuting family of $\bold i$-boxes}
\label{subsec:ex_matrix}
Let $\cartan=(\sfc_{i,j})_{i,j\in I}$ be a symmetrizable Cartan matrix with a symmetrizer $(d_i)_{i\in I}\in\Z_{>0}^I$: $d_i\sfc_{i,j}=d_j\sfc_{j,i}$. 
Let $[a,b]$ be an interval and $\bold i=(i_a,\ldots, i_b)$
be a sequence in $I$.
Let $\F$ be a \mcf 
associated with $\bold i$.

Let $\hB(\F)=(b_{\ci,\ci'})_{(\ci,\ci')\in\F\times\F}$
be the skew-symmetrizable exchange matrix with  index  set
$\F$ (without frozen indices)
together with the skew-symmetrizer $d_{[x,y]} = d_{i_x} $  for $[x,y] \in \F$
 whose positive entries 
are given as follows:
\eq \label{eq:ex_matrix} 
 &&\hs{3ex} b_{[x,y],[x',y']}=   \\
&&\hs{7ex}\begin{cases} 
1  & \text{if ($x=x'$ and $y'=y_-$) or ($y=y'$ and $x'=x_-$),}  \\
-\sfc_{i_x,i_{x'}}  & \text{if  $\sfc_{{i_x,i_{x'}}}<0$ and one of the following conditions $(a)$--$(d)$ is satisfied:} \nonumber\\
\end{cases}
\eneq
\begin{enumerate}[(a)]
\item  $[x,y_+] \in \F $,  \ $x$ is the \efe of $[x,y]$,   \ $x'_- <x<x'$,  \ $y'< y_+<y'_+$,
\item $[x,y_+] \in \F $,   \  $y'$ is the \efe of $[x',y']$,   \ $x'_- <x$,  \ $y<y'< y_+<y'_+$,
\item $[x'_-,y'] \in \F $, \   $y'$ is the \efe of $[x',y']$,   \ $x_-<x'_- <x$,   \ $y<y'< y_+$,
\item $[x'_-,y'] \in \F$,  \   $x$ is the \efe of $[x,y]$,   \ $x_-<x'_- <x<x'$,  \ $y'< y_+$.
\end{enumerate}

Such a matrix $\hB(\F)$ exists since
$ b_{[x,y],[x',y']}>0$ and $ b_{[x',y'],[x,y]}>0$
cannot happen simultaneously.
We set
$$\tB(\F)=(b_{\ci,\ci'})_{(\ci,\ci')\in\F\times\Fex}.$$

Let $Q(\F)$ be the quiver associated with $\tB(\F)$; that is, a quiver  with $\F$ as the set of vertices and with the set of the arrows given as follows:
\eqn
&&\ba{ll}
\qt{$\bullet$\ horizontal arrows} : \quad &[x,y] \To[\ d_{i_x}\ ] [x',y'] \\ 
 &\text{if $x=x'$ and $y'=y_-$ or $y=y'$ and $x'=x_-$,}  \\[2ex]
\qtext{$\bullet$\ vertical arrows} : \quad &[x,y] \To[\ -d_{i_x}\sfc_{i_x,i_{x'}}\ ] [x',y']   \\
&\text{if  $\sfc_{{i_x,i_{x'}}}<0$ and one of the above conditions $(a)$--$(d)$ is satisfied}.
\ea
\eneqn
We denote a \emph{horizontal arrow} simply by $[x,y] \To{} [x,'y']$
 when there is no afraid of confusion.

\Rem
Note that we have always
\eqn
b_{[x,y],[x',y']}\in\st{1,-1,-\sfc_{i_x,i_{x'}}, \sfc_{i_x,i_{x'}},0}
\qt{for any $[x,y],[x',y']\in \F$.}\eneqn
  \enrem

\subsection{Example} \label{subsec:example_exch_matrix}

Let $\cartan$ be the Cartan matrix of type $C_3$:

$$\cartan = \left(\begin{array}{rrr}
2 & -1 & 0 \\
-1 & 2 & -2 \\
0 & -1 & 2
\end{array}\right).$$

Consider the sequence $\bold i$ in $I=\{1,2,3\}$ 
$$\bold i =\left(1, 3, 2, 1, 3, 3, 3, 3, 1, 2, 1, 3, 3, 2, 3, 3, 2, 1, 3, 1\right).$$

Let  $\F$ be a maximal  commuting family of \iboxes
$$
\begin{aligned}
\F=& \{ \left[15, 15\right], \left[14, 14\right], \left[15, 16\right], \left[14, 17\right], \left[13, 16\right], \left[18, 18\right], \left[13, 19\right], \left[12, 19\right], \left[18, 20\right], \\
 & \left[11, 20\right],  \left[10, 17\right], \left[9, 20\right], \left[8, 19\right], \left[7, 19\right], \left[6, 19\right], \left[5, 19\right], \left[4, 20\right], \left[3, 17\right], \left[2, 19\right], \left[1, 20\right] \}.
 \end{aligned}
$$

Enumerating the boxes in $\F$ by $\ci_k$ for $1\le k\le 20$ in the order described above,  we obtain 
 an admissible chain $\frakC= \{\ci_k\}_{1\le k\le l}$    of $\bold i$-boxes with extent $[1,20]$,
 associated with the pair 
$$(15, ( \LL, \RR, \RR, \LL, \RR, \RR, \LL, \RR, \LL, \LL, \LL, \LL, \LL, \LL, \LL, \LL, \LL, \LL, \LL)).$$

The following is obtained by stacking the $\ci_k$'s  from bottom to top, with their colors on the right:
\begin{center}
\begin{tikzpicture}[scale=0.4]
  \draw (15, 0) rectangle (16, 1);
  \node at (15.5000000000000, 1/2) {15};
  \node at (15.5000000000000, 1/2) {15};
  \node[right] at (22, 1/2) {3};
  \draw (14, 1) rectangle (15, 2);
  \node at (14.5000000000000, 3/2) {14};
  \node at (14.5000000000000, 3/2) {14};
  \node[right] at (22, 3/2) {2};
  \draw (15, 2) rectangle (17, 3);
  \node at (15.5000000000000, 5/2) {15};
  \node at (16.5000000000000, 5/2) {16};
  \node[right] at (22, 5/2) {3};
  \draw (14, 3) rectangle (18, 4);
  \node at (14.5000000000000, 7/2) {14};
  \node at (17.5000000000000, 7/2) {17};
  \node[right] at (22, 7/2) {2};
  \draw (13, 4) rectangle (17, 5);
  \node at (13.5000000000000, 9/2) {13};
  \node at (16.5000000000000, 9/2) {16};
  \node[right] at (22, 9/2) {3};
  \draw (18, 5) rectangle (19, 6);
  \node at (18.5000000000000, 11/2) {18};
  \node at (18.5000000000000, 11/2) {18};
  \node[right] at (22, 11/2) {1};
  \draw (13, 6) rectangle (20, 7);
  \node at (13.5000000000000, 13/2) {13};
  \node at (19.5000000000000, 13/2) {19};
  \node[right] at (22, 13/2) {3};
  \draw (12, 7) rectangle (20, 8);
  \node at (12.5000000000000, 15/2) {12};
  \node at (19.5000000000000, 15/2) {19};
  \node[right] at (22, 15/2) {3};
  \draw (18, 8) rectangle (21, 9);
  \node at (18.5000000000000, 17/2) {18};
  \node at (20.5000000000000, 17/2) {20};
  \node[right] at (22, 17/2) {1};
  \draw (11, 9) rectangle (21, 10);
  \node at (11.5000000000000, 19/2) {11};
  \node at (20.5000000000000, 19/2) {20};
  \node[right] at (22, 19/2) {1};
  \draw (10, 10) rectangle (18, 11);
  \node at (10.5000000000000, 21/2) {10};
  \node at (17.5000000000000, 21/2) {17};
  \node[right] at (22, 21/2) {2};
  \draw (9, 11) rectangle (21, 12);
  \node at (9.50000000000000, 23/2) {9};
  \node at (20.5000000000000, 23/2) {20};
  \node[right] at (22, 23/2) {1};
  \draw (8, 12) rectangle (20, 13);
  \node at (8.50000000000000, 25/2) {8};
  \node at (19.5000000000000, 25/2) {19};
  \node[right] at (22, 25/2) {3};
  \draw (7, 13) rectangle (20, 14);
  \node at (7.50000000000000, 27/2) {7};
  \node at (19.5000000000000, 27/2) {19};
  \node[right] at (22, 27/2) {3};
  \draw (6, 14) rectangle (20, 15);
  \node at (6.50000000000000, 29/2) {6};
  \node at (19.5000000000000, 29/2) {19};
  \node[right] at (22, 29/2) {3};
  \draw (5, 15) rectangle (20, 16);
  \node at (5.50000000000000, 31/2) {5};
  \node at (19.5000000000000, 31/2) {19};
  \node[right] at (22, 31/2) {3};
  \draw (4, 16) rectangle (21, 17);
  \node at (4.50000000000000, 33/2) {4};
  \node at (20.5000000000000, 33/2) {20};
  \node[right] at (22, 33/2) {1};
  \draw (3, 17) rectangle (18, 18);
  \node at (3.50000000000000, 35/2) {3};
  \node at (17.5000000000000, 35/2) {17};
  \node[right] at (22, 35/2) {2};
  \draw (2, 18) rectangle (20, 19);
  \node at (2.50000000000000, 37/2) {2};
  \node at (19.5000000000000, 37/2) {19};
  \node[right] at (22, 37/2) {3};
  \draw (1, 19) rectangle (21, 20);
  \node at (1.50000000000000, 39/2) {1};
  \node at (20.5000000000000, 39/2) {20};
  \node[right] at (22, 39/2) {1};

\end{tikzpicture}
\end{center}

One can read the \efe of each $\ci_k$ from bottom to top as follows:
$$15, 14, 16,17,13,18,19,12,11,10,9,8,7,6,5,4,3,2,1$$

By \eqref{eq:ex_matrix},  the exchange matrix  $\hB(\F)$ is given as follows:

\begin{center}
\resizebox{.82\textwidth}{!}{
$
\left(\begin{array}{rrrrrrrrrrrrrrrrrrrr}
0 & 0 & -1 & 0 & 0 & 0 & 0 & 0 & 0 & 0 & 0 & 0 & 0 & 0 & 0 & 0 & 0 & 0 & 0 & 0 \\
0 & 0 & 2 & -1 & 0 & 0 & 0 & 0 & 0 & 0 & 0 & 0 & 0 & 0 & 0 & 0 & 0 & 0 & 0 & 0 \\
1 & -1 & 0 & 0 & 1 & 0 & 0 & 0 & 0 & 0 & 0 & 0 & 0 & 0 & 0 & 0 & 0 & 0 & 0 & 0 \\
0 & 1 & 0 & 0 & -2 & 0 & 2 & -2 & 1 & -1 & 1 & 0 & 0 & 0 & 0 & 0 & 0 & 0 & 0 & 0 \\
0 & 0 & -1 & 1 & 0 & 0 & -1 & 0 & 0 & 0 & 0 & 0 & 0 & 0 & 0 & 0 & 0 & 0 & 0 & 0 \\
0 & 0 & 0 & 0 & 0 & 0 & 0 & 0 & -1 & 0 & 0 & 0 & 0 & 0 & 0 & 0 & 0 & 0 & 0 & 0 \\
0 & 0 & 0 & -1 & 1 & 0 & 0 & 1 & 0 & 0 & 0 & 0 & 0 & 0 & 0 & 0 & 0 & 0 & 0 & 0 \\
0 & 0 & 0 & 1 & 0 & 0 & -1 & 0 & 0 & 0 & -1 & 0 & 1 & 0 & 0 & 0 & 0 & 0 & 0 & 0 \\
0 & 0 & 0 & -1 & 0 & 1 & 0 & 0 & 0 & 1 & 0 & 0 & 0 & 0 & 0 & 0 & 0 & 0 & 0 & 0 \\
0 & 0 & 0 & 1 & 0 & 0 & 0 & 0 & -1 & 0 & -1 & 1 & 0 & 0 & 0 & 0 & 0 & 0 & 0 & 0 \\
0 & 0 & 0 & -1 & 0 & 0 & 0 & 2 & 0 & 1 & 0 & 0 & 0 & 0 & 0 & -2 & -1 & 1 & 0 & 0 \\
0 & 0 & 0 & 0 & 0 & 0 & 0 & 0 & 0 & -1 & 0 & 0 & 0 & 0 & 0 & 0 & 1 & 0 & 0 & 0 \\
0 & 0 & 0 & 0 & 0 & 0 & 0 & -1 & 0 & 0 & 0 & 0 & 0 & 1 & 0 & 0 & 0 & 0 & 0 & 0 \\
0 & 0 & 0 & 0 & 0 & 0 & 0 & 0 & 0 & 0 & 0 & 0 & -1 & 0 & 1 & 0 & 0 & 0 & 0 & 0 \\
0 & 0 & 0 & 0 & 0 & 0 & 0 & 0 & 0 & 0 & 0 & 0 & 0 & -1 & 0 & 1 & 0 & 0 & 0 & 0 \\
0 & 0 & 0 & 0 & 0 & 0 & 0 & 0 & 0 & 0 & 1 & 0 & 0 & 0 & -1 & 0 & 0 & -1 & 1 & 0 \\
0 & 0 & 0 & 0 & 0 & 0 & 0 & 0 & 0 & 0 & 1 & -1 & 0 & 0 & 0 & 0 & 0 & -1 & 0 & 1 \\
0 & 0 & 0 & 0 & 0 & 0 & 0 & 0 & 0 & 0 & -1 & 0 & 0 & 0 & 0 & 2 & 1 & 0 & 0 & 0 \\
0 & 0 & 0 & 0 & 0 & 0 & 0 & 0 & 0 & 0 & 0 & 0 & 0 & 0 & 0 & -1 & 0 & 0 & 0 & 0 \\
0 & 0 & 0 & 0 & 0 & 0 & 0 & 0 & 0 & 0 & 0 & 0 & 0 & 0 & 0 & 0 & -1 & 0 & 0 & 0
\end{array}\right)$
}
\end{center}

\
The quiver $Q(\F)$ is given as follows (for simplicity, the arrows with  label $1$ are left unlabelled):
\begin{center}
\resizebox{\textwidth}{!}{
$
\xymatrix@C=0.6cm{
[18,18] & [18,20] \ar[r] \ar[l] & [11,20]\ar[r] \ar[dl]& [9,20] \ar[r] & [4,20]\ar[r] \ar[dll]  & [1,20]&&&& \\
[14,14]  \ar[dr]|2  & [14,17]\ar[r]\ar[l]   \ar[u] \ar[drr]|2 & [10,17]\ar[r] \ar[u]\ar[drr]|2  & [3,17] \ar[drrrrr]|2 \ar[ur]&  & &&&& \\
[15,15]  & [15,16 ] \ar[r]|2 \ar[l]|2  & [13,16] \ar[ul]|2 &[13,19] \ar[r]|2 \ar[l]|2 &[12,19] \ar[r]|2 \ar[ulll]|2 &[8,19]\ar[r]|2 &[7,19]\ar[r]|2 &[6,19]\ar[r]|2 &[5,19]\ar[r]|2 \ar[ullllll]|2 &[2,19] 
}
$
}
\end{center}

\section{Vertical arrows} \label{sec:vertical_arrows}
Let $\F$ be a \mcf.
In this section,  we fix an $\bold i$-box
$$\text{$\ci_0\seteq[x,y]\in \F_\ex$  with color $i$.  }$$
By Lemma~\ref{lem:equicolored}, we have either $[x,y_+]\in\F$ or $[x_-,y]\in\F$; hence there is a horizontal arrow
adjacent to $\ci_0$.
We shall analyze vertical arrows adjacent to $\ci_0$
according to the configuration of the horizontal arrows
adjacent to $\ci_0$.

\smallskip
Recall that $\F_j\seteq\st{\ci\in\F\mid i_\ci=j}$ for $j\in I$.
Define
\eqn
&&\Vo_j=  \set{\ci' \in \F_j }{ b_{\ci',\ci_0}<0},\\
&&\Vi_j= \set{\ci'\in \F_j }{ b_{\ci',\ci_0}>0},
\eneqn
and  set
\eqn &&\Vo=\bigsqcup_{j;\;\sfc_{i,j}<0}\Vo_j, \qquad
\Vi=\bigsqcup_{j;\;\sfc_{i,j}<0}\Vi_j.
\eneqn

\subsection{Case: $[x_+,y] \to [x,y ] \to  [x_-,y]$} \label{subsec:>>} 
\quad
\hfill

Suppose that $[x,y]\in \Fex$ with $[x_+,y], [x_-,y]\in \F$.  
 Set $i:=i_x$. 
In this case, $x$ is the \efe of $\ci_0=[x,y]$,
and $x_-$ is the \efe of $[x_-,y]$.
Moreover $[x, y_+]\not\in\F$.

\begin{lemma} \label{lem:>>Varrows}
Assume that $j$ satisfies  $\sfc_{i,j}<0$.
Then  we have
\eqn
\Vo_j=\Vo_j^e \sqcup \Vo_j^o ,  \qt{and} \  \Vi_j=\Vi_j^e \sqcup \Vi_j^o,
\eneqn
where 
\eqn 
&&\Vo_j^e =  \set{[x',y']\in \F_j}{ [x'_-,y']\in \F,  \ \text{ $y'$ is the \efe of $[x',y']$},   \ x_-< x'_- < x' <x }, \\
&&\Vo_j^o =  \set{[x', y']\in \F_j}{[x'_-,y']\in \F,  \ x_-< x'_- < x <x' },\\
&&\Vi_j^e = \set{[x',y']\in \F_j}{ [x',y'_+]\in \F,  \ \text{$x'$ is the \efe of $[x',y']$},  \  x_-< x'_- < x' <x }, \\
&&\Vi_j^o =  \set{[x',y']\in \F_j}{\text{$x'$ is the \efe of $[x',y']$ },    \ x'_-< x_- < x' <x } .
\eneqn
\end{lemma}
\begin{proof}
(1)\ By the definition,  we have 
\eqn
\Vo_j&&=\set{[x',y']\in\F_j}{   [x'_-,y'] \in \F,  \text{$y'$ is the \efe of $[x',y']$},    \ x_- < x'_- < x , \ y< y' < y_+ } \\
&& \quad \cup \set{[x',y']\in\F_j}{  [x'_-,y'] \in \F,   \text{$x$ is the effective end of $[x,y]$},  \  x_- <x'_-<x<  x', \  y'< y_+  }  \\
&&=\set{[x',y']\in\F_j}{[x'_-,y'] \in \F,  \text{$y'$ is the \efe of $[x',y']$},    \ x_- < x'_- < x , \ y< y' } \\
&& \quad \cup \set{[x',y']\in\F_j}{  [x'_-,y'] \in \F,    \  x_- <x'_-<x<  x' }  \\
&&=\set{[x',y']\in\F_j}{   [x'_-,y'] \in \F,  \text{$y'$ is the \efe of $[x',y']$},    \ x_- < x'_- <x'< x, \ y < y' }\\
&& \quad  \sqcup \set{[x',y']\in\F_j}{  [x'_-,y'] \in \F,    \  x_- <x'_-<x<  x' }\\
&&=\set{[x',y']\in\F_j}{  \text{$y'$ is the \efe of $[x',y']$},    \ x_- < x'_- <x'< x}\\
&& \quad  \sqcup \set{[x',y']\in\F_j}{ [x'_-,y'] \in \F,    \  x_- <x'_-<x<  x'
 }.
\eneqn
Indeed,  $x_-$ is the \efe of $[x_-,y]\in\F$.
and  $x'_-$ is the \efe of,
$[x'_-,y'] \in \F$ , the inequality $x_-< x'_-$ implies $y'< y_+$
by Lemma  \ref{lem:ydyp}. Hence the second equality follows.

The third equality follows by dividing the cases $x'<x$ and $x<x'$.

The fourth equality follows from
Lemma~\ref{lem:oppefe}\;\eqref{item:yy'}:
$x'<x$ implies $y<y'$
since $x$ is the \efe of $[x,y]$ and $y'$ is the \efe of $[x',y']$.
Hence we obtain $\Vo_j=\Vo_j^e\sqcup \Vo_j^o$.

 \mnoi
 (2)
By the definition we have
\eqn
\Vi_j&&=\set{[x',y']\in\F_j}{  [x',y'_+] \in \F,  \text{$x'$ is the \efe of $[x',y']$},    \ x_- < x' < x , \ y< y'_+ < y_+ } \\
&& \quad \cup \set{[x',y']\in\F_j}{[x_-,y] \in \F,   \text{$x'$ is the effective end of $[x',y']$},  \  x'_- <x_-<x'<  x, \  y< y'_+  }  \\
&&=\set{[x',y']\in\F_j}{ [x',y'_+] \in \F,  \text{$x'$ is the \efe of $[x',y']$},    \ x_- < x'_-< x' < x , \  y'_+ < y_+ } \\
&& \quad \sqcup \set{  [x',y']\in\F_j}{ \text{$x'$ is the effective end of $[x',y']$},  \  x'_- <x_-<x'<  x}.
\eneqn
Here the first equality follows from the fact that $y$ is not the \efe of $[x,y]$.
Since $x'$ is the \efe of $[x',y'] \in \F$,
 $x$ is the \efe of $[x,y] \in \F$  and $x'< x$,
Lemma~\ref{lem:ydyp} implies that $y< y'_+$, which implies the second equality.

\medskip
Assume that $[x',y'_+] \in \F$.
Since  $[x',y'_+]$ and $[x_-,y]$ commute, $x_-\le x'_-$ implies that
$y'_+ < y_+$.
Hence we obtain $\Vi_j=\Vi_j^e\sqcup \Vi_j^o$,  as desired.
\end{proof}

\begin{lemma} \label{lem:>>Vej} 
Let $\sfc_{i,j}<0$.  We have
\eqn
\Vo^e_j=&&\set{[x',y']\in \F_j}{ \text{$[x',y']$ is in the right corner and} \ x_-(j)^+ <x' \le x(j)^-},\\
\Vi^e_j=&&\set{[x',y']\in \F_j}{\text{$[x',y']$ is in the left corner and} \ x_-(j)^+ <x' \le x(j)^-} \\ 
&&\hs{3ex}\sqcup \set{[x',x']\in \F_j }{[x',x'_+]\in \F,  \ x_-(j)^+ <x' \le x(j)^-}.
\eneqn
\end{lemma}

\begin{proof}
First note that  $i\neq j$  and hence
\eq \label{eq:xddx}
\text{$x_-(j)^+ <x' \le x(j)^-$ is equivalent to $x_-<x'_-<x'<x$.}
\eneq

\snoi
(1) Assume that $[x',y']\in \Vo^e_j$. 
 Then $y<y'$ as shown in the proof of Lemma \ref{lem:>>Varrows}.
 Hence we have $x'<x< y< y'$ so that that $x'<y'$.  
 Because $y'$ is the \efe of $[x',y']$,  we obtain that $[x',y'_-]\in \F$ by
Corollary~\ref{cor:ibpm} and hence $[x', y']$ is in the right corner. 
Conversely if $[x', y']$  is  in the right corner of $\F_j$,  then 
$[x'_-,y']\in \F$  and $y'$ is the \efe of $[x',y']$. 
Hence the assertion for $\Vo^e_j$ follows from Lemma \ref{lem:>>Varrows}.

\snoi
(2) Note that $[x', y'_+] \in \F$ and $[x',y']$ has the \efe $x'$ if and only if 
either $[x', y']$ is in the left corner or $y'=x'$ and  $[x',x'_+]\in \F$.
Hence the assertion for $\Vi^e_j$ follows from Lemma \ref{lem:>>Varrows}.
\end{proof}

\begin{lemma} \label{lem:>>Voj}
We have
\eqn
&&\Vo^o_j=
\set{ [x(j)^+,y']\in \F }  {\text{$[x(j)^-,y'] \in \F $ and $x_-(j)^+ < x$}}, 
 \\
&&\Vi^o_j=
\set{ [x_-(j)^+,y'] \in \F} { \text{  $x_-(j)^+$ is the \efe of $[x_-(j)^+,y']$ and $x_-(j)^+ < x$ }}.
\eneqn
We have $|\Vo^o_j| \le 1$ and $|\Vi^o_j| \le 1$.
\end{lemma}
\begin{proof}
The first assertion follows  from Lemma \ref{lem:>>Varrows} together with that
\eq \label{eq:xdxd}
x_-< x'_- < x <x'  \  \text{is equivalent to}  \ x'=x(j)^+ \  \text{and} \ x_-(j)^+ < x,
\eneq
and
\eq \label{eq:dxdx}
x'_-< x_- < x' <x  \  \text{is equivalent to}  \ x'=x_-(j)^+ \  \text{and} \ x_-(j)^+< x.
\eneq

 If $[x(j)^+,y'], [x(j)^-,y'] \in \F $,  then $x(j)^-$ is the \efe of $[x(j)^-,y'] $ so that $y'$ is determined by Lemma \ref{lem:efebijab}.
Hence the second assertion follows.
\end{proof}

Note that  $\Vi_j=\Vo_j=\emptyset$ if $x_-(j)^+>x$.

\begin{proposition} \label{prop:>>Varrow_generic} 
Assume that  ${x_-(j)}^+ < x$ .
Then we have the following: 
\bnum
\item  There exist $w$ and $z$ such that
  $ x(j)^-\le w  \le z \le b$,  $[x_-(j)^+, z] \in \F$ with \efe $x_-(j)^+$, and 
 $[x(j)^-, w] \in \F$ with  \efe $x(j)^-$.

\item 
We have
\eqn
\Vi^o_j =\{[x_-(j)^+, z]\}. 
\eneqn

\item
We have
\eq \label{eq:Vooj}
\Vo^o_j =
\begin{cases}
\{ [ x(j)^+,  w]\}  & \text{if} \ x(j)^+\le w,\text{equivalently,}  \ x(j)^-<w. \\
 \emptyset & \text{if} \ x(j)^+> w ,\text{equivalently,}  \ x(j)^-=w.
\end{cases}
\eneq

\item If $w=z$, then $\Vo^e_j=\Vi^e_j=\emptyset$.

\item Assume that $w<z$. Then $\Vo^e_j \neq \emptyset$.
There is a bijection between $\Vo^e_j$ and $\Vi^e_j$ sending $[x',y']$ to $[x',y'']$, where $y''$ is the smallest element such that $[x',y'']\in \F$.

More precisely, 
if \eq
\Vo^e_j=\set{[x^{(k)},y^{(k)}] }{1\le k \le t}\eneq
with a strictly decreasing sequence $\st{[x^{(k)},y^{(k)}]}_{1\le k \le t}$
{\rm(}see Lemma~\ref{lem:same}{\rm)},
then we have  
\eq
\Vi^e_j=\set{[x^{(k)},y^{(k+1)}] }{1\le k \le t-1} \cup \{[x^{(t)},w]\}.
\eneq
 Moreover, $y^{(1)}=z$. 
\end{enumerate}
\end{proposition}
\begin{proof}
By the assumption, we have $x_-(j)^+ \le x(j)^-$. 

  \noi

  (i) \ Note that $x\neq y$ since $[x_+,y] \in \F$. Hence  there exist $w,z$ such that
  $[x_-(j)^+, z] \in \F$ with \efe $x_-(j)^+$, and 
  $[x(j)^-, w] \in \F$ with \efe $x(j)^-$
by Lemma \ref{lem:xj+zefe},
  since $[x,y] \in \F$ with \efe $x$ and $x_-(j)^+ \le x(j)^-\le x$.
We have $z \ge w$ by Lemma~\ref{lem:same}.

\mnoi
(ii) and (iii)
follow from (i) and  Lemma \ref{lem:>>Voj}. 

\mnoi
(iv)  If $x_-(j)^+ =x(j)^-$, then $\Vo^e_j=\Vi^e_j=\emptyset$ by Lemma \ref{lem:>>Vej}. 

Assume that $x_-(j)^+ <x(j)^-$. Then for any $\bold i$-box $[x',y']$ in $\F_j$  with $x_-(j)^+< x'\le  x(j)^-$ we have  $y'=z=w$.
Moreover $[x'_-,z] \in \F$ since $[x_-(j)^+,y']$, $[x',y']\in\F$,
and $x'$ is the \efe of $[x',y']$ since $x(j)^-$ is the \efe of $[x(j)^-,y']$.

Hence $[x',y']$ is not  in the right corner.
Moreover $[x',y_+] \not\in \F$
by Lemma~\ref{lem:same} and $[x'_-,y']\in\F$. Hence $\Vo^e_j=\Vi^e_j=\emptyset$ by Lemma \ref{lem:>>Vej}.

\medskip
(v) 
Note that $x_-(j)^+<x(j)^-\le w<z$.
Since $x_-(j)^+ $ is the \efe of $[x_-(j)^+,z]$, there exists $x^{(1)}$ such that $x_-(j)^+ < x^{(1)}$  and  $[x^{(1)},z] \in \F$.  
Take the largest $x^{(1)}$ among them. Since $w< z$ and $[x(j)^-,w]\in \F$,  $[x^{(1)},z]$ lies in the right corner.
 Note that $[x(j)^-,w] \subset [x^{(1)},z]$ since $w<z$. It follows that $x^{(1)}\le x(j)^-$  and hence $[x^{(1)},z] \in \Vo^e_j$ by Lemma \ref{lem:>>Vej}.

\smallskip
If $[x',y']\in \Vo^e_j$,  then  by Lemma \ref{lem:>>Vej}  it is in the right corner. 
Let $y''$ be the smallest element such that $[x',y'']\in \F$ and $y''<y$. By Lemma  \ref{lem:rightleft} such $y''$ exists.  
Moreover, either $[x',y'']$ in the left corner or $y''=x'$.  In the both cases,  we have $[x',y''] \in \Vi^e_j$.

Conversely,  assume that  $[x',y'']\in \Vi^e_j$.  
Let $y'$ be the the  largest element such that $[x',y']\in \F$ and $y'<y''$. 
By Lemma  \ref{lem:rightleft} such $y'$ exists. 
 Moreover either $[x',y']$ is in the right corner or $[x',y']\in\F_\ex$.
 Because $x_-(j)^+  < x'$,  $[x',y']$ does not contain the box $[x_-(j)^+, z]\in \F$ so that $[x',y']$ is not a frozen.  Thus $[x',y']$ is in the right corner.  It follows that $[x',y'] \in \Vo^e_j$ by Lemma 
\ref{lem:>>Vej},    as desired.

\smallskip

Let \eqn
\Vo^e_j=\set{[x^{(k)},y^{(k)}] }{1\le k \le t}\eneqn
such that 
\eqn 
[x^{(1)},y^{(1)}]  \supsetneq  [x^{(2)},y^{(2)}]  \supsetneq  \cdots  \supsetneq   [x^{(t)},y^{(t)}].  
\eneqn
Then we have  
$y^{(1)}=z$ since $[x^{(1)},z]$ is the largest $\bold i$-box in $\Vo^e_j$,
and
it is easy to see that $y^{(k+1)}$ is the smallest element such that   $[x^{(k)}, y^{(k+1)}]\in \F_j$ for $1\le k\le t-1$.

Let  $v$ be the smallest element such that $[x^{(t)},v]\in \F$. 
Then $x^{(t)}$ is the \efe of $[x^{(t)},v]$.
We claim that $v=w$. 
Assume that  $v\neq w$. Since $x(j)^-$ is the \efe of  $[x(j)^-,w]$ and  $x^{(t)}$ is the \efe of  $[x^{(t)},  v]$, we have $x^{(t)} \neq x(j)^-$ so that $x^{(t)} < x(j)^-$. 
Hence  
$ [x(j)^-,w]\subset [x^{(t)}, v]$ so that $w<v$.
It follows that there exists $x^{(t)} \le x'\le x(j)^-$ such that $[x',v]$ is in the right corner, which contradicts the choice of $x^{(t)}$.
\end{proof}

\subsection{Case: $[x,y_-] \leftarrow [x,y] \leftarrow [x,y_+]$}  \label{subsec:<<}
\hfill

Suppose that  $[x,y]\in \F$ with $[x,y_-], [x,y_+]\in \F$.  Set  $i\seteq i_x$.
We will omit the proofs of the following propositions in this subsection, since they are similar to those of Lemma \ref{lem:>>Varrows}, Lemma \ref{lem:>>Vej}, Lemma \ref{lem:>>Voj}, and Proposition \ref{prop:>>Varrow_generic}.
\begin{proposition} \label{prop:<<Varrows}
Let  $\sfc_{i,j}<0$.
We have
\eqn
\Vi_j=\Vi_j^e \sqcup \Vi_j^o ,  \qt{and} \  \Vo_j=\Vo_j^e \sqcup \Vo_j^o,
\eneqn
where 
\eqn 
&&\Vi_j^e= \set{[x',y']\in \F_j}{ [x',y'_+]\in \F,  \ \text{ $x'$ is the \efe of $[x',y']$},   \ y< y'<y'_+<y_+ } \\
&&\qquad = \set{[x',y']\in \F_j}{  \ \text{$[x',y']$ is in the left corner},   \ y(j)^+\le y'<y_+(j)^- }, \\
&&\Vi_j^o =  \set{[x', y']\in \F_j}{[x',y'_+]\in \F,  \ y'<y<y'_+<y_+ } \\
&&\qquad =\set{[x',y(j)^-]\in \F}{[x',y(j)^+]\in \F,  \ \text{and} \ y < y_+(j)^-}\\
&&\Vo_j^e = \set{[x',y']\in \F_j}{ [x'_-,y']\in \F,  \ \text{$y'$ is the \efe of $[x',y']$},  \  y<y'<y'_+<y_+ } \\
&&\qquad = \set{[x',y']\in \F_j}{  \ \text{$[x',y']$ is in the right corner},   \ y(j)^+\le y'<y_+(j)^- } \\
&&\qquad \sqcup\set{[x',x']}{[x_-',x']\in \F, \ y(j)^+\le y'<y_+(j)^-  } \\
&&\Vo_j^o =  \set{[x',y']\in \F_j}{\text{$y'$ is the \efe of $[x',y']$ },    \ y< y'< y_+<y'_+} \\
&&\qquad =\set{[x',y_+(j)^-]\in \F}{\text{$[x',y_+(j)^-]$ has the \efe $y_+(j)^-$} \ \text{and} \ y < y_+(j)^-}.
\eneqn
In particular, if $y> y_+(j)^-$, then $\Vi_j=\Vo_j=\emptyset$.
\end{proposition}

\begin{proposition}\label{propV:2main}
Assume that 
$y< y_+(j)^-$.
\bnum
\item There exist $a \le z\le w\le b$ such that 
\eqn
&&[z,y_+(j)^-]\in \F \ \text{with  \efe $y_+(j)^-$}, \text{and} \\
&& [w, y(j)^+] \in \F \ \text{with  \efe $y(j)^+$}.
\eneqn
\item We have
$\Vo_j^o=\{[z,y_+(j)^-]\}$.
\item We have
\eqn
\Vi_j^o=\begin{cases}
\{ [w,y(j)^-] \} &  \text{if} \ w < y(j)^+ \\
\emptyset  & w= y(j)^+.
\end{cases}
\eneqn 
\item 
If $z=w$, then
$\Vi_j^e=\Vo_j^e=\emptyset$.
\item
Assume that $z<w$. Then $\Vi_j^e\neq \emptyset$.
If
 \eq
\Vi_j^e=\set{[x^{(k)},y^{(k)}] }{1\le k \le t}
\eneq
such that 
\eqn
[x^{(1)},y^{(1)}]  \supsetneq  [x^{(2)},y^{(2)}] \supsetneq \cdots \supsetneq  [x^{(t)},y^{(t)}], 
\eneqn
then we have  
\eq
\Vo^e_j=\set{[x^{(k+1)},y^{(k)}] }{1\le k \le t-1} \cup \{[w,y^{(t)}]\}.
\eneq
\end{enumerate}
\end{proposition}

\subsection{Case: $[x,y_-] \leftarrow [x,y ] \rightarrow [x_-,y]$} \label{subsec:<>}
\hfill

Suppose that $[x,y] \in \Fex$ with $[x,y_-], [x_-,y]\in \F$. 
Let $i=i_x$ and fix $j \in I$ such  that  $\sfc_{i,j} <0$.
In this case, $y$ is the \efe of $[x,y]$ and $x_-$ is the \efe of $[x_-,y]$.
\begin{proposition} \label{prop:<>_vertical}
\bnum
\item  We have
\eqn 
\Vo_j&&= \Vo_j^e\sqcup \Vo_j^o,
\eneqn
where
 \eqn
 &&\Vo_j^e=\set{[x',y']\in \F_j}{\text{$[x',y']$ is in the right corner,
     $x_-(j)^+ <x'\le x(j)^-$ and $y(j)^+\le  y'$}} \\
 &&\text{and}\\
 &&\Vo_j^o=\set{[x(j)^+,y']\in \F_j}{[x(j)^-,y']\in \F, \quad  x_-(j)^+ <x , \quad  y(j)^+\le  y'}. 
\eneqn

\item
$[x',y'] \in\Vi_j$ if and only if exactly one of the following holds: 
\bna
\item $[x',y'_+]\in \F$, $[x',y']$ has the \efe $x'$,  $x_-<x'_-<x'<x$, $y<y'<y'_+<y_+$
\item $[x',y'_+]\in \F$,   $x_-<x'_- $, $y'<y<y'_+<y_+$. Equivalently, 
$[x',y']=[x',y(j)^-]$ such that $[x',y(j)^+]\in \F$ and $x_-(j)^+<x'$,  $y(j)^+<y_+$.
\item $x'_- <x_-<x'$, $y'<y<y'_+ $.  Equivalently,  $[x',y']=[x_-(j)^+,y(j)^-]$.   
\item $[x',y']$ has the \efe $x'$, $x'_-<x_-<x'<x$, $y<y'$.
Equivalently,  $[x',y']=[x_-(j)^+,y']$ has the \efe $x_-(j)^+$ , $x_-(j)^+ < x$, $y<y'$. 
\ee
\ee
\end{proposition}
Hence we have$$\Vi_j=\Vi_j(a) \sqcup \Vi_j(b) \sqcup \Vi_j(c) \sqcup \Vi_j(d)$$
where $\Vi_j(a), \Vi_j(b),\Vi_j(c)$ and $\Vi_j(d)$ are the
subsets that satisfy the above conditions (a)--(d), respectively.
\begin{proof}
(i)
 We have
\eqn 
\Vo_j=\set{[x',y']\in \F_j}{[x'_-,y']\in \F, \, \text{$y'$ is the \efe of  $[x',y']$}, \, x_-<x'_-<x, \, y < y'  <y_+}.
\eneqn
If $[x',y']\in \Vo_j$, then $[x'_-,y']\in \F$ with \efe $x'_-$.
Since $[x_-,y] \in \F$ with \efe $x_-$ and $x_-<x'_-$, Lemma \ref{lem:ydyp}
implies that $y'< y_+$.

Set
\eqn
\Vo_j^e\seteq\set{[x',y']\in \F_j}{[x'_-,y']\in \F, \, \text{$y'$ is the \efe of  $[x',y']$}, \, x_-<x'_-<x'<x, \, y < y'}.
\eneqn
Assume that $[x',y']\in \Vo_j^e$.
Then $x'<x \le y <y'$ so that $x'\neq y'$.  Hence $[x',y']$ is in the right corner by Lemma \ref{lem:rightleft}.
Since $x_-<x'_-<x'<x$ is equivalent to $x_-(j)^+ < x'\le x(j)^-$ and  $y<y'$ is equivalent to $y(j)^+\le y$,  we get the desired description of $\Vo_j^e$.

Set
\eqn
\Vo_j^o=\set{[x',y']\in \F}{[x'_-,y']\in \F, \, \text{$y'$ is the \efe of  $[x',y']$}, \, x_-<x'_-<x<x', \, y < y'}.
\eneqn
and assume that $[x',y']\in \Vo_j^o$. 
Then $x'=x(j)^+$ and $x'_- = x(j)^-$. 
Since $ x_-< x(j)^- $ is equivalent to $x_-(j)^+ <x$, we get the desired description of $\Vo_j^o$.

\mnoi
(ii) We have
\eqn
\Vi=&&\set{[x',y']}{\text{$[x',y'_+]\in \F$, $[x',y']$ has the \efe $x'$, $x_-<x'<x$, $y<y'_+<y_+$}}  \\
&&\cup\set{[x',y']}{\text{$[x',y'_+]\in \F$, $x_-<x'$, $y'<y<y'_+<y_+$}} \\
&&\cup\set{[x',y']}{\text{ $x'_- <x_-<x'$, $y'<y<y'_+$}} \\
&&\cup\set{[x',y']}{\text{ $[x',y']$ has the \efe $x'$, $x'_-<x_-<x'<x$, $y<y'_+$}} \\
=&&\set{[x',y']}{\text{$[x',y'_+]\in \F$, $[x',y']$ has the \efe $x'$, $x_-<x'_-<x'<x$, $y<y'_+<y_+$}}  \\
&&\cup\set{[x',y']}{\text{$[x',y'_+]\in \F$,  $x_-<x'_-<x'$, $y'<y<y'_+<y_+$}} \\
&&\cup\set{[x',y']}{\text{ $x'_- <x_-<x'$,  $y'<y<y'_+$}} \\
&&\cup\set{[x',y']}{\text{ $[x',y']$ has the \efe $x'$, $x'_-<x_-<x'<x$, $y<y'_+$}} \\
=&&\bl \set{[x',y']}{\text{$[x',y'_+]\in \F$, $[x',y']$ has the \efe $x'$, $x_-<x'_-<x'<x$, $y<y'_+<y_+$}}  \\
&&\cup\set{[x',y']}{\text{$[x',y'_+]\in \F$, $x_-<x'_-<x'$, $y'<y<y'_+<y_+$}} \br \\
&&\sqcup\set{[x',y']}{\text{ $x'_- <x_-<x'$,  $y'<y<y'_+$ }} \\
&&\sqcup\set{[x',y']}{\text{ $[x',y']$ has the \efe $x'$, $x'_-<x_-<x'<x$, $y<y'<y'_+$}} \\
 = &&\set{[x',y']}{\text{$[x',y'_+]\in \F$, $[x',y']$ has the \efe $x'$, $x_-<x'_-<x'<x$, $y<y'<y'_+<y_+$}}  \\
&&\sqcup\set{[x',y']}{\text{$[x',y'_+]\in \F$, $x_-<x'$, $y'<y<y'_+<y_+$}} \\
&&\sqcup\set{[x',y']}{\text{ $x'_- <x_-<x'$, $y'<y<y'_+$ }} \\
&&\sqcup\set{[x',y']}{\text{ $[x',y']$ has the \efe $x'$, $x'_-<x_-<x'<x$, $y<y'$}}. 
\eneqn
\end{proof}

\begin{corollary} \label{cor:<>x-j+>y}
If there is no $x'$ such that $x_-<x'<y$ and $i_{x'}=j$, equivalently,  $x_-(j)^+>y$, then $\Vo_j=\Vi_j=\emptyset$ and $x_-(j)^+ > y(j)^-$.
\end{corollary}

 \begin{corollary}
If there is no $y'$ such that $y<y'\le b$ and $i_{y'}=j$, equivalently $y(j)^+ > b$, then $\Vo_j=\emptyset$ and $\Vi_j(a)=\Vi(b)_j=\Vi_j(d)=\emptyset$.
\end{corollary}

\begin{lemma} \label{lem:<>largest_smallest}
 \hfill
\bnum
\item
If $y(j)^+\le b$, then there exists $u$ such that  $[u,y(j)^+]\in \F$ with \efe $y(j)^+$ and $u_-<x$.
\item 
If  $x_-(j)^+<y$, then there exists $z <y_+$ 
such that $[x_-(j)^+,z]\in \F$ with \efe $x_-(j)^+$.
\ee
\end{lemma}
\Proof
(i) follows from  Lemma \ref{lem:xj+zefe} (ii) together with the assumption  that $[x,y]$ has the \efe $y$. The inequality $u_-<x$ follows from
Lemma~\ref{lem:ydyp} (ii) since
$y$ is the \efe of $[x,y]$ and $y\le y(j)^+$.

(ii) follows from  Lemma \ref{lem:xj+zefe} (i) since $x_-$ is the \efe of $[x_-,y]$ and  $x_-(j)^+<y$.
The inequality $z<y_+$ follows from Lemma~\ref{lem:ydyp} (i).
\QED

\begin{corollary} \label{lem:<>_Vijd}
Assume that $x_-(j)^+ < y$
and let  $z < y_+$ be the element such that $[x_-(j)^+,z]\in \F$ with \efe $x_-(j)^+$. 
Then we have
\eqn
\Vi_j(d)=
\begin{cases}
\{[x_-(j)^+,z]\} & \text{if} \ x_-(j)^+<x, \ y<z, \\
\emptyset & \text{otherwise.} 
\end{cases}
\eneqn
\end{corollary}

\begin{lemma} \label{lem:<>inF}
If $[x_-(j)^+,y(j)^-]\in \F$, then $\Vo_j=\emptyset$ and $\Vi_j=\{ [x_-(j)^+,y(j)^-]\}$.
\end{lemma}
\begin{proof}
(i) Let $[x',y'] \in \Vo_j$. Then  $y< y'$ implies $y(j)^+ \le y'$ and $x_-<x'_-$ implies $x_-(j)^+ \le x'_- $.
Hence   $[x',y']$ do not commute with $[x_-(j)^+, y(j)^-]$.

\smallskip
(ii) Let $[x',y'] \in \Vi_j$ and assume that  either (a) or (b) hold. Then  $x_-<x'_-$ implies $x_-(j)^+\le x'_-$ and $y<y'_+$  implies $y(j)^+ \le y'_+$.
Hence $[x',y'_+]\in \F$ do not commute with  $[x_-(j)^+, y(j)^-]$. 

\smallskip
(iii) Let $[x',y'] \in \Vi_j$ and assume that  (d) holds.  Then
$x' < x \le y <y'$ so that $x'<y'$.  Because $[x',y']$ has the \efe $x'$,  we conclude that 
\eqn
[x'_+,y']\in \F.
\eneqn 
On the other hand, $y<y'$ implies $y(j)^+ \le y'$ and $x_-< x'$ implies  $x_-(j)^+ \le x'$.
Hence
Hence $[x'_+,y']\in \F$ do not commute with  $[x_-(j)^+, y(j)^-]$. 

\medskip
By (i), (ii), (iii) we obtain the assertion.
\end{proof}

\begin{lemma}  \label{lem:<>_notin}
Assume that $x_-(j)^+ < y$ and let 
$z<y_+$  be the element such that $[x_-(j)^+,z]\in \F$ with  \efe $x_-(j)^+$.
If $[x_-(j)^+,y(j)^-]\notin \F$,
then we have 
\eqn\text{$y(j)^-< z$ or  equivalently $y<z$.}\eneqn
\end{lemma}
\Proof
Since $[x_-,y]\in \F$ with \efe $x_-$, 
 if  $z<y(j)^-$, then $[x_-(j)^+,y(j)^-]\in \F$
by Lemma \ref{lem:xj+yj-}.
Hence $y(j)^-\le z$.  Since $[x_-(j)^+,z] \in  \F$ and $[x_-(j)^+,y(j)^-] \notin \F$, we have $y(j)^-\neq z$, as desired.
\QED

\begin{lemma} \label{lem:<> x<x-j+<y}
Assume that $x_-(j)^+ < y$. 
If $x<x_-(j)^+$, then $\Vo_j=\emptyset$ and $\Vi_j=\Vi_j(c)=\{[x_-(j)^+,y(j)^-]\}$. 
\end{lemma}
\begin{proof}
Since $x_-(j)^+>x$, there is no $x_- < x'' < x$ with $i_{x''}=j$.
Hence we have $\Vi_j(a)=\Vi_j(d)=\Vo_j=\emptyset$.

Assume that $[x',y'] \in \Vi_j(b)$. Then $y'=y(j)^{-}$.  Because $x_-< x'_- $ and  $x < x_-(j)^+$, we have $x<x'_-$. 
Thus $[x', y'_+]=[x', y(j)^+] \in \F$ do not commute with $[x,y_-]$, a contradiction.
Hence $\Vi_j(b)=\emptyset$.

Since $y$ is the \efe of $[x,y]$ and $x<x_-(j)^+$,
Lemma~\ref{lem:oppefe} (i) implies that
$z<y$.
Hence by the Lemma \ref{lem:<>_notin}, we have $[x_-(j)^+,y(j)^-]\in \F$, as desired.
\QED

\begin{proposition} \label{prop:<>generic}  
Assume   $x_-(j)^+<x$
 and  $[x_-(j)^+,y(j)^-]\not\in \F$.
\bna
\item There exists $y<z<y_+$ such that $[x_-(j)^+,z]\in \F$  with \efe $x_-(j)^+$. 
There exists $u$ such that $[u,y(j)^+]\in\F$  with  \efe $y(j)^+$ and $u_-<x$.
Moreover we have
\eqn 
x_-(j)^+ <u.
\eneqn

\item  $\Vi_j(d)=\{[x_-(j)^+,z]\}$, $\Vi_j(c)=\emptyset$.

\item
We have
\eqn
\Vi_j(b)=
\begin{cases}
\{[u,y(j)^-]\}  & \text{if} \ u< y(j)^+, \\
\emptyset  & \text{if} \ u=y(j)^+ .\\
\end{cases}
\eneqn

\item Let $x^{(1)}$ be the integer such that $[x^{(1)},z]$ is the smallest $\bold i$-box in $\F$ of the form $[x',z]$. 
Then
either $x^{(1)}=z$ or $[x^{(1)},z]$ is in the right corner.

Moreover, we have
\eqn
[u,y(j)^+] \subset [x^{(1)},z] \subsetneq [x_-(j)^+,z].
\eneqn

\item Assume that $x^{(1)}=z$. 
 Then we have
$x(j)^+=y(j)^+=u=z$, 
$\Vo_j^e=\Vi_j(a)=\Vi_j(b)=\emptyset$ and 
$\Vo_j^o=\{[x(j)^+,x(j)^+]\}$.

\item Assume that $x^{(1)}<z$, equivalently, $[x^{(1)},z]$ is in the right corner.
Let 
$RC$ be the set, ordered by inclusion, of the \iboxes $[x',y']\in \F_j$ which lies in the right corner, contain $[u,y(j)^+]$ and are contained in $[x_-(j)^+,z]$. 
Let $t$ be the cardinality of $RC$ and $[x^{(1)}, y^{(1)}],  [x^{(2)}, y^{(2)}], \ldots ,  [x^{(t)}, y^{(t)}]$ its elements enumerated in decreasing order. 

\bnum
\item 
We have $y^{(1)}=z$, $x^{(t)}=u$ and $y(j)^+\le y^{(t)}$.

\item
We have
\eqn
x^{(k)} \le y(j)^-  \qtext{for} \quad 1\le k< t.
\eneqn

\item We have
\eqn
\Vo_j=RC=\{[x^{(1)}, y^{(1)}],  [x^{(2)}, y^{(2)}], \ldots ,  [x^{(t)}, y^{(t)}] \}
\eneqn
and
\eqn
\Vi_j(a)=\{[x^{(1)}, y^{(2)}],  [x^{(2)}, y^{(3)}], \ldots ,  [x^{(t-1)}, y^{(t)}] \}. 
\eneqn
\ee
\ee

\end{proposition}

\Proof
(a) 
Since $x_-(j)^+<x \le y$, we get the first assertion by Lemma \ref{lem:<>largest_smallest} (ii) and Lemma \ref{lem:<>_notin}.
 In particular, we have  
$y(j)^+\le z\le b$.   Hence we get the  second assertion by Lemma  \ref{lem:<>largest_smallest}  (i).

\smallskip
Note that $x_-(j)^+<x\le y<z$ and hence $|[x_-(j)^+,z]|_\phi\ge 2$ .
Assume that $y(j)^+<z$. Then $[u,y(j)^+]$ is properly contained in $[x_-(j)^+,z]$. It follows that $x_-(j)^+ <u$, since 
$[x_-(j)^+,z]$ has the \efe $x_-(j)^+$.
Assume that $z=y(j)^+$. Since  $[u,y(j)^+]$ has the \efe $y(j)^+=z$, and $[x_-(j)^+,z]$ has the \efe $x_-(j)^+$, we conclude that  $[u,y(j)^+]$ is properly contained in $[x_-(j)^+,z]$, because $|[x_-(j)^+,z]|_\phi\ge 2$.  Hence $x_-(j)^+ <u$,  as desired.

\medskip
(b) 
We have $\Vi_j(c)=\emptyset$ since $[x_-(j)^+,y(j)^-]\not\in \F$.
By Lemma \ref{lem:<>_notin}, we have  $y<z$. 
It follows that  $\Vi_j(d)=\{[x_-(j)^+,z]\}$ by Lemma \ref{lem:<>_Vijd}.

\medskip
(c)
 By Proposition \ref{prop:<>_vertical}, we have 
\eqn
\Vi_j(b)=
\begin{cases}
\{[u,y(j)^-]\} & \text{ if $x_-(j)^+< u$ and $ y(j)^+<y_+$,} \\
\emptyset &\text{otherwise.}
\end{cases}
\eneqn 
Since $y<z<y_+$, we get  $ y(j)^+<y_+$.
Hence we obtain the assertion by (a).

\medskip
(d)
By (a), there exists such an element $x^{(1)}$ and $[x^{(1)},z]$ has the \efe $z$.
Since $[x_-(j)^+,z]$ has the \efe $x_-(j)^+$ and $x_-(j)^+<x\le y<z$,  
 we have
\eq \label{eq:<>x-j+<x1}
x_-(j)^+ <x^{(1)}.
\eneq
Note that $y(j)^+\le z$ since $y<z$.
If $y(j)^+=z$, then $[u,y(j)^+]=[x^{(1)},z]$ since  $ [u,y(j)^+]$ has the \efe $y(j)^+$ and $[x^{(1)},z]$ has the \efe $z$.
If $y(j)^+<z$, then $x_-(j)^+<x^{(1)}\le u$. In the both cases, we get the assertion.

\medskip 
(e) 
By (d), we have $$x^{(1)}=u=y(j)^+=z.$$
Hence $\Vi_j(b)=\emptyset$ by (c).
Every box in $\F_j$ smaller than $[x_-(j)^+,z]$ is of the form $[x',z]$  with \efe $x'$. Hence $[x',z]$ is not in the right corner or in the left corner. 
Hence $\Vo_j^e=\Vi_j(a) =\emptyset$ by Proposition \ref{prop:<>_vertical}.
Since $u_-<x$ and $x \le y<z=u$, we have $u=x(j)^+$. 
Then the $\bold i$-box $[x(j)^+,x(j)^+]$ satisfies the conditions in Proposition \ref{prop:<>_vertical} so that $\Vo_j^o=\{[x(j)^+,x(j)^+]\}$.

\medskip
(f-i) 
Since $[x_-(j)^+,z]$ has the \efe $x_-(j)^+$, the $\bold i$-box  $[x^{(1)}, y^{(1)}]$, which is  the largest one in the right corner contained in  $[x_-(j)^+,z]$, is of the form $[x'', z]$. 

Since $[u, y(j)^+]$ has the \efe $y(j)^+$, the $\bold i$-box  $[x^{(t)}, y^{(t)}]$, which is  the smallest one in the right corner containing $[u, y(j)^+]$, is of the form $[u,y'']$ with $y'' \ge y(j)^+$.  
\medskip

(f-ii)
Let $1\le k \le t$. If $ y(j)^- < x^{(k)}$, then $y(j)^+ \le x^{(k)}$. 
Since $x^{(k)}\le u \le y(j)^+$, we have $$x^{(k)}= u =y(j)^+.$$
It follows that $k=t$ since $x^{(k)} < x^{(k+1)}$ for $1\le k< t$. 
 \medskip

(f-iii)
First note that
\eqn 
\Vo_j^e&&=\set{[x',y']\in \F_j}{\text{$[x',y']$ is in the right corner},\  x_-(j)^+ < x'\le x(j)^-  ,  \ y(j)^+\le y'\le z} \\
&&=\set{[x',y']\in \F_j}{ \text{$[x',y']$ is in the right corner},\   [x,y] \subset [x',y'] \subset [x_-(j)^+,z]}.
\eneqn
Indeed if $[x',y']$ is in $\F_j$ and $x_-(j)^+ <x'$, then
$[x',y'] \subset [x_-(j)^+,z]$ and hence $y'\le z $. Thus the first equality follows from Proposition \ref{prop:<>_vertical} and the second equality is obvious.

\smallskip
Note that 
\eqn
RC&&= \set{[x',y']\in \F_j}{\text{$[x',y']$ is in the right corner}, \ [u,y(j)^+] \subset [x',y'] \subset [x_-(j)^+,z]} \\
&&= \set{[x',y']\in \F_j}{\text{$[x',y']$ is in the right corner}, \ x_-(j)^+< x' \le u, \ y(j)^+\le y'\le z },
\eneqn
where the inequality $x_-(j)^+< x' \le u$ comes from \eqref{eq:<>x-j+<x1}.

If $[x',y']\in \Vo_j^e$, then $y(j)^+\le y'$ and hence $[u,y(j)^+]\subset [x',y']$ because $[u,y(j)^+]$ has the \efe $y(j)^+$ and $[x',y']$ has the \efe $y'$.
It follows that 
\eqn
\Vo_j^e\subset RC.
\eneqn

We shall show that $\Vo_j = RC$ by dividing the cases into two.

 (Case 1) Assume that $x<u$.
Then $u_-<x$ implies that $u=x(j)^+$  so that 
\eqn
&&RC=\\
&& \set{[x',y']\in \F_j}{\text{$[x',y']$ is in the right corner}, \ x_-(j)^+< x' \le x(j)^-, \ y(j)^+\le y'\le z }\sqcup \{[x^{(t)}, y^{(t)}]\} \\
&&= \Vo_j^e\sqcup \{[x^{(t)}, y^{(t)}]\}.
\eneqn
Since $[x^{(t)}, y^{(t)}] =[u, y^{(t)}]=[x(j)^+, y^{(t)}]$ is in the right corner, we have $[x(j)^-, y^{(t)}] \in \F$. 
Since $y< y^{(t)} $, we conclude that $\Vo_j^o=\{[x^{(t)}, y^{(t)}] \}$ by Proposition \ref{prop:<>_vertical}.
Hence $RC=\Vo_j$.

\smallskip
 (Case 2) Assume that $u<x$. 
If $[x(j)^-,y']\in \F$ for some $y'$,  then $y' <y(j)^+$ since $[u,y(j)^+]$ has the \efe $y(j)^+$.
Thus $\Vo_j^o=\emptyset$ by Proposition \ref{prop:<>_vertical}. 
Since $u<x$, we have $[x,y] \subset [u,y(j)^+].$  It follows that $RC \subset \Vo_j^e$ and hence $RC= \Vo_j^e=\Vo_j$.

\smallskip

Note  that
\eqn 
\Vi_j(a)&&=\set{[x',y']\in \F_j}{\text{$[x',y']$ is in left corner},\  x_-(j)^+ < x'\le x(j)^-  ,  \ y(j)^+\le y'< z}\\
&&=\set{[x',y']\in \F_j}{\text{$[x',y']$ is in left corner},\ [x,y] \subset [x',y'] \subsetneq [x_-(j)^+,z] }.
\eneqn
Indeed, 
assume that $[x',y']\in \Vi_j(a)$.
Since $x'<x\le y<y'$, we have  $x'\neq y'$ so that $[x',y']$ is in the left corner.
Hence we have $y'<z$ because $[x_-(j)^+,z]$ contains $[x',y']$ properly and has the \efe $x_-(j)^+$.
Because $z<y_+$, one may replace the condition $y<y'<y'_+<y_+$ with  $y(j)^+\le y'<z$
to obtain the first equality.  The second equality is obvious.

We claim that 
\eqn
\Vi_j(a)=\{[x^{(1)}, y^{(2)}],  [x^{(2)}, y^{(3)}], \ldots ,  [x^{(t-1)}, y^{(t)}] \}. 
\eneqn
Indeed, the largest box in $\Vi_j(a)$ should be $[x^{(1)},y^{(2)}]$, since $[x_-(j)^+,z] \notin \Vi_j(a)$. 
Let  $y'$ be the smallest element such that $[x^{(t)},y']=[u,y']  \in\F$.  Then $y'\le y(j)^+$, since $[u,y(j)^+]\in \F$.
Because $[u,y']$ has the \efe $u$ and $[u,y(j)^+]$ has the \efe $y(j)^+$, we have $y' < y(j)^+$.
It follows that $[u,y'] \notin \Vi_j(a)$.
Thus the $\bold i$-box $[x^{(t-1)},y^{(t)}]$ is the smallest \ibox contained in $ \Vi_j(a)$, so that the claim follows.
\QED

\subsection{Case: $[x_+,y] \rightarrow [x,y ] \leftarrow [x,y_+]$} \label{subsec:><}
\hfill

Suppose that $[x,y]\in \Fex$ with $[x_+,y] , [x,y_+] \in \F$.
Then we have $x < x_+ \le  y$ so that $x<y$.
Let $i=i_x$ and fix $j \in I$ such  that  $\sfc_{i,j} <0$.

We omit the proofs of the following propositions in this subsection, since they are similar to those in subsection \ref{subsec:<>}.

\begin{proposition} \label{prop:><_horizontal}
\bnum
\item  We have
\eqn 
\Vi_j
= \Vi_j^e\sqcup \Vi_j^o,
\eneqn
where
 \eqn
&&\Vi_j^e=\set{[x',y']\in \F_j}{\text{$[x',y']$ is in the left corner}, \quad  y(j)^+ \le y' < y_+(j)^- , \quad x' \le x(j)^- } \\
&&\qtext{and}\quad \Vi_j^o=\set{[x',y(j)^-]\in \F_j}{[x',y(j)^+]\in \F, \quad   y < y_+(j)^-, \quad  x' \le x(j)^-}. 
\eneqn
\item
$[x',y'] \in\Vo_j$ if and only if exactly one of the following holds: 
\bna
\item $[x'_-,y']\in \F$, $[x',y']$ has the \efe $y'$,  $y<y'<y'_+<y_+$,  $x_-<x'_-<x'<x$.
\item $[x'_-,y']\in \F$,   $y'_+<y_+$, $x_-<x'_-<x<x'$. Equivalently, 
$[x',y']=[x(j)^+,y']$ such that $[x(j)^-,y']\in \F$ and $y'<y_+(j)^-$, $x_-<x(j)^-$.
\item  $y'<y_+<y'_+$, $x'_-<x<x'$. Equivalently, $[x',y']=[x(j)^+,y_+(j)^-]$. 
\item $[x',y']$ has the \efe $y'$, $y<y'<y_+<y'_+$, $x'<x$,  
equivalently, $[x',y']=[x' ,y_+(j)^-]$ has the \efe $y_+(j)^-$,  $y<y_+(j)^-$, $x'<x$.
\ee
\ee
\end{proposition}
Hence we have$$\Vo_j=\Vo_j(a) \sqcup \Vo_j(b) \sqcup \Vo_j(c) \sqcup \Vo_j(d),$$
where $\Vo_j(a), \Vo_j(b),\Vo_j(c)$ and $\Vo_j(d)$ are
subsets that satisfy each of the above conditions.

\begin{proposition} \hfill
\bna
\item If $y_+(j)^-<x$, then $\Vi_j=\Vo_j=\emptyset$ and $y_+(j)^-< x(j)^+$.
\item If $x<y_+(j)^-<y$, then $\Vi_j=\emptyset$ and $\Vo_j=\{[x(j)^+,y_+(j)^-]\}$.
\item If $[x(j)^+,y_+(j)^-]\in \F$, then $\Vi_j=\emptyset$ and $\Vo_j=\{[x(j)^+,y_+(j)^-]\}$.
\ee
\end{proposition}

\begin{proposition} \label{prop:><generic}
Assume   $y<y_+(j)^-$  and  $[x(j)^+,y_+(j)^-] \not\in \F$.
\bna
\item There exists $x_-<z<x$ such that $[z,y_+(j)^-]\in \F$  with \efe $y_+(j)^-$. 
There exists $u$ such that $[x(j)^-,u]\in\F$  with  \efe $x(j)^-$ and $y< u_+$.
Moreover we have
\eqn 
u < y_+(j)^-.
\eneqn

\item  $\Vo_j(d)=\{[z,y_+(j)^-]\}$, $\Vo_j(c)=\emptyset$.

\item
We have
\eqn
\Vo_j(b)=
\begin{cases}
\{[x(j)^+,u]\}  & \text{if} \ x(j)^-<u,\\
\emptyset  & \text{if} \ x(j)^-=u. \\
\end{cases}
\eneqn

\item Let $y^{(1)}$ be the element such that $[z,y^{(1)}]$ is the smallest $\bold i$-box in $\F$ of the form $[z,y']$. 
Then
either $y^{(1)}=z$ or $[z,y^{(1)}]$ is in the left corner.

We have
\eqn
[x(j)^-,u] \subset [z,y^{(1)}] \subsetneq [z, y_+(j)^-].
\eneqn

\item Assume that $y^{(1)}=z$. 
 Then we have
$x(j)^-=y(j)^-=u=z$, 
$\Vi_j^e=\Vo_j(a)=\Vo_j(b)=\emptyset$ and 
$\Vi_j^o=\{[y(j)^-,y(j)^-]\}$.

\item Assume that $z<y^{(1)}$, equivalently, $[z,y^{(1)}]$ is in the left corner.

Let 
$LC$ be the set, ordered by inclusion, of the \iboxes $[x',y']\in \F_j$ which lies in the left corner, contain $[x(j)^-,u]$ and are contained in $[z,y_+(j)^-]$. 
Let $t$ be the cardinality of $LC$ and $[x^{(1)}, y^{(1)}],  [x^{(2)}, y^{(2)}], \ldots ,  [x^{(t)}, y^{(t)}]$ its elements enumerated in decreasing order. 

\bnum
\item 
We have $x^{(1)}=z$, $y^{(t)}=u$ and $x^{(t)} \le x(j)^-$.

\item
We have
\eqn
x(j)^+  \le y^{(k)}  \qtext{for} \quad 1\le k< t.
\eneqn

\item We have
\eqn
\Vi_j=LC=\{[x^{(1)}, y^{(1)}],  [x^{(2)}, y^{(2)}], \ldots ,  [x^{(t)}, y^{(t)}] \}
\eneqn
and
\eqn
\Vo_j(a)=\{[x^{(2)}, y^{(1)}],  [x^{(3)}, y^{(2)}], \ldots ,  [x^{(t)}, y^{(t-1)}] \}. 
\eneqn
\ee
\ee

\end{proposition}

\subsection{Case:  $[x,x] \rightarrow [x_-,x]$} \label{subsec:>}
\hfill

Suppose that $[x,x]\in \Fex$ with $[x_-,x]  \in \F$. 
Let $i=i_x$ and fix $j \in I$ such  that  $\sfc_{i,j} <0$.

\begin{proposition} \label{prop:>_vertical}
\bnum
\item  We have
\eqn 
\Vo_j
= \Vo_j^e\sqcup \Vo_j^o,
\eneqn
where
 \eqn
&&\Vo_j^e=\set{[x',y']\in \F_j}{\text{$[x',y']$ is in the right corner}, \quad x_-(j)^+<x'\le x(j)^-, \ x<y'} \\
&&\qtext{and}\quad \Vo_j^o=\set{[x(j)^+,y']\in \F_j}{[x(j)^-,y']\in \F, \ x_-(j)^+<x}. 
\eneqn

\item
$[x',y'] \in\Vi_j$ if and only if exactly one of the following  holds: 
\bna
\item $[x',y'_+]\in \F$, $[x',y']$ has the \efe $x'$,  $x_-<x'_-<x'<x$, $x<y'<y'_+<x_+$.
\item $[x',y'_+]\in \F$,   $x_-<x'_-  $, $y'<x<y'_+<x_+$. Equivalently, 
$[x',y']=[x',x(j)^-]$ such that $[x',x(j)^+]\in \F$ and $x_-(j)^+<x'$,  $x(j)^+<x_+$.
\item $x'_- <x_-<x'$, $y'<x<y'_+ $. Equivalently,  $[x',y']=[x_-(j)^+,x(j)^-]$. 
\item $[x',y']$ has the \efe $x'$, $x'_-<x_-<x'<x$, $x<y'$, 
equivalently, $[x',y']=[x_-(j)^+,y']$ has the \efe $x_-(j)^+$ , $x_-(j)^+ < x$, $x<y'$. 
\ee
\ee
\end{proposition}
Hence we have$$\Vi_j=\Vi_j(a) \sqcup \Vi_j(b) \sqcup \Vi_j(c) \sqcup \Vi_j(d)$$
where $\Vi_j(a), \Vi_j(b),\Vi_j(c)$ and $\Vi_j(d)$ are
subsets that satisfy each of the above conditions.

\Proof
(i) 
We have
\eqn
&&\Vo_j = \set{[x',y']\in \F}{[x'_-,y']\in \F, \  \text{$[x',y']$ has the \efe $y'$}, \ x_-<x'_-<x, \  x<y'<x_+} \\
&&\cup \set{[x',y']\in \F}{[x'_-,y']\in \F, \  x_-<x'_-<x<x', \ y'<x_+} \\
&&= \set{[x',y']\in \F}{[x'_-,y']\in \F, \  \text{$[x',y']$ has the \efe $y'$}, \ x_-<x'_-<x'<x , \  x<y'<x_+} \\
&&\sqcup \set{[x',y]\in \F}{[x'_-,y']\in \F, \  x_-<x'_-<x<x', \ y'<x_+} \\
&&= \set{[x',y']\in \F}{[x'_-,y']\in \F, \  \text{$[x',y']$ in the right corner}, \ x_-<x'_-<x'<x , \  x<y'<x_+} \\
&&\sqcup \set{[x(j)^+,y']\in \F}{[x(j)^-,y']\in \F, \  x_-(j)^+< x, \ y'<x_+}. 
\eneqn
Note that if $[x',y']\in \Vo_j$, then $[x'_-,y']\in \F_j$ and $x_-<x'$. Because $[x_-,x] \in \F$  with \efe $x_-$, we have 
$y'<x_+$. 
Since $x_-<x'_-<x'<x$ is equivalent to $x_-(j)^+<x'\le x(j)^-$, we get the assertion.

(ii) We omit the proof since it is similar to the one in Proposition \ref{prop:<>_vertical} (ii).
\QED

 We omit the proof of the following proposition since it is similar to those in subsection \ref{subsec:<>}.
\begin{proposition} \label{prop:>_exceptional}
\hfill
\bna
\item If $x_-(j)^+>x$, then $\Vo_j=\Vi_j=\emptyset$ and $x(j)^- <x_-(j)^+ $.
\item If $[x_-(j)^+,x(j)^-]\in \F$, then $\Vo_j=\emptyset$ and $\Vi_j=\{[x_-(j)^+,x(j)^-]\}$.
\ee
\end{proposition}

\begin{proposition} \label{prop:>_generic}
Assume   $x_-(j)^+<x$
 and  $[x_-(j)^+,x(j)^-]\not\in \F$.
\bna
\item There exists $x<z<x_+$ such that $[x_-(j)^+,z]\in \F$ with \efe $x_-(j)^+$.

\item  $\Vi_j(d)=\{[x_-(j)^+,z]\}$, $\Vi_j(c)=\emptyset$.

\item 
We have
\eqn
\Vo^e_j
=\set{[x',y'] \in \F_j}{\text{$[x',y']$ is in the right corner}, \ [x,x] \subsetneq [x',y' ] \subsetneq [x_-(j)^+,z]},
\eneqn
and
\eqn
\Vi_j(a)=\set{[x',y'] \in \F_j}{\text{$[x',y']$ is in the left corner},  \ [x,x] \subsetneq [x',y' ] \subsetneq [x_-(j)^+,z]}.
\eneqn

\item
Assume that 
$[x(j)^-,w] \in \F$ for some $w > x(j)^-$ with  \efe $x(j)^-$.
Then we have 
$\Vo_j^o=\{[x(j)^+,w]\}$, and
$\Vi_j(b)=\emptyset$.

Moreover, if $\Vo_j^e= \{[x^{(1)}, y^{(1)}],  [x^{(2)}, y^{(2)}], \ldots ,  [x^{(t)}, y^{(t)}] \}$
for some $t\ge 1$ such that $[x^{(k)}, y^{(k)}] \supsetneq [x^{(k+1)}, y^{(k+1)}] $ for $1\le k< t$, 
then $y^{(1)}=z$ and
\eqn  \Vi_j(a)=  \{[x^{(1)}, y^{(2)}],  [x^{(2)}, y^{(3)}], \ldots ,  [x^{(t-1)}, y^{(t)}] \} \cup \{[x^{(t)},w]\}.
\eneqn

\item
Assume that 
$[u, x(j)^+] \in \F$ for some $u < x(j)^+$  with  \efe $x(j)^+$.
Then we have 
$\Vo_j^o=\emptyset$, and
$\Vi_j(b)=\{[u,x(j)^-] \}$.

Moreover, if $\Vo_j^e= \{[x^{(1)}, y^{(1)}],  [x^{(2)}, y^{(2)}], \ldots ,  [x^{(t)}, y^{(t)}] \}$
for some $t\ge 1$ such that $[x^{(k)}, y^{(k)}] \supsetneq [x^{(k+1)}, y^{(k+1)}] $ for $1\le k< t$, 
then $y^{(1)}=z$, $x^{(t)}=u$ and  
\eqn  \Vi_j(a)=  \{[x^{(1)}, y^{(2)}],  [x^{(2)}, y^{(3)}], \ldots ,  [x^{(t-1)}, y^{(t)}] \}.
\eneqn
\ee
\end{proposition}

\Proof 
(a)
Since  $[x,x],[x_-,x]\in \F$,  there exists $z<x_+$ such that $[x_-(j)^+,z]\in \F$  with \efe $x_-(j)^+$ by Lemma \ref{lem:xj+zefe}. 

Since $[x_-,x]\in \F$ with \efe $x_-$, 
 if  $z<x(j)^-$, then $[x_-(j)^+,x(j)^-]\in \F$ by Lemma \ref{lem:xj+yj-}.
Hence $x(j)^-\le z$.  Since $[x_-(j)^+,z] \in  \F$ and $[x_-(j)^+,x(j)^-] \notin \F$, we have $x(j)^-\neq z$ so that $x<z$ as desired.

\mnoi
(b) follows from (a) and Proposition \ref{prop:>_vertical}.

\mnoi
(c)\ 
If $[x',y'] \in  \Vo^e_j$, then we have $x_-(j)^+<x'$ so that $[x',y' ] \subsetneq [x_-(j)^+,z]$, since $x'$ and $x_-(j)^+$ have the same color, and it is clear that $[x,x] \subsetneq [x',y' ]$. 
Conversely, assume that $[x',y']$ is in the right corner and $ [x,x] \subsetneq [x',y' ] \subsetneq [x_-(j)^+,z]$. 
Since $[x',y']$ is in the right corner and it is properly contained  in $[x_-(j)^+,z]$, we have $x_-(j)^+ <x$. The condition $[x,x] \subsetneq [x',y' ] $ implies that $x'<x$ and $x<y'$,  since $[x,x]$ and $[x',y' ]$ have different colors.  
\smallskip

If  $[x',y'] \in  \Vi_j(a)$, then $x'<x<y'$ and  hence $[x',y']$ is in the left corner.  The condition $x_-<x'_-<x'<x$ is equivalent to $x_-(j)^+<x'\le x(j)^-$. Thus $[x_-(j)^+,z]$ properly contains $[x',y']$. Moreover, because  $[x',y']$ is in the left corner and  $[x_-(j)^+,z]$ has the \efe $x_(j)^+$, we conclude that $y'<z$.  It follows that 
\eqn
\Vi_j(a)=\set{[x',y'] \in \F_j}{\text{$[x',y']$ is in the left corner},  \ x_-(j)^+<x'\le x(j)^-, \ x(j)^+ \le y'<z},
\eneqn
which is equivalent to the description in the proposition.

\medskip 
(d)
By the assumption we have $[x(j)^+,w] \in \F$.
By Proposition \ref{prop:>_vertical}, we have $\Vo_j^o=\{[x(j)^+,w]\}$. 
Since $[x(j)^+,w]$ does not commute with any box of the form $[x',x(j)^-]$, we have $\Vi_j(b)=\emptyset$ by 
by Proposition \ref{prop:>_vertical}.
\smallskip

Let $\Vo_j^e= \{[x^{(1)}, y^{(1)}],  [x^{(2)}, y^{(2)}], \ldots ,  [x^{(t)}, y^{(t)}] \}$. Since $[x_-(j)^+,z]$ has the \efe $x_-(j)^+$, we have $y^{(1)}=z$. 
By the assumption , $[x(j)^-,w]$ is the smallest box in $\F_j$ that contains $[x,x]$. 
Hence a box in $\F_j$ contains $[x,x]$ if and only if it contains $[x(j)^-,w]$.
Note that  $[x^{(k)}, y^{(k+1)}]$ $(1\le  k < t)$  is the largest box in the left corner,  which  is contained in $[x^{(k)},y^{(k)}]$.

Let  $v$ be the smallest element such that $[x^{(t)},v]\in \F$. 
Then $x^{(t)}$ is the \efe of $[x^{(t)},v]$.
Since $x(j)^-$ is the \efe of  $[x(j)^-,w]$ and  $x^{(t)}$ is the \efe of  $[x^{(t)},  v]$, we have $x^{(t)} \neq x(j)^-$ so that $x^{(t)} < x(j)^-$. 
Hence  
$ [x(j)^-,w]\subset [x^{(t)}, v]$ so that $w\le v$.
We claim that $v=w$. 
Indeed if $v\neq w$, then $w<v$ and hence 
 there exists $x^{(t)} \le x'\le x(j)^-$ such that $[x',v]$ is in the right corner, which is a contradiction to the choice of $x^{(t)}$. It follows that $[x(j)^-,w]$ is the smallest box in the left corner in $\F_j$ and contains $[x,x]$.
Thus we obtain that 
$ \Vi_j(a)=  \{[x^{(1)}, y^{(2)}],  [x^{(2}, y^{(3)}], \ldots ,  [x^{(t-1)}, y^{(t)}] \} \cup \{[x^{(t)},w]\}$.

\medskip 
(e)
By the assumption we have $[u,x(j)^-] \in \F$.
Since no box of the form $[x(j)^+,y']$ commutes with $[u, x(j)^-] $, we obtain that $\Vo_j^o=\emptyset$.

Note that $x_-(j)^+<x<z$ and hence $x_-(j)^+<z$.
Assume that $x(j)^+<z$. Then $[u,x(j)^+]$ is properly contained in $[x_-(j)^+,z]$. It follows that $x_-(j)^+ <u$, since 
$[x_-(j)^+,z]$ has the \efe $x_-(j)^+$.
Assume that $z=x(j)^+$. Since  $[u,x(j)^+]$ has the \efe $x(j)^+=z$, and $[x_-(j)^+,z]$ has the \efe $x_-(j)^+$, we conclude that  $[u,x(j)^+]$ is properly contained in $[x_-(j)^+,z]$.  Hence $x_-(j)^+ <u$. 
It follows that $\Vi_j(b)=\{ [u,x(j)^-]\}$ by Proposition \ref{prop:>_vertical} together with (a).
\smallskip

Let $\Vo_j^e= \{[x^{(1)}, y^{(1)}],  [x^{(2)}, y^{(2)}], \ldots ,  [x^{(t)}, y^{(t)}] \}$.
Since $[x_-(j)^+,z]$ has the \efe $x_-(j)^+$, we have $y^{(1)}=z$. 
By the assumption , $[u,x(j)^+]$ is the smallest box in $\F_j$ that contains $[x,x]$. 
Hence a box in $\F_j$ contains $[x,x]$ if and only if it contains $[u,x(j)^+]$.
Note that  $[x^{(k)}, y^{(k+1)}]$ $(1\le  k < t)$  is the largest box in the left corner,  which  is contained in $[x^{(k)},y^{(k)}]$.

Since $[u, x(j)^+]$ has the \efe $x(j)^+$, the $\bold i$-box  $[x^{(t)}, y^{(t)}]$, which is  the smallest one in the right corner containing $[u, x(j)^+]$, is of the form $[u,y'']$ with $y'' \ge x(j)^+$.  In particular, we have $u=x^{(t)}$.
Let  $y'$ be the smallest element such that $[x^{(t)},y']=[u,y']  \in\F$.  Then $y'\le x(j)^+$, since $[u,x(j)^+]\in \F$.
Because $[u,y']$ has the \efe $u$ and $[u,x(j)^+]$ has the \efe $x(j)^+$, we have $y' < x(j)^+$.
It follows that $[u,y'] \notin \Vi_j(a)$.
Thus $[x^{(t-1)},y^{(t)}]$ is the smallest box contained in $ \Vi_j(a)$. It follows that 
$\Vi_j(a)=  \{[x^{(1)}, y^{(2)}],  [x^{(2)}, y^{(3)}], \ldots ,  [x^{(t-1)}, y^{(t)}] \}.$
\QED

\subsection{Case:  $[x,x] \leftarrow [x,x_+]$} \label{subsec:<}
\hfill

Suppose that $[x,x]\in \Fex$ with $[x,x_+]  \in \F$.
Let $i=i_x$ and fix $j \in I$ such  that  $\sfc_{i,j} <0$.

 We omit the proofs of the following propositions since they are similar to those in subsection \ref{subsec:><}.

\begin{proposition} \label{prop:<_vertical}
\bnum
\item  We have
\eqn 
\Vi_j
= \Vi_j^e\sqcup \Vi_j^o,
\eneqn
where
 \eqn
&&\Vi_j^e=\set{[x',y']\in \F_j}{\text{$[x',y']$ is in the left corner}, \quad x(j)^+ \le y'< x_+(j)^-, \quad x'<x} \\
&&\qtext{and}\quad \Vi_j^o=\set{[x', x(j)^-]\in \F_j}{{[x',x(j)^+]} \in\F , \ x< x_+(j)^-}. 
\eneqn

\item
$[x',y'] \in\Vo_j$ if and only if exactly one of the following holds: 
\bna
\item $[x'_-,y]\in \F$, $[x',y']$ has the \efe $y'$,  $x<y'<y'_+< x_+$, $x_-<x'_-< x'<x$.
\item $[x'_-,y']\in \F$,   $y'_+< x_+$, $x_-<x'_- < x<x'$. Equivalently, 
$[x',y']=[x(j)^+,y']$ such that $[x(j)^-,y']\in \F$ and $y'< x_+(j)^-$,  $ x_-< x(j)^-$.
\item $y'<x_+<y'_+$, $x'_-<x<x'$.  Equivalently, $[x',y']=[x(j)^+,x_+(j)^-]$.
\item $[x',y']$ has the \efe $y'$, $x<y'<x_+<y'_+$, $x'<x$, 
equivalently, $[x',y']=[x', x_+(j)^-]$ has the \efe $x_+(j)^-$ , $x < x_+(j)^-$, $ x' < x$. 
\ee
\ee
\end{proposition}
Hence we have
$$\Vo_j=\Vo_j(a) \sqcup \Vo_j(b) \sqcup \Vo_j(c) \sqcup \Vo_j(d)$$
where $\Vo_j(a), \Vo_j(b),\Vo_j(c)$ and $\Vo_j(d)$ are
subsets that satisfy each of the above conditions.

\begin{proposition}
\hfill
\bna
\item If $x_+(j)^-<x$, then $\Vi_j=\Vo_j=\emptyset$. 
\item If $[x(j)^+, x_+(j)^-]\in \F$, then $\Vi_j=\emptyset$ and $\Vo_j=\{[x(j)^+, x_+(j)^-]\}$.
\ee
\end{proposition}

\begin{proposition}
Assume $x_+(j)^->x$ 
 and  $[x(j)^+,x_+(j)^-]\not\in \F$.
\bna
\item There exists $x_-<z<x$ such that $[z, x_+(j)^-]\in \F$ with  \efe $x_+(j)^-$.

\item  $\Vo_j(d)=\{[z,x_+(j)^-]\}$, $\Vo_j(c)=\emptyset$.

\item 
We have
\eqn
\Vi^e_j
=\set{[x',y'] \in \F_j}{\text{$[x',y']$ is in the left corner}, \ [x,x] \subsetneq [x',y' ] \subsetneq [z,x_+(j)^-]},
\eneqn
and
\eqn
\Vo_j(a)=\set{[x',y'] \in \F_j}{\text{$[x',y']$ is in the right corner},  \ [x,x] \subsetneq [x',y' ] \subsetneq [z, x_+(j)^-]}.
\eneqn

\item
Assume that 
$[w,x(j)^+] \in \F$ for some $w<x(j)^+$  with \efe $x(j)^+$.
Then we have 
$\Vi_j^o=\{[w,x(j)^-]\}$, and
$\Vo_j(b)=\emptyset$.

Moreover, if $\Vi_j^e= \{[x^{(1)}, y^{(1)}],  [x^{(2)}, y^{(2)}], \ldots ,  [x^{(t)}, y^{(t)}] \}$
for some $t\ge 1$ such that $[x^{(k)}, y^{(k)}] \supsetneq [x^{(k+1)}, y^{(k+1)}] $ for $1\le k< t$, 
then $x^{(1)}=z$ and
\eqn  \Vo_j(a)=  \{[x^{(2)}, y^{(1)}],  [x^{(3)}, y^{(2)}], \ldots ,  [x^{(t)}, y^{(t-1)}] \} \cup \{[w, y^{(t)}]\}.
\eneqn

\item
Assume that 
$[x(j)^-,u] \in \F$ for some $u>x(j)^-$ with  \efe $x(j)^-$.
Then we have 
$\Vi_j^o=\emptyset$, and
$\Vo_j(b)=\{[x(j)^+,u] \}$.

Moreover, if $\Vi_j^e= \{[x^{(1)}, y^{(1)}],  [x^{(2)}, y^{(2)}], \ldots ,  [x^{(t)}, y^{(t)}] \}$
for some $t\ge 1$ such that $[x^{(k)}, y^{(k)}] \supsetneq [x^{(k+1)}, y^{(k+1)}] $ for $1\le k< t$, 
then $x^{(1)}=z$, $y^{(t)}=u$,  and
\eqn  \Vo_j(a)=  \{[x^{(2)}, y^{(1)}],  [x^{(3)}, y^{(2)}], \ldots ,  [x^{(t)}, y^{(t-1)}] \}.
\eneqn
\ee
\end{proposition}

\subsection{Example} 
 We consider the example in Subsection \ref{subsec:example_exch_matrix}.

Let $[x,y]=[14,17]$. Since $[x,y_-]=[14,14]$ and $[x_-,y]=[10,17]$ both lie in $\F$, the vertical arrows adjacent to $[x,y]$ are governed by Subsection \ref{subsec:<>}. 

Fix $j=3$. Then $x_-(j)^+=12$ and $y(j)^+=15$ so that $x_-(j)^+<x$ and $[x_-(j)^+,y(j)^+] \notin \F$.
Therefore  $\Vi_j$ and $\Vo_j$ are determined by Proposition \ref{prop:<>generic}.

Set $z=19$ and $u=13$. Then  
$[12,19]=[x_-(j)^+,z] \in \F$ with \efe $12=x_-(j)^+$
 and $[13,19]=[u,y(j)^+] \in\F$ with \efe$19=y(j)^+$.
 Note that  $u_-=12<14=x$ and $x_-(j)^+=12 < 13=u$.
Applying Proposition \ref{prop:<>generic} (b), 
$$\Vi_j(d)=\{[12,19]\} \qtext{and} \ \Vi_j(c)=\emptyset.$$ 
Since $u=13<16=y(j)^- $, we have   $$\Vi_j(b) =\{[13,16] \}$$
 by Proposition \ref{prop:<>generic} (c).
 
Finally $[13,19]$ is the smallest element in $\F$ whose left end is $z=19$. 
Hence by Proposition \ref{prop:<>generic} (f) (iii), we have
$$\Vi_j(a) =\emptyset$$
and
\eqn
\Vo_j&&=\set{[x',y']\in \F_j }{[x',y'] \ \text{is in the right corner and} \ [13,19] \subset [x',y'] \subset [12,19]} \\
&&= \{[13,19] \}. 
\eneqn
Hence, the arrows between the box $[14,17]$ and those in $\F_3$ are 
$$[12,19] \To[2] [14,17], \ [13,16] \To[2] [14,17] \qtext{and} \quad [14,17] \To[2] [13,19],$$ 
confirming those calculated in Subsection~\ref{subsec:example_exch_matrix}.

\section{Monoidal categorifications}
\subsection{Category $\cat_w$ and $\cat_\g^{[a,b]}$}
We shall review the monoidal categories $\cat_w$ for the quiver Hecke algebra case
and $\cat_\g^{[a,b]}$ for the quantum affine algebra case. See  \cite{KKOP18, KKOP22} for more details. 
\subsubsection{$\cat_w$}
Let $\bg$ be a symmetrizable Kac-Moody algebra associated with a symmetrizable Cartan matrix $\cartan=(\sfc_{i,j})_{i,j\in I_\bg}$;
 i.e., $\mathsf D \cartan$ is symmetric for a diagonal matrix $\mathsf D = {\rm diag}(d_i \mid i \in I)$ with $d_i\in\Z_{>0}$.
We fix a set $\set{\al_i}{i\in I} $ of simple roots of $\g$ and a weight lattice $\wl$. Take a $\Q$-valued symmetric bilinear form $( \cdot,\cdot)$  on $\wl$ satisfying $(\al_i,\al_j)=d_i\sfc_{i,j}$ for any $i,j\in I_\bg$.
Let $\mathsf Q^{+}\seteq\soplus_{i\in I_\bg} \Z_{\ge 0} \al_i$ be the  positive root lattice of $\bg$ and let $\weyl_\bg$  be  the Weyl group of $\bg$.

Let 
$I_\bg^\beta=\set{\nu=(\nu_1,\ldots, \nu_{\height{\beta}})}{ \al_{\nu_1}+\cdots +\al_{\nu_{\height{\beta}}}=\beta }$ for each $\beta \in \mathsf Q^+$, where $\height{\beta}=\sum_{i\in I_\bg} |b_i|$ for $\beta=\sum_{i\in I_\bg} b_i\al_i$.

Let $R^\bg(\beta)$ be the quiver Hecke algebra of type $\bg$  at $\beta $ for each $\beta \in \mathsf Q^+$ over a base field $\cor$.
 For the definition, we refer, for example, to \cite[Definition 1.8]{KKOP18}.
We denote by $R^\bg(\beta)\gmod$ the category of finite-dimensional graded $R^\bg(\beta)$-modules.
Then the category $R^\bg\gmod \seteq \soplus_{\beta\in \mathsf Q^+} R^\bg(\beta)\gmod$ becomes a monoidal category whose  tensor product is given by the \emph{convolution product}:  for $X\in R^\bg(\beta)\gmod$ and $Y\in R^\bg(\gamma)\gmod$,  $X\conv Y\seteq  R^\bg(\beta+\gamma)e(\beta,\gamma) \tens_{R^\bg(\beta)\tens_\cor R^\bg(\gamma)} (X\tens_\cor Y)$
where
$$e(\beta,\gamma)\seteq\sum_{\substack{\nu\in I^{\beta+\gamma}\\
\sum_{k=1}^{\height{\beta}} \al_{ \nu_k} =\al,\ \sum_{k=1}^{\height{\gamma}} \al_{ \nu_{k+\height{\beta}}}
 =\gamma}}\hs{-3ex}e(\nu)\hs{1ex}\in R^\bg(\beta+\gamma),$$
and $e(\nu)$ denotes the standard idempotent generator of $R^\bg$.

{\em From now on,  we  denote the convolution product $ \conv$ by  $\tens$ for the sake of  simplicity. }
We say that two simple modules $X$ and $Y$  in $R^\bg \gmod$ \emph{strongly commute} if $X\tens Y$ is again a simple module. A simple module $X$ is called \emph{real} if $X\tens X$ is simple.
If $X\in R^\bg(\beta) \gmod$, we set $\wt(X)\seteq-\beta \in \mathsf Q^-$, where $ \mathsf Q^-$ denotes the negative root lattice $-\mathsf Q^+$.

Let $\Aqn_{\Z[q^{\pm 1}]}$ be the integral form of the   \emph{unipotent  quantum coordinate ring $\Aqn$ associated with $\g$} 
 (for the precise definitions, see \cite[Section 1]{KKOP18}).
Then there exists a $\Z[q^{\pm 1}]$-algebra isomorphism  $K(R^\bg\gmod) \isoto \Aqn_{\Z[q^{\pm 1}]}$ (\cite{KL09, R08}). 
For each pair of integral weight $\la$ and $\mu$ of $\bg$
such that $\la=w\eta$ and $\mu=w'\eta$
for some dominant integral weight $\eta$ and $w,w'\in\weyl_\bg$ with $w'\le w$, 
there exists a member    $D(\la,\mu)$  of the upper global basis of  $\Aqn_{\Z[q^{\pm 1}]}$, 
called the \emph{unipotent quantum minor}.
It was shown in \cite[Proposition 4.1]{KKOP18}
that there exists a unique self-dual simple module $M(\la,\mu)$,
called the \emph{determinantial module}  in $R^\bg\gmod$ whose isomorphism class corresponds to $D(\la,\mu)$ under the isomorphism above. 
Note that $M(\la,\mu)$ is a real simple module and it admits an \emph{affinization of degree $d$}, where $d\in 2\Z_{\ge 0}$
and  $2(\al_i,\eta) \in d\Z$   for all $i\in I$  (see, \cite[Theorem 3.26]{KKOP19A}).
We say that a simple module is  \emph{\afr} if  it is real and  admits an affinization. 
Hence  $M(\la,\mu)$  is \afr.

For an $R^\bg(\beta)$-module $M$, we define
\eqn
&&\W(M)\seteq \set{\gamma \in \mathsf Q^+ \cap (\beta - \mathsf Q^+ )}{e(\gamma,\beta-\gamma)M\neq 0}, \\
&&\W^*(M)\seteq \set{\gamma \in \mathsf Q^+ \cap (\beta - \mathsf Q^+ )}{e(\beta-\gamma,\gamma)M\neq 0}.
\eneqn

We denote by $\cat_w$ the full subcategory $R^\bg\gmod$ whose objects $M$ satisfy $\W(M)\subset \mathsf Q^+ \cap w \mathsf Q^-$. 

Let $w$ be an element of the Weyl group $\weyl_\bg$.
For a reduced expression $\underline w=s_{i_1}s_{i_2}\cdots s_{i_l}$ of $w$, define
$ w_{\le k}\seteq s_{i_1}s_{i_2}\cdots s_{i_k}$  and $w_{<k}\seteq  w_{\le k-1}  $ for $1\le k\le l$. 
The family  in $R^\bg\gmod$
$$\set{S^{\bg,\underline w}_k\seteq  M(w_{\le k}\La_{i_k}, w_{<k}\La_{i_k})}{1\le k \le l}$$ is called the set of the  \emph{cuspidal modules  associated with $\underline w$}.
We have $\wt(S^{\bg,\underline w}_k)=-\beta_k$, where $\beta_k=s_{i_1}\cdots s_{i_{k-1}}\al_{i_k}$ for $1\le k \le l$. 
Then the category $\cat_w$ is the smallest full subcategory of $R^\bg \gmod$ stable under taking convolution product $\tens$, subquotients, extensions, grading shifts, and containing the cuspidal modules $\set{S^{\bg,\underline w}_k}{1 \le k\le l}$. 
 Recall that the \emph{quantum unipotent coordinate ring} $\Anw$ is the subalgebra of $\Aqn$ generated by the \emph{dual PBW generators} $\set{\iota(f^*_{\beta_k})}{1\le k \le l}$ (for the precise definition, see \cite[Section 7.1]{GY21}). The unipotent quantum minor $D(w_{\le k}\La_{i_k}, w_{<k}\La_{i_k})$ is equal to the dual PBW generator $\iota(f^*_{\beta_k})$ up to a power of $q$ (\cite[Lemma 7.6]{GY21}, \cite[Proposition 7.4]{GLS13S}) and  hence the Grothendieck ring  $K(\cat_w)$ is isomorphic to the  integral form $\Anw_{\Z[q^{\pm1}]}$ of $\Anw$.

 In the sequel, we neglect grading shifts in the category $R^\bg \gmod$.

We set $\bold i =(i_1,i_2,\ldots,i_l) $.  For each $\bold i$-box $[x,y]$ in $[1,l]$,   we set
\eqn
M^{\bg,\underline w}[x,y]\seteq \hd(S^{\bg,\underline w}_y \tens S^{\bg,\underline w}_{y_-}  \tens \cdots \tens S^{\bg,\underline w}_{x_+} \tens S^{\bg,\underline w}_{x}   ) \in \cat_w.
\eneqn

 Note that $$M^{\bg,\underline w}[x,y] = M(w_{\le y} \La_{i_y},  w_{<x}\, \La_{i_x})
\qtext{for} \ 1\le x\le y\le l.$$ 
In particular, $M^{\bg,\underline w}[x,y]$ is a real simple module in $\Cw$ with affinization of degree $(\al_{i_x},\al_{i_x})$ (see \cite[Theorem 3.26]{KKOP19A}).

\Prop[{\cite[Corollary 5.8]{KKOP24}}] \label{prop:ibox_commute}
If two $\bold i$-boxes $[x_1,y_1]$ and $[x_2,y_2]$ commute, then $M^{\bg,\underline w}[x_1,y_1]$ and 
$M^{\bg,\underline w}[x_2,y_2]$ strongly commute. 
\enprop

The short exact sequence in the the next proposition is called the \emph{T-system}.
\Prop[{\cite[Proposition 4.4 (a)]{KKOP23B}}]\label{prop:T}
 For an \ibox $[x,y]$ in $[1,l]$ with a color $i\in I_\bg$ such that
$x<y$, 
we have the following short exact sequence in $\cat_w$ $($up to grading shifts$)$:
\eqn
0\to  \tens_{j\in I_\bg\setminus\st{i}} M^{\bg,\underline w}[x(j)^+, y(j)^-]^{\tens -\sfc_{j,i}} \to M^{\bg,\underline w}[x_+,y] \tens M^{\bg,\underline w}[x, y_{-}] \to M^{\bg,\underline w}[x_+,y_-]\tens M^{\bg,\underline w}[x,y] \to 0.
\eneqn
\enprop

A pair of modules $(M,N)$ in $R^\bg \gmod$ is called \emph{unmixed} if $\W^*(M)\cap \W(N) \subset \{0\}$.
The pair $(S^{\bg, \underline w}_u, \, S^{\bg,\underline w}_v)$ is unmixed whenever $u>v$ (\cite[Lemma 2.14]{KKOP18}) and hence so is $(M^{\bg \underline w}[x',y'], M^{\bg \underline w}[x,y])$ for $1\le x\le y<x'\le y'\le l$.

\begin{lemma} \label{lem:unmixed_simple_hd_qH}
Let $X\in R(\beta)\gmod$ and $Y\in R(\gamma)\gmod$.
\bna
\item If $(X,Y)$ is an unmixed pair of objects in $R^\bg \gmod$,
  $X$ is \afr,  and $Y$ has a simple head, 
then $X \tens Y $ has a simple head.
\item If $(X,Y)$ is an unmixed pair of objects in $R^\bg \gmod$,  $Y$ is \afr,  and $X$ has a simple head, 
then $X \tens Y $ has a simple head.
\ee
\end{lemma}
\Proof
(a)
Since $(X,Y)$ is unmixed, so is $(X,\hd Y)$. It follows that $\La(X,Y) =-(\beta,\gamma)=\La(X,\hd Y)$,
where $\La$ denotes the degree of the R-matrix $\rmat{X,Y}$ (\cite[Section 2.3]{KP18}).
 Hence the assertion follows from \cite[Proposition 2.5]{KKOP23A}.

(b) Since $(X,Y)$ is unmixed, so is $(\hd X, Y)$. It follows that $\La(X,Y) =-(\beta,\gamma)=\La(\hd X,Y)$. 
 Hence the assertion follows from the opposite version of \cite[Proposition 2.5]{KKOP23A}.
\QED

For $1\le a\le b\le l$, let  $\cat_{w_{\le b}, w_{<a}}$ be the smallest full subcategory of $\cat_w$ that is stable under taking convolution product, subquotients, extensions, grading shifts, and contains the cuspidal modules $\set{S^{\bg,\underline w}_k}{a \le k\le b}$.
Note that $\cat_{w_{\le b}, w_<a}$  is the full subcategory $R^\bg\gmod$ whose objects $M$ satisfy $\W(M)\subset \mathsf Q^+ \cap w_{\le b} \mathsf Q^-$ and
$\W^*(M)\subset \mathsf Q^+ \cap w_{<a} \mathsf Q^+$. 
 Let $w'\seteq s_{i_a}s_{i_{a+1}\cdots s_{i_b}}$ and 
$\underline w'=(i_a,\ldots, i_b)$.
 Then there is a $\Q(q)$-algebra isomorphism 
 \eq\T_{w_{<a}}\seteq\T_{i_1}\circ\cdots\circ \T_{i_{a-1}} \cl   K(\cat_{w'}) \isoto K(\cat_{w_{\le b},w_{<a}}),\label{eq:Tw}
 \eneq
where  $\T_{i_k}$ are \emph{Lusztig's braid symmetries} (see \cite[Chapter 37]{Lu93} and also \cite{KKOP24.2})
such that
\eq \label{eq:Tibox} 
\T_{w_{<a}} \bl\bigl[M^{\bg,\underline w'}[x,y]\bigr]\br = \bigl[M^{\bg,\underline w} [x+a-1,y+a-1]\bigr] \quad
 \text{for $1\le x\le y\le b-a+1$}.\eneq 
For  \eqref{eq:Tibox}, see, for example, \cite[Proposition 7.1]{GLS13S}. 

\subsubsection{Root modules in $\cat_w$} \label{subsubsec:root_modules} 

Recall that if $M,N \in R^\bg\gmod$ are simple modules and
one of them is \afr, then there exists a non-zero homogeneous morphism $\rmat{M,N}\cl M\tens N \to N\tens M$ in $R^\bg \gmod$, which is unique up to
a constant multiple. 
We call it the \emph{R}-matrix.
Let $\La(M,N)$ be the homogeneous degree of $\rmat{M,N}$.
We set
\eqn
\de(M,N)&&\seteq \dfrac{1}{2}(\La(M,N)+\La(N,M)),\\
\tLa(M,N)&&\seteq \dfrac{1}{2}\bl\La(M,N)+(\wt(M),\wt(N)) \br.
\eneqn
They are non-negative integers.
For modules $M,N \in R^\bg\gmod$, we denote the head of $M\tens N$ by $M\hconv N$.
\Prop \label{prop:d=wtwt/2}
Let $L$ be an \afr simple module and $M$ a simple module in $R^\bg\gmod$. Set $d= (\wt(L),\wt(L))/2  $. 
Then we have
\bnum
\item If $\de(L,M) >0$, then $\de(L, L\hconv M) < \de(L,M)$,
\item $\tLa(L,L\hconv M) = \tLa(L,M) +d$,
\item $\tL(M,L)\le \tLa(L\hconv M,L) \le \tLa(M,L)+d$,
\item if $\de(L,M) >0$, then we have $\tLa(L\hconv M,L)< \tLa(M,L)+d$
\ee 
\end{proposition}
\Proof
(i) follows from \cite[Corollary 3.18]{KKOP19A}. 
(ii) follows from $\La(L,L\hconv M) =\La(L,M)$. 
The first inequality in (iii) comes from \cite[Theorem 2.11 (ii)]{KKOP23A} and the second inequality comes from $\La(L\hconv M, L) \le \La(M,L)$.  
(iv) follows from $\La(L\hconv M,L) < \La(M,L)$ in \cite[Lemma 3.17]{KKOP19A}. 
\QED

\Def
We say that a real simple module $L$ in $R\gmod$ is a \emph{root module} if 
$\mathsf d_L \seteq (\wt(L),\wt(L))/2\in \Z_{>0}$ and $L$ admits an affinization of degree  $2 \, \mathsf d_L= (\wt(L),\wt(L))$.
\end{definition}
Hence if $L$ is a root module, then  by \cite[Lemma 3.11]{KKOP19A} we have
\eq \label{eq:root_de_tLa}
\de(L,M), \quad \tLa(L,M) \in \Z_{\ge 0}\mathsf d_L \qtext{for any simple module $M$}.
\eneq

\Prop \label{prop:root_module}
Let $L$ be a root module, and $M$ a simple module in $R^\bg\gmod$. If $\de(L,M)>0$, then we have 
\eqn
&&\tLa(L,L\hconv M) = \tLa(L,M)+ \mathsf{d}_L, \\
&&\tLa(L\hconv M,L) = \tLa(M,L) \qtext{and} \quad \La(L\hconv M, L) = \La(M,L)-2 \mathsf{d}_L,\\
&& \de(L \hconv M, L) = \de(M,L) -\mathsf{d}_L.
\eneqn
\end{proposition}
\Proof
The first assertion follows from the definition. The second follows from Proposition \ref{prop:d=wtwt/2} (iii), (iv) and \eqref{eq:root_de_tLa}.
The third follows from the second.
\QED

\Lemma \label{lem:M_root}
For any $w\in \weyl_\bg$ and $i \in I_\bg$ such that $ws_i > w$, the determinantial module $M(ws_i \La_i, w\La_i)$ is a root module. 
\enlemma
\Proof
By  \cite[Theorem 3.26]{KKOP19A}, the determinantial module $L=M(ws_i \La_i, w\La_i)$  admits an affinization of degree  $(\al_i,\al_i)= (\wt(L),\wt(L))$. Hence the assertion follows.
\QED

\begin{lemma} \label{lem:de_Sy+} 
Let $\underbar w=s_{i_1}\cdots s_{i_l}$ be a reduced expression of $w$ and let $1\le x \le y < y_+\le l $. Then we have 
$$\de(S^{\bg,\underline w}_{y_+}, M^{\bg,\underline w}[x,y])=d_{i_x}.$$
\end{lemma}
\Proof
We write
$\fM[x,y]=M^{\bg,\underline w}[x,y]$ and $\bS_{y_+}=S^{\bg,\underline w}_{y_+}$
for simplicity.
We shall first show that
\eq \label{eq:SM}
\de(\bS_{y_+}, \fM[x,y])>0.\eneq
We have
$
\de(\fM[x,y],\fM[x_+,y_+])>0$ by  the  T-system (Proposition~\ref{prop:T}).
Since $\fM[x_+,y_+]\simeq \bS_{y_+}\hconv\fM[x_+,y]$, we have
\eqn
\de(\fM[x,y],\fM[x_+,y_+])\le
\de(\fM[x,y],\bS_{y_+})+\de(\fM[x,y],\fM[x_+,y])
=\de(\fM[x,y],\bS_{y_+}),
\eneqn
where the last equality follows from Proposition~\ref{prop:ibox_commute}.
Hence we have obtained \eqref{eq:SM}.

Since $\bS_{y_+}$ is a root module by Lemma \ref{lem:M_root}, 
and  $$d_{i_x} = (\al_{i_x},\al_{i_x})/{2}=\bl\wt(\bS_{y_+}),
  \wt(\bS_{y_+})\br/{2},$$
we have 
$$\de(\bS_{y_+}, \fM[x,y])-d_{i_x} =\de(\bS_{y_+}, \bS_{y_+}\hconv\fM[x,y])=\de (\bS_{y_+},\fM[x,y_+])$$
by Proposition \ref{prop:root_module}. 
Since $\de (\bS_{y_+},\fM[x,y_+])=0$ by Proposition \ref{prop:ibox_commute}, we get the desired equality.
\QED

\subsubsection{$\cat_\g^{[a,b]}$}
Let $\g$ be an affine Kac-Moody algebra and $\uqpg$ the
corresponding quantum affine algebra.
Let  $\cat_\g$ be the category of finite-dimensional integrable modules over $\uqpg$ and 
let $\cat_\g^0 \subset \cat_\g$ be the Hernandez-Leclerc category  of  finite-dimensional $\uqpg$ modules (see \cite[Section 2.2]{KKOP22}).  
 Note that $\cat_\g$ and $\cat_\g^0$ are rigid monoidal categories.
 We denote by $\RD(M)$ the right dual of $M$.

 For non-zero $M,N \in \cat_\g$ such that the universal R-matrix $\Runiv_{M,N_z}$  is \emph{rationally renormalizable}, we have a distinguished non-zero morphism $\rmat{M,N}: M\tens N \to N\tens M$ in $\cat_\g$, called the \emph{R-matrix}, and an integer valued invariant $\La(M,N)$, which plays a similar role of the degree of R-matrix for quiver Hecke algebra modules  (see \cite[Section 2.2]{KKOP19C} for details).

If $M$ and $N$ are simple, then we set $$\de(M ,N ) = \dfrac{1}{2}(\La(M,N)+\La(\RD^{-1}(M),N))=\dfrac{1}{2}(\La(M,N)+\La(N,M)). $$

Let $\cartan = (\sfc_{i,j})_{i,j\in I_\bg}$ be the Cartan matrix of a  finite-type  simply-laced simple  Lie algebra $\bg$.
Let $\mathcal D\seteq \set{\sL_i}{i\in I_\bg}$ be a \emph{strong duality datum in $\cat_\g^0$ associated with $\bg$}; 
that is, $\mathcal D$ is a family of real simple modules in $\cat_\g^0$ such that 
\bnum
\item   $\de(\sL_i, \RD^k (\sL_i)) = \delta(k=\pm 1)$ for any $k\in \Z$,
\item $\de(\sL_i, \RD^k (\sL_j)) = -\delta(k=0) \sfc_{i,j}$ for $i,j\in I_\bg$ with $i\neq j$ and for any $k\in \Z$.
\ee
Then there exists a faithful, exact, monoidal functor $\F_{\mathcal D}\cl R^\bg\gmod \to \cat_\g^0$ such that $\F_{\mathcal D}(L(i)) \simeq \sL_i$ for $i\in I_\g$, where $L(i)$ denotes a unique one-dimensional simple module in $R^\bg(\al_i)\gmod$.

Recall that a  triple $\mathcal Q=(\Delta_{\gf}, \sigma,\xi)$  is called a \emph{$\rm Q$-datum for $\g$ } where
\begin{enumerate}
 \item $\Delta_{\gf}$  is the Dynkin diagram of $\gf$ with the set of vertices $\If$, where $\gf$ is the  simply-laced finite type Lie algebra canonically associated to $\g$ (\cite{KKOP20A}),
\item $\sigma$ is an automorphism on $\Delta_{\gf}$ which yields the subdiagram $\Delta_{\g_0}$ inside the Dynkin diagram $\Delta_\g$ of $\g$, which is obtained by removing the $0$-node,
\item  $\xi$ is a function from $\If$ to $\Z$, called a \emph{height function} on $(\Delta, \sigma)$ satisfying certain conditions (see \cite[Definition 6.1]{KKOP22}).
\end{enumerate}
For the precise definition and properties of $\rm Q$-data, see \cite[Section 6]{KKOP22}. 
For each $\rm Q$-data of $\g$, we obtain  a natural strong duality datum of $\mathcal D_{\mathcal Q}$ in $\cat_\g^0$ associated with $\gf$ (see \cite[Theorem 6.12]{KKOP22}).
We say that a \emph{duality datum $\mathcal D$ arises from a $\rm Q$-datum $\mathcal Q$} if $\mathcal D=\mathcal D_{\mathcal Q}$.

Let $\mathcal D$ be an arbitrary strong duality datum in $\cat_\g^0$. 
We choose a reduced expression $\underline w_0=s_{i_1}s_{i_2}\cdots s_{i_r}$ of the longest element $w_0$ of the Weyl group $\weyl_\bg$ of $\bg$.
Let $\widehat {\underline w}_0$ be a unique extension of the function $\underline w_0: [1,r] \to I_\bg$ to $\Z=(-\infty, \infty)$  satisfying $\widehat {\underline w}_0(k+r)=(\widehat {\underline w}_0(k))^*$ for any $k\in \Z$, where $*$ denotes the involution of $I_\bg$ defined by $\al_{\im^*}=-w_0(\al_\im)$. 

For each $k\in \Z$ we define
\bnum
\item $S^{\mathcal D,\widehat {\underline w}_0 }_k: = \F_{\mathcal D}(S^{\bg, \underline w_0}_k)$ for  $1\le k\le r$,
\item $S^{\mathcal D,\widehat {\underline w}_0 }_{k+r}\seteq  \RD(S^{\mathcal D,\widehat {\underline w}_0 }_k)$ for $k\in\Z$.
\ee
The modules $S^{\mathcal D,\widehat {\underline w}_0 }_k$ are called the \emph{ affine cuspidal modules in $\cat_\g^0$ corresponding to $(\mathcal D, \widehat {\underline w}_0)$}.

For an interval $[a,b]$ in $\Z$, we define the category $\cat^{\mathcal D,\widehat {\underline w}_0,[a,b] }_\g$ as  the smallest full subcategory of $\cat_\g^0$ that is stable under taking tensor product $\tens$, subquotients, extensions, and contains the affine cuspidal modules $\set{S^{\mathcal D,\widehat {\underline w}_0 }_k}{a \le k\le b}$.  If there is no risk of confusion, we may simply denote it by $\cat^{[a,b] }_\g $.

Set $$\bold i \seteq (\widehat {\underline w}_0(k))_{k\in \Z}.$$

For each $\bold i$-box $[x,y]$,   we set
\eqn
M^{\mathcal D,\widehat {\underline w}_0 }[x,y]\seteq \hd(S^{\mathcal D,\widehat {\underline w}_0 }_y \tens S^{\mathcal D,\widehat {\underline w}_0 }_{y_-}  \tens \cdots \tens S^{\mathcal D,\widehat {\underline w}_0 }_{x_+} \tens  S^{\mathcal D,\widehat {\underline w}_0 }_{x}   ) \in \cat_\g^0.
\eneqn

\Prop[{\cite[Theorem 4.21]{KKOP22}}]\hfill \label{prop:ibox_commute_qa}
\bna
\item For any $\bold i$-box $[x,y]$,  $M^{\mathcal D,\widehat {\underline w}_0 }[x,y]$
  is a real simple module in $\cat_\g^0$.
\item If two $\bold i$-boxes $[x_1,y_1]$ and $[x_2,y_2]$ commute, then $M^{\mathcal D,\widehat {\underline w}_0 }[x_1,y_1]$ and 
$M^{\mathcal D,\widehat {\underline w}_0 }[x_2,y_2]$ strongly commute. 
\ee
\enprop

The following short exact sequence is called a \emph{T-system of $\bold i$-boxes}.
\Prop[{\cite[Theorem 4.25]{KKOP22}}]
For an \ibox with color $i$ such that $x<y$,  we have the following short exact sequence in $\cat_\g^0$:
\eqn
&&0\to  \tens_{j\in I_\bg\setminus\st{i}} M^{\mathcal D,\widehat {\underline w}_0 } [x(j)^+, y(j)^-]^{\tens -\sfc_{j,i}} \to
M^{\mathcal D,\widehat {\underline w}_0 }[x_+,y] \tens M^{\mathcal D,\widehat {\underline w}_0 }[x, y_{-}]\\*
&&\hs{45ex}\to M^{\mathcal D,\widehat {\underline w}_0 }[ x_+,y_-]\tens M^{\mathcal D,\widehat {\underline w}_0 }[x,y] \to 0.
\eneqn
When $y=x_+$, we understand $M^{\mathcal D,\widehat {\underline w}_0 }[x_+,y_-]=\one$.
\enprop

Recall that a pair  $(X,Y)$ of  modules  in $\cat_\g^0$ is called \emph{unmixed} (respectively,  \emph{strongly unmixed})  if 
$$\de(\RD(X),Y)=0 \quad  (\text{respectively,} \  \de(\RD^k(X),Y)=0 \qtext{for } \ k \in \Z_{\ge 1}).$$
The pair $(S^{\mathcal D,\widehat {\underline w}_0 }_u, S^{\mathcal D,\widehat {\underline w}_0 }_v)$ is strongly unmixed if $u>v$ (\cite[Proposition 5.7 (ii)]{KKOP20C}) and the pair $(M^{\bg \underline w}[x',y'], M^{\bg \underline w}[x,y])$ is strongly unmixed  for $x\le y<x'\le y'$ (\cite[Proposition 4.15]{KKOP22}).

\begin{lemma} \label{lem:unmixed_simple_hd_qa}
Let $X$ and $Y_i$ $(i=1,\ldots,n)$ be simple objects in $\cat_\g$.
\bnum 
\item If $(X,Y_i)$ is unmixed for $1\le i\le n$  
and $Y_1\tens \cdots \tens Y_n$ has a simple head, then  $X \tens Y_1\tens \cdots \tens Y_n$ has a simple head.

\item If $(Y_i,X)$ is \ unmixed for $1\le i\le n$  and $Y_1\tens \cdots \tens Y_n$ has a simple head, then  $Y_1\tens \cdots \tens Y_n\tens X$ has a simple head.
\ee
\end{lemma}
\Proof
(i) By \cite[Proposition 4.5(ii)]{KKOP19C} and \cite[Proposition 4.1 (ii)]{KKOP19C}, it is enough to show that $\sum_{k=1}^n\La(X,Y_k)=\La(X,\hd( Y_1\tens \cdots \tens Y_n)) $.

Recall that $\La(X,S)=\La(S,\D X)$ for any simple $S$ (\cite[Proposition 3.18]{KKOP19C}).
Hence we have
\eqn
\La\bl X,\hd( Y_1\tens \cdots \tens Y_n)\br
&&=\La\bl\hd( Y_1\tens \cdots \tens Y_n),\D X\br\\
&&=\sum_{k=1}^n\La(Y_k,\D X)=\sum_{k=1}^n\La(X,Y_k).
\eneqn
Here, the second equality follows from the fact that $Y_k$ commutes with $\D X$ for any $k$, by \cite[Lemma 4.3]{KKOP19C}. 

\snoi
(ii) can be proved in a similar way as (i). 
\QED

\subsubsection{Mutations }
{}From now on, we take
one of the following choices 
for a monoidal category $\cat$, a sequence $\bold{i}$, and an interval $[a,b]$

\eq
&&\left\{ \parbox{70ex}{
\bna
\item Let $\mathbf g$ be a symmetrizable Kac-Moody algebra and
  $w$ an element of the Weyl group of $\mathbf g$.
  Let  $\cartan=(\sfc_{i,j})_{i,j\in I}$ be the Cartan matrix of $\mathbf g$.
  Let $R$ be a quiver Hecke algebra associated with $\cartan$.
  Fix a reduced expression $\underline w=s_{i_1}s_{i_2}\cdots s_{i_l}$ of $w$, and we  
take $\bold i \seteq (i_1,i_2,\ldots, i_l) $. 
\item
  We take  $a,b$ such that $1\le a\le b\le l$. 
\item  
  Let $\cat$ be the monoidal category
   $\cat_{w_{\le b}, w_{<a}}\subset R\gmod$. 
\item Set  $S_k \seteq S^{\bg,\underline w}_k$ for $1\le k\le l$ and $\bM[x,y]\seteq M^{\bg,\underline w}[x,y]$ for an $\bold i$-box $[x,y]$ in  $[a,b]$.
\ee
} \right.
\label{eq:Cw}
\eneq

\eq
&&\left\{ \parbox{70ex}{
\bna
\item  Let $\g$ be an affine Kac-Moody algebra and $\uqpg$ the corresponding quantum affine algebra.
Let $\mathbf g$ be a  simply-laced finite type complex simple Lie algebra. Let  $\cartan=(\sfc_{i,j})_{i,j\in I}$ be the Cartan matrix of $\mathbf g$.
Fix  a  reduced expression  $\underline w_0 =s_{\im_1} s_{\im_2} \cdots s_ {\im_r}$  of the longest element $w_0$ of $\mathbf g$ 
and let $\widehat {\underline w_0}$ be the  extension of $\underline w_0$ to $\Z$  satisfying $\widehat {\underline w}_0(k+r)=(\widehat {\underline w}_0 (k))^*$ for $k\in \Z$. 
We take $\bold i = \widehat {\underline w}_0$.
Fix a strong duality data $\mathcal D=\set{\sL_i}{i\in I_\bg}$ associated with $\mathbf g$
that arises from a $Q$-data $\mathcal Q$ for $\g$. 
\item $[a,b]$ is an arbitrary interval in $\Z$.
\item Let $\cat$ be the monoidal category
  $\cat^{\mathcal D, \widehat {\underline w}_0,[a,b]}_{\g}$.
\item Set $S_k \seteq S^{\mathcal D,\widehat {\underline w}_0}_k$ for $1\le k\le r$ and $\bM[x,y] \seteq M^{\mathcal D,\widehat {\underline w}_0}[x,y]$ for an $\bold i$-box $[x,y]$ in  $[a,b]$.
\ee} 
\right.
\label{eq:Cab} 
\eneq

{\em In the sequel,
  for  an $\bold i$-box $[x,y]$, we write  $[x,y]$ instead of the corresponding module $\bM[x,y]$ for simplicity of notation.}

\medskip
Let us recall the (quantum) cluster algebra structure of the Grothendieck ring $K(\cat)$.
For the definition of cluster algebra , see for example, \cite[Section 7.1]{KKOP22}. 

We shall  denote by $[M]$ the isomorphism class in the Grothendieck ring of a module $M$ in $\cat$, and by $[\F]$ the set of isomorphism classes of modules in a family $\F$.

\medskip
Assume \eqref{eq:Cw}. 
First assume that $a=1$ and $b=l$ so that $\cat=\cat_w$. Then the Grothendieck ring $K(\cat)$ is isomorphic to the quantum unipotent subgroup $\Anw$.   
In \cite{GY17}, it is shown that the quantum unipotent subgroup $\Anw$ has a quantum cluster algebra structure (for the definition of quantum cluster algebra, see for example \cite[Section 2.4]{KKOP24}).  
The {\em initial quantum seed} in  \cite[Section 10.1]{GY17} (see also \cite[Theorem 7.3]{GY21}) can be described in the following way:
for  a reduced expression $\underline w=s_{i_1}\cdots s_{i_l}$ of $w$, set $\bold i= (i_1,\ldots, i_l)$.
Let $\F^{\bold i}_+ \seteq \st{\{1,k]}_{k\in[1,l]}$ be the maximal commuting family of $\bold i$-boxes corresponding to the admissible chain with  extent $[a,b]$  associated with $(1,(\RR,\RR,\ldots, \RR))$.
 Let  $\La^{\seed^{w,\bold i}_+}$ be the $l\times l$ skew-symmetric matrix  such that $(\La^{\seed^{w,\bold i}_+})_{k,j}=\La(\{1,k ], \{1, j])$ for $1\le k, j \le l$. 
Then $\seed^{\underline w,\bold i}_+\seteq(\F^{\bold i}_+, \tB(\F^{\bold i}_+) ;\F^{\bold i}_+,(\F^{\bold i}_+)_\ex)$ is a  monoidal seed in $\cat_w$  (see Definition \ref{def: monoidal seed} below) and 
the triple  $( [\F^{\bold i}_+],-\La^{\seed^{w,\bold i}_+},\tB(\F^{\bold i}_+) )$ gives a quantum seed 
for the  quantum cluster algebra $\Anw$ under the isomorphism $\Anw\simeq \Q(q^{\pm 1/2})\tens_{\Z[q^{\pm1}]}K(\cat_w)$.

 Now assume that  $1\le a\le b\le l$. 
 Let $w'\seteq s_{i_a}s_{i_{a+1}}\cdots  s_{i_b}  $ and 
$\underline w'\seteq(i_a,\ldots, i_b)$.
Via   the  $\Q(q)$-algebra isomorphism $\T_{w_{<a}}$ in \eqref{eq:Tw},
the Grothendieck ring $\Q(q^{\pm 1/2})\tens_{\Z[q^{\pm1}]}
K(\cat_{w_{\le b}, w_{<a}})$ has the quantum cluster algebra structure coming from the one of $\Q(q^{\pm 1/2})\tens_{\Z[q^{\pm1}]} K(\cat_{w'})=A_q(\n(w'))$.
By \eqref{eq:Tibox} the initial quantum cluster  $[\F^{\bold i'}_+]$ of $A_q(\n(w'))$ maps to the quantum cluster $[\F^{\bold i}_{[a,b]}]$ in  $K(\cat_{w_{\le b},w_{<a}})$, 
where $\F^{\bold i}_{[a,b]}\seteq\st{\{a,k]\mid a\le k\le b}$ 
is the maximal commuting family of $\bold i$-boxes corresponding to the admissible chain with  extent $[a,b]$ associated with the pair  $(a,(\RR,\RR,\ldots, \RR))$.
Hence $\bl[\F^{\bold i}_{[a,b]}],\tB(\F^{\bold i}_{[a,b]})\br$
a seed of the cluster algebra $K(\cat)$. We call it the initial seed of $K(\cat)$.

\smallskip

It is known that $\cat_w$  is a monoidal categorification of $\Anw$
when the quiver Hecke algebra is symmetric (and hence the Cartan matrix is symmetric)
\cite{KKKO18, Qin17}.
In particular, every (quantum) cluster monomial of $\Anw$ corresponds to a real simple module in $\cat_w$. 
However, it is still an open problem for a non-symmetric case.

\bigskip
Assume \eqref{eq:Cab}.
Let $\F_-\seteq \st{[k,b\}}_{k\in[a,b]}$ be the family of $\bold i$-boxes corresponding to the admissible chain associated with the pair  $(b,(\LL,\LL,\ldots, \LL))$.
Then $\seed^{[a,b],\Dd,\hhw}_-\seteq(\F_-, \tB(\F_-))$ is a
 monoidal seed in $\cat^{\mathcal D, \widehat {\underline w}_0,[a,b]}_{\g}$ and 
the Grothendieck ring $K(\cat^{\mathcal D, \widehat {\underline w}_0,[a,b]}_{\g})$ is isomorphic to the cluster algebra $\A([\seed^{[a,b],\Dd,\hhw}_-])$ associated with the seed
$[\seed^{[a,b],\Dd,\hhw}_-] \seteq ( [\F_-],\tB(\F_-) )$
(\cite[Theorem 8.1]{KKOP22}).
We call this seed {\em the initial seed.}

One of the main
results of \cite{KKOP22} is 
that $\cat$ is a monoidal categorification of $K(\cat)$.
In particular, for any maximal commuting family $\F$ of \iboxes in $[a,b]$,
there exists a skew-symmetrizable exchange matrix $\tB$
such that $([\F],\tB)$ is a seed in the cluster algebra $K(\cat)$.

\bigskip

In the sequel we ignore the grading shifts in case \eqref{eq:Cw}.
To be uniform, we say that a module in $\shc$ in case \eqref{eq:Cab}
is \afr if it is real.

For objects $M,N \in \cat$, we denote the head of $M\tens N$ by $M\hconv N$.
If both of them are simple and one of them is \afr, then $M\hconv N$ is simple.

For any simple module $X\in \cat$, there exists a unique
$(m_k)_{k\in [a,b]} \in\Z_{\ge^0}^{[a,b]}$ such that 
\eqn
X\simeq \hd(S_b^{\tens m_b} \tens S_{b-1}^{\tens m_{b-1}}  \tens   \cdots  \tens  S_a^{\tens m_a} ).
\eneqn
by \cite[Proposition 2.18]{KKOP18} and \cite[Theorem 6.10]{KKOP20C}.
We say that $S_u$ is a \emph{cuspidal component of $X$} if $m_u\neq 0$.

\Lemma\label{lem:comad}
Let $X_1,X_2\in \cat$ be simple modules such that
$$X_k\simeq\hd\bl S_b^{\tens m_{k,b}} \tens S_{b-1}^{\tens m_{k,b-1}}
\tens   \cdots  \tens  S_a^{\tens m_{k,a}}\br\qt{for $k=1,2$.}$$
Then the simple module 
$\hd\bl S_b^{\tens( m_{1,b}+m_{2,b})} \tens S_{b-1}  ^{\tens( m_{1,b-1}+m_{2,b-1})}\tens   \cdots  \tens  S_a^{\tens( m_{1,a}+m_{2,a})}\br$ appears once as  a composition factor of $X_1\tens X_2$, and the other composition factors of $X_1\tens X_2$  are of the form $\hd\bl S_b^{\tens n_b} \tens S_{b-1}  ^{\tens n_{b-1}}\tens   \cdots  \tens  S_a^{\tens  n_{a }}\br$ such that 
$(n_a,\cdots,n_b)$ is strictly less than $(m_{1,a}+m_{2,a},\ldots, m_{1,a}+m_{2,a} )$ in the bi-lexicographic order.
 In particular,  if $X_1$ and $X_2$ strongly commute, then
  $$X_1\tens X_2\simeq\hd\bl S_b^{\tens( m_{1,b}+m_{2,b})} \tens S_{b-1}
  ^{\tens( m_{1,b-1}+m_{2,b-1})}\tens   \cdots  \tens  S_a^{\tens( m_{1,a}+m_{2,a})}\br.$$
\enlemma

\Proof
For $\bold m = (m_a,\ldots, m_b) \in \Z^{[a,b]}_{\ge 0}$, set $P(\bold m)\seteq
 S_b^{\tens m_{b}} \tens S_{b-1}^{\tens m_{b-1}} \tens \cdots \tens S_a^{\tens m_a} $ and $V(\bold m)\seteq\hd\bl P(\bold m)\br$.
Then  in the Grothendieck ring $K(\cat)$ 
we have 
$$[P(\bold m)] = [V(\bold m)]+ \sum_{\bold m' \prec_{\rm bi} \bold m} c_{\bold m,\bold m'} [V(\bold m')] \quad \text{for some}\ c_{\bold m,\bold m'} \in \Z_{\ge 0}$$
by \cite[Theorem 6.12]{KKOP18} and   \cite[Proposition 2.15]{KKOP18}, where $\prec_{\rm bi} $ denotes the bi-lexicographic order on $\Z^{[a,b]}_{\ge 0}$.

Thus we have
$$[V(\bold m)] = [P(\bold m)]+ \sum_{\bold m' \prec_{\rm bi} \bold m} f_{\bold m,\bold m'} [P(\bold m')] \quad \text{for some}\ f_{\bold m, \bold m'} \in \Z.$$
Set $f_{\bold m, \bold m}=1$.
Set $\bold m_k=(m_{k,a},\ldots, m_{k,b})$ so that  $X_k=V(\bold m_k)$ for $k=1,2$.
Then we have
\eq \label{eq:tens_comp}
[V(\bold m_1) \tens V(\bold m_2) ] = [P(\bold m_1+\bold m_2) ]+ \sum_{\bold m_1', \bold m_2'}
f_{\bold m_1,\bold m_1'} f_{\bold m_2,\bold m_2'}  [P(\bold m_1'+\bold m'_2)],\eneq
where
the sum runs over the pairs $(\bold m'_1,\bold m'_2)$ such that
$\bold m'_1 \preceq_{\rm bi} \bold m_1$,  $\bold m'_2 \preceq_{\rm bi} \bold m_2$ and $ (\bold m_1' , \bold m_2') \neq (\bold m_1 , \bold m_2) $.
Note that for such pairs, we have
$$\bold m_1' +\bold m_2' \prec_{\rm bi} \bold m_1 +\bold m_2.$$
Hence the right hand side of  \eqref{eq:tens_comp} is of the form 
$$[V(\bold m_1+\bold m_2)] +\text{(a $\Z$-linear combination of $[V(\bold m)]$s such that   $\bold m \prec_{\rm bi} \bold m_1+\bold m_2$)},$$
which implies the assertion.
\QED

\begin{lemma} \label{lem:T_long}
Let $[x,y],[x',y']$ be $\bold i$-boxes in $[a,b]$ with the same color.
If $x\le x' \le y\le y'$, then we have\eqn
[x',y'] \hconv[x,y] \simeq [x',y]\tens [x,y'] \simeq [x,y']\tens [x',y].
\eneqn
\end{lemma}
\Proof
If $y=y'$, then   two $\bold i$-boxes commute, and hence the assertion is trivial.
Assume that $y<y'$.
Then we have
$$[x',y']\simeq[y_+,y']\hconv[x',y]
\qtq {} [x,y']\simeq [y_+,y']\hconv[x,y].$$

Hence we obtain
\eqn
&&[x',y']\tens [x,y]\monoto{} [x',y]\tens [y_+,y']\tens [x,y]
\epito [x',y]\tens[x,y'].\eneqn
Since $[y_+,y']$ is simple the composition does not vanish.
Since $[x',y]\tens[x,y']$ is simple, it is an epimorphism.
\QED

\begin{definition} \hfill \label{def: monoidal seed}
\bnum
\item Let $\K$ be an index set. We say that a family  of simple  modules
  $ \{ M_j \}_{j \in K}$ in $\cat$ is a
\emph{commuting family} if
$$M_i \tens M_j \simeq M_j \tens M_i\qt{for any $i,j\in K$.}$$

\item Let $\{ M_i\}_{i\in\K}$ be a commuting family in $\cat$ and let $\tB =
(b_{ij})_{(i,j)\in\K\times\Kex}$ be a a skew-symmetrizable exchange matrix. 
If every $M_i$ is \afr, then $\seed \seteq (\{ M_i\}_{i\in \K },\tB \KK)$  is called a \emph{monoidal seed in $\cat$}.
\item We say that a monoidal seed $\seed$ in $\cat$ with $\tB =(b_{ij})_{(i,j)\in\K\times\Kex}$ is  \emph{admissible} if it satisfies the following conditions:
\bna 
\item
  for each
$k\in\Kex$, there exists a simple object $M'_k$ of $\cat$  such that
 there is an exact sequence in $\cat$
\begin{align} \label{eq:mut_ses}
0 \to  \dtens_{b_{ik} >0} M_i^{\tens  b_{ik}} \to M_k \otimes M_k'
\to \dtens_{b_{ik} <0} M_i^{\tens  (-b_{ik})} \to 0,
\end{align}
\item  moreover, 
  $M_k'$ commutes with $M_i$ for any $i\in\K\setminus\st{k}$.
\ee
\item 
An admissible monoidal seed  $\seed=(\{M_i\}_{i\in\K},\tB\KK)$ in $\cat$ is called \emph{$\La$-admissible } 
if $\de(M_k,M'_k)=d_{k}$ for all $k\in \Kex$, 
where $M'_k$ is the object in \eqref{eq:mut_ses} and $(d_k)_{k \in \Kex}$ denotes the skew-symmetrizer of $\tB$.
\end{enumerate}
{\em Note that we do not assume that $
  [\seed]\seteq(\{[ M_i]\}_{i\in \K },\tB \KK)$
is a seed in the cluster algebra $K(\shc)$ in {\rm (ii)--(iv)}.} 

\end{definition}

For an admissible monoidal seed
$\seed=(\{M_i\}_{i\in\K},\tB\KK)$ in $\cat$,
we set
$$\mu_k(\seed)
\seteq  \bl\{M_i\}_{i\neq k}\cup\{M_k'\},\mu_k(\tB) \KK\br$$
and call it the {\em   mutation} of $\seed$ in direction $k$.
If $M_k'$  $(k\in \Kex)$ is \afr, then $\mu_k(\seed)$
is a monoidal seed in $\cat$. 

If moreover $[\seed]\seteq (\{[M_i]\}_{i\in\K},\tB\KK)$ is a seed in $K(\cat)$,
then $[\mu_k(\seed)]$ is equal to the mutation $\mu_k([\seed])$
of $[\seed]$.

 Let $\F$ be a \mcf and let $\tB(\F)$ be the skew-symmetrizable exchange matrix associated with $\F$ defined in \S\,\ref{subsec:ex_matrix}. 
 Then $(\F, \tB(\F); \F, \F_\ex)$  is a monoidal seed in $\cat$ by Proposition \ref{prop:ibox_commute} and  Proposition \ref{prop:ibox_commute_qa}.
 
We shall see in Theorem~\ref{thm:main} below that 
$([\F], \tB(\F); \F, \F_\ex)$  is a seed in the cluster algebra $K(\cat)$.

\medskip
  Recall that there is a canonical way to associate a simply-laced finite type root system to the category $\cat_\g^0$
(\cite{KKOP20A,KKOP20C}): 
for a simple  module $M \in \Ca_\g$, set $ \wt(M) \in
\Hom_\Set(\sig, \Z)$ by
\begin{align*}
\wt(M)(i,a)= \Li(M, V(\varpi_i)_a) \qquad \text{ for } (i,a) \in
\sig,
\end{align*}
 Here
$
 \sig \seteq (I_0 \times \cor^\times) / \sim, 
$
where $I_0$ denotes the index set of simple roots of the underlying finite-dimensional simple Lie algebra of $\g$, and  the equivalence relation $\sim$ is given by 
$ (i,x) \sim (j,y) $ if and only if $V(\varpi_i)_x \simeq V(\varpi_j)_y$.
The integer $\La^\infty(X,Y)$ is an invariant defined for every pair of simple modules $M, N\in \cat_\g$ (see \cite[Definition 3.6]{KKOP19C}).
Then  we have
\eq
&&\hs{5ex}\wt(S)=\wt(M)+\wt(N)
\eneq
for any simples $M$, $N$ and any simple subquotient $S$ of $M\tens N$  by \cite[Lemma 3.10]{KKOP19C}.

The lemma below is analogous to \cite[Lemma 7.13]{KKOP22}.
\Lemma \label{lem:ses_solution}
Let $\{ M_i\}_{i\in\K}$ be a commuting family of \afr simples in $\cat$.
Let $k\in K$ and assume that there exists a simple $X \in \cat$
and an exact sequence
\begin{align*}
0 \to A \to M_k \tens X \to B \to 0,
\end{align*}
such that \bna
\item $X$ strongly commutes with $M_j$ for all $j \in \K \setminus \{k\}$,
\item $\de(M_k,X)=d$ for some $d\ge 1$,
\item $A=\dtens_{i \in \K} M_i^{\tens m_i}$,
$B=\dtens_{i \in \K} M_i^{\tens n_i}$ for some
$m_i,n_i\in\Z_{\ge0}$. \ee 
Then we have 
\eqn
 \sum_{i\in \K}\wt(M_i)(m_i-n_i)=0
\qtq 
\sum_{i\in \K} \La(M_j,M_i)(m_i-n_i)=-2 d  \delta_{j,k} \qt{for any $j\in\K$.}
\eneqn 
\enlemma

\Proof 
We have
\eqn \sum_{i\in \K}m_i  \wt(M_i)=
\wt(A)=\wt(M_k)+\wt(X)=\wt(B)=\sum_{i\in \K}n_i\wt(M_i)
\eneqn 
so that 
\eqn
\sum_{i\in \K}\wt(M_i)(m_i-n_i)=0.
\eneqn

For any $j\in\K$, we have 
\eqn\sum_{i\in
\K}  \La(M_i,M_j) m_i&&=\La(A,M_j)=\La(X\hconv M_k, M_j)
=\La(X, M_j)+\La(M_k, M_j),\\
\sum_{i\in \K} \La(M_j,M_i)  n_i&&=\La(M_j,B)=\La(M_j,M_k\hconv
X)=\La(M_j, M_k)+\La(M_j,X). \eneqn
 Hence we have
 \eqn - \sum_{i\in \K}\La(M_j,M_i) (m_i-n_i)&&=
\La(X, M_j)+\La(M_k, M_j)+\La(M_j, M_k)+\La(M_j,X)\\*
&&=2\bl\de(M_j,X)+\de(M_j, M_k)\br=2 d \delta_{j,k}, \eneqn 
as desired.
\QED

The proposition below is analogous to \cite[Lemma 7.12]{KKOP22}.
\Prop  \label{prop:unique_solution}
Let $\seed=(\{M_i\}_{i\in\K},\tB,\KK)$ be a $\La$-admissible monoidal seed  in $\cat$,
and set $\La^\seed_{ij}=\Lambda(M_i,M_j)$. 
 Assume that $\K$ is a finite set. 
Then we have $\dim_\Q(\sum_{i\in \K}\Q\wt(M_i))\le |\Kfr|$,  and
for any $k\in\Kex$, $(b_{ik})_{i\in \K}$ is a solution 
$(v_i)_{i\in \K}$ in $\Q^\K$ of the equations
 \eq  \label{eq:unique_solution}
 \sum_{i\in
\K}\wt(M_i)v_i=0\qtq \sum_{i\in
\K}(\La^\seed)_{ji}v_i=-2 d_{k}  \delta_{j,k} \qt{for any $j\in\Kex$.}
\eneq 
Assume further that 
 $\dim_\Q(\sum_{i\in \K}\Q\wt(M_i))=|\Kfr|$.
Then  $(b_{ik})_{i\in \K}$ is a unique solution of \eqref{eq:unique_solution} for each $k\in \Kex$.
\end{proposition}
\Proof
First note that   $(b_{ik})_{i\in \K}$ is a solution of \eqref{eq:unique_solution}  for each $k\in \Kex$ by applying Lemma  \ref{lem:ses_solution} with $m_i=\max{(0, b_{i,k})}$ and $n_i=\max{(0, -b_{i,k})}$.

Let $f: \Q^\K\to \Q^{\Kex} \soplus(\Q \otimes_{\Z}  \mathsf Q)$ be the $\Q$-linear map given by 
$e_i \mapsto (\sum_{j\in \Kex} (\La^\seed)_{ji} e_j, \wt(M_i))$,
where $\{e_i\}_{i\in \K}$ denotes the standard basis of $\Q^\K$.
Since $(b_{ik})_{i\in \K}$ is a solution of \eqref{eq:unique_solution} for each $k\in \Kex$, $\Im(f)$ contains  $\Q^{\Kex} \soplus 0$. 
Moreover,  the image of the composition $ \Q^\K\To[f] \Q^{\Kex} \soplus(\Q \otimes_{\Z}  \mathsf Q) \epito \Q \otimes_{\Z}  \mathsf Q$ is $\sum_{i\in \K} \Q\wt(M_i)$, and hence 
 we have $\Im(f)= \Q^{\Kex} \soplus \bl \sum_{i\in \K} \Q\wt(M_i)\br$.
It follows that 
$\dim_\Q \bl \sum_{i\in \K} \Q\wt(M_i) \br=\dim_\Q \Im(f)-|\Kex| \le |\K|-|\Kex|=|\Kfr|$.

If  $\dim_\Q(\sum_{i\in \K}\Q\wt(M_i))=|\Kfr|$, then $f$ is injective and hence  $(b_{ik})_{i\in \K}$ is a unique solution of   \eqref{eq:unique_solution}  for each $k\in \Kex$.
\QED

The lemma below is analogous to \cite[Lemma 7.16]{KKOP22}.

\Lemma\label{lem:smallfrozen}
Assume \eqref{eq:Cw}.  Let $\F$ be a \mcf.
 Then we have
$$\dim_\Q\Bigl(\sum_{[x,y]\in \F}\Q\wt\bl [x,y]\br\Bigr)=|\F_\fr|=\big|\st{\im_s \in I_\bg \mid 
s\in[a,b]}\big|.$$ \enlemma 
\Proof
By Lemma~\ref{lem:Fj}, we have 
$$\sum_{[x,y]\in \F_j}\Q\wt\bl [x,y]\br
=\sum_{x\in[a,b], i_x=j}\Q \wt([x])\qt{for any $j\in I_\bg$.}$$
  Hence we have 
$$\sum_{[x,y]\in \F}\Q\wt\bl [x,y]\br
=\sum_{x\in[a,b]}\Q \wt([x]),$$
 whose dimension is $\big|\st{\im_s \in I_\bg \mid s\in[a,b]}\big|$.
\QED

The next lemma is an analogue of \cite[Proposition 7.17]{KKOP22}. 
\begin{lemma} \label{lem: b is mu}
Let $\mathfrak{C}=(\ci_k)_{1\le k\le l}$ be an admissible chain of
$\bold i$-boxes in  $[a,b]$  associated with $\bold i$, 
and let $\F=\set{\ci_k}{1\le k \le l}$ be the corresponding maximal commuting family of $\bold i$-boxes. 
Assume that $\seed\seteq\bl \F,\tB \KK\br$ is a  \Lad  monoidal seed in $\cat$ for some exchange matrix $\tB=(\mathsf b_{i,j})_{i\in \K, j\in \Kex}$ 
with $\K=\F$ and $\Kex=\F_\ex$,  a skew-symmetrizer $(d_{\ci_k})_{k\in \Kex}$ is given by $d_{k}=d_{\im_{\ci_k}}$, and $|\Kfr|= \displaystyle \dim_\Q \sum_{1\le k\le l} \Q \wt(\ci_k) $. 
If $k_0 \in \Kex$
and $\ci_{k_0}$  is a movable \ibox such that
$\tc_{k_0+1}=\ci_{k_0+1}=[x,y]$, 
then the mutation of $\seed$ in direction $k_0$ is given by the box move $B_{k_0}$ at $k_0$.  That is, we have
$$ \ci_{k_0}' =\begin{cases}
[x,y_-]  & \text{if $\ci_{k_0}=[x_+,y]$,} \\
[x_+,y]  & \text{if $\ci_{k_0}=[x,y_-]$.}
\end{cases}
$$
\end{lemma}

\begin{proof}
Assume that $\ci_{k_0}=[x_+,y]$.

Recall the T-system  \eqn
0\to  \tens_{j\in I_\bg} [x(j)^+, y(j)^-]^{\tens -\sfc_{j,i_x}} \to [x_+,y] \tens [x, y_{-}] \to [x_+,y_-]\tens [x,y] \to 0.
\eneqn
First note that 
$$\de([x_+,y],[x,y_-])=d_{i_x}.$$
Indeed, $ \tens_{j\in I_\bg} [x(j)^+, y(j)^-]^{\tens -\sfc_{j,i_x}} \not\simeq [x_+,y_-]\tens [x,y] $ implies that
 $[x_+,y] \tens [x,y_-]$ is not simple. In turn, we have 
 $\de([x_+,y],[x,y_-])>0$  by  \cite[Corollary 3.2.3]{KKKO18},  since $[x_+,y]$ and $[x,y_-]$ are \afr. 
Hence we have
\eqn
&&0< \de([x_+,y],[x,y_-])= \de(S_y\hconv [x_+,y_-],[x,y_-])\le  \de(S_y, [x,y_-] )+ \de([x_+,y_-],[x,y_-]) \\
&&= \de(S_y, [x,y_-] )+0=d_{i_y},
\eneqn
where the last equality comes from Lemma \ref{lem:de_Sy+}. Now the assertion follows from \cite[Lemma 3.11]{KKOP19A}.

By Lemma \ref{lem:ses_solution} together with the T-system above,
we conclude that $(m_k-n_k)_{k\in \K}$ is a solution of \eqref{eq:unique_solution} for $k_0$, where
\eqn
&\begin{aligned}
m_k=\begin{cases}
-\sfc_{j,i_x} & \hs{-1ex}\text{if $\ci_k=[x(j)^+,  y(j)^{-}]$ for some $j \in I_\bg \setminus \{i_x\}$,} \\
0 &\text{otherwise}
\end{cases},
&&
n_k=\begin{cases}
1 &\hs{-1ex} \text{if $\ci_k=[x_+,y_-]$ or $[x,y]$,} \\
0 &\text{otherwise}.
\end{cases}
\end{aligned}
\eneqn

Because $|\Kfr|= \dim_\Q \sum_{1\le k\le l} \Q \wt(\ci_k) $, 
we have $m_k-n_k=\mathsf b_{k,k_0}$ for any $k\in \K$ by  Proposition \ref{prop:unique_solution}.
Since $m_kn_k=0$, we have $n_k=\max(0,-\mathsf b_{k,k_0})$ for $k\in \K$.
 Hence we obtain the equality in the middle of the following:
\eqn
[x_+,y]\hconv [x,y_-]\simeq  [x_+,y_-]\tens [x,y] = \tens_{\mathsf b_{[x',y'],[x_+,y]<0}}  [x',y']^{\tens -\mathsf b_{[x',y'],[x_+,y]}} \simeq [x_+,y] \hconv [x_+,y]' .
\eneqn
Here, the first isomoprhism comes from the $T$-system, 
and  $[x_+,y]'$ denotes the mutation of $[x_+,y]$  in the admissible seed $\seed=\bl \F,\tB \KK\br$ so that the last isomorphism follows from  \eqref{eq:mut_ses}.
It follows that  $[x,y_-] \simeq [x_+,y]'$, 
as desired.
A similar proof works for the case $\ci_{k_0}=[x,y_-]$.
\end{proof}

\begin{lemma} \label{lem:unique_solution_K}
  Assume \eqref{eq:Cw}. Let $\seed=(\{x_i\}_{i\in \Kex\sqcup \Kfr}, L,\tB)$ be a quantum seed in the quantum cluster algebra $K(\cat_w)$.  Then the matrix
$\tB=(b_{i,k})_{i\in\K,k\in\Kex}$  is a unique solution of 
 \eq  \label{eq:unique_solution_al}
 \sum_{i\in
\K}\wt(x_i)b_{i,k}=0\qtq \sum_{i\in
\K}L_{ji}b_{i,k}=2 d_{k}  \delta_{j,k} \qt{for any $j,k\in\Kex$.}
\eneq 
\end{lemma}
\begin{proof}
The first equation  follows from the mutation relation in cluster algebra together with that every cluster variable in $K(\cat_w)\simeq \Anw$ is homogeneous with respect to the $\mathsf Q$-grading on $\Anw = \soplus_{\beta\in \mathsf Q^-} \Anw_\beta $. 
The second equation follows from that $(L,B)$ is a compatible pair; that is,
$ \sum_{i \in \K} l_{ji} b_{ik}  = 2d_k\delta_{j,k}$ for any $i\in \Kex$ and $j\in\K$.

Now observe that the space $\sum_{i\in \K}\Q\wt(x_i)$ is invariant under the mutation, and hence we have
$\dim_\Q(\sum_{i\in \K}\Q\wt(x_i))=|\Kfr|$ since it holds for the initial seed.

The uniqueness follows from the same argument in Proposition \ref{prop:unique_solution}.
\end{proof}

In the remainder of this section, we will  prove our main theorem.
\Th \label{thm:main}
 Let $\F$ be a maximal commuting family of $\bold i$-boxes and let $\tB(\F)=(b_{[x',y'],[x'',y'']})$ be the matrix in  \eqref{eq:ex_matrix}  associated with $\F$. Let $[x,y] \in \F_\ex$.
\bnum
\item
 There exists a simple object $\mu([x,y])\in \cat$ such that $\de([x,y], \mu([x,y]))=d_{i_x}$ and it fits into the following exact sequence in $\cat$.
  \label{main_i}
  \eqn  \label{eq:mutses}
&&0 \to  \tens_{b_{[x',y'],[x,y]}>0} [x',y']^{\tens b_{[x',y'],[x,y]}} \to [x,y] \tens \mu([x,y]) \to \tens_{b_{[x',y'],[x,y]}<0}[x',y']^{\tens -b_{[x',y'],[x,y]}}\to 0.
\eneqn
\item   The pair $([\F], \tB(\F))$ is a seed of the cluster algebra $K(\cat)$. \label{main_ii}
\item The simple object $\mu([x,y])$ in
  \eqref{main_i} is real and strongly commutes with all $[x'',y''] \in \F \setminus \{[x,y]\}$. \label{main_iii}
  \ee
\enth

We shall prove  \eqref{main_ii} and \eqref{main_iii} assuming \eqref{main_i}.
\Proof[Proof of \eqref{main_ii} and \eqref{main_iii}]
(A)\hs{1.5ex} First, assume that
there exists a skew-symmetrizable exchange matrix $\tB$
such that $([\F], \tB \KK)$  with $\K=\F$, $\Kex=\F_\ex$
is a seed in
the cluster algebra $K(\cat)$.

Since $K(\cat)$ is factorial, the cluster variable $[m]$ is prime for every $m\in \F$ by \cite[Theorem 1.3(ii)]{GLS13}.

 Let $[x'',y''] \in \F \setminus \{[x,y]\}$. Assume that  $\mu([x,y])$ and $[x'',y'']$ do not strongly commute.  Then the length of $ \mu([x,y])\tens [x'',y''] $ is $2$, since we know that the length of $[x,y]\tens \mu([x,y])\tens [x'',y''] $ is equal to $2$  by applying  the exact functor $-\tens [x'',y'']$ to the exact sequence in \eqref{main_i}.
Thus there exist  simple objects $U,V\in \cat$ and a non-split short exact sequence 
$$0 \to U \to  \mu ( [x,y] )  \tens [x'',y''] \to V \to 0.$$
Hence we have
$$0 \to [x,y] \tens U \to [x,y] \tens  \mu([x,y])\tens [x'',y''] \to [x,y] \tens  V \to 0.$$
By applying  the exact functor $-\tens [x'',y'']$ to the exact sequence in \eqref{main_i}, 
we deduce that $ [x,y] \tens \mu([x,y])\tens [x'',y''] $ has  length $2$,
and $[x,y]\tens U$ and $[x,y]\tens V$ are simple. Hence
we obtain
\eq \label{eq:xyU}
 \bl \tens_{b_{[x',y'],[x,y]}>0} [x',y']^{\tens b_{[x',y'],[x,y]}} \br \tens [x'',y'']\simeq [x,y]\tens U\ \text{or} \ [x,y]\tens V.
\eneq

Hence \eqref{eq:xyU} yields a contradiction, since the left hand side of \eqref{eq:xyU} is a tensor product of $\bold i$-boxes  belonging to $\F \setminus{[x,y] }$, and the class of any $[x',y']\in \F$
is prime in $K(\cat)$.

It follows that $\mu([x,y])$ and $[x'',y'']$  strongly commute.

Hence the quadruple $(\F, \tB(\F)  \KK)$ is a $\La$-admissible monoidal seed in $\cat$.
It follows that  the matrix $\tB(\F)$ is a unique solution of  \eqref{eq:unique_solution} by Proposition \ref{prop:unique_solution}.
Hence we obtain $\tB=\tB(\F)$ by Lemma \ref{lem:unique_solution_K} and \cite[Proposition 7.14]{KKOP22}, which implies \eqref{main_ii}.

Let us show that $\mu([x,y])$ is real.  
By applying $-\tens \mu([x,y])$ to the exact sequence in \eqref{main_i},
we obtain that $[x,y]\tens \mu([x,y]) \tens \mu([x,y])$ has length $2$.  
Thus $ \mu([x,y]) \tens \mu([x,y])$ has length less than or equal to $2$. Assume that $ \mu([x,y]) \tens \mu([x,y])$ has length $2$ with simple composition factors $Y,Z \in \cat$. 
Then  in $K(\cat)$ we have 
\eqn
&&\big[[x,y]\tens Y\big] +  \big[[x,y]\tens Z \big]     = \big[[x,y] \tens  \mu([x,y]) \tens \mu([x,y]) \big]  \\
&&= 
 \Big[ \tens_{b_{[x',y'],[x,y]}>0} [x',y']^{\tens b_{[x',y'],[x,y]}} \tens \mu([x,y])\Big] + \Big[ \tens_{b_{[x',y'],[x,y]}<0} [x',y']^{\tens -b_{[x',y'],[x,y]}} \tens \mu([x,y])\Big]. 
\eneqn
It follows that $[x,y]\tens Y$ is simple and isomorphic to
the tensor product of members in $\F\setminus\st{[x,y]}$ and $\mu([x,y])$.  It contradicts the fact that the classes of  $[x',y']\in\F$ and $\mu([x,y])$ are cluster variables and hence they are prime elements. 
Thus we conclude that $\mu([x,y])$ is real and hence we obtain \eqref{main_iii}.

Thus we  have shown  \eqref{main_ii},\eqref{main_iii} under the assumption that $\tB$ exists.

\mnoi
(B)\ Now, let us prove
that there exists a skew-symmetrizable exchange matrix $\tB$
such that $([\F], \tB \KK)$  with $\K=\F$, $\Kex=\F_\ex$
is a seed in
the cluster algebra $K(\cat)$.

Since any $[\F]$ is obtained by a succession of box moves
from the initial seed, 
we may assume that
$\F$ is obtained by the box move $B_{k_0}$ at $k_0$ 
for some $k_0\in \Kex$ from another family of $\bold i$-boxes $\F'$
such that  $([\F'], \tB(\F'))$  is a seed in $K(\cat)$.
Since $(\F', \tB(\F'))$ is a $\La$-admissible monoidal seed by (A),
we conclude that
 $([\F], \mu_{k_0}(\tB(\F')))$ is the mutation of $([\F'], \tB(\F'))$ in direction $k_0$ by Lemma~\ref{lem: b is mu} and Lemma~\ref{lem:smallfrozen}. 
 In particular, $(|\F|, \mu_{k_0}(\tB(\F')))$ is a seed in $K(\cat)$.
Hence by (A), we obtain \eqref{main_ii} and \eqref{main_iii}. 
\QED
\Rem
We conjecture that $\mu([x,y])$ in Theorem~\ref{thm:main}
has an affinization.
It is not known in case \eqref{eq:Cw} with a non-symmetric quiver Hecke algebra.
\enrem

\subsection{Strategy of the proof
  of Main Theorem}

{}The remainder of this section is devoted to proving  \eqref{main_i} in
Main Theorem \ref{thm:main}.

Let $[x,y]\in\Fex$. Set
\eqn
\fM^O&&\seteq\bigotimes_{b_{[x',y'],[x,y]}<0}[x',y']^{\tens -b_{[x',y'],[x,y]}},\\
\fM^I&&\seteq\bigotimes_{b_{[x',y'],[x,y]}>0} [x',y']^{\tens b_{[x',y'],[x,y]}}.
\eneqn

In the cases of the subsections \ref{subsec:>>}, \ref{subsec:<>}, and \ref{subsec:>}, we have $x_-<x'$ for all \iboxes $[x',y']$ appearing in $\fM^O$ or in $\fM^I$, except $[x_-,y]$. It follows  that $(\fM^I, S_{x_-})$ is unmixed but  $(\fM^O, S_{x_-})$ is not.
Similarly,  in the cases of the subsections \ref{subsec:<<}, \ref{subsec:><}, and \ref{subsec:<}, the pair $(S_{x_+}, \fM^O)$ is unmixed but $(S_{x_+}, \fM^I)$ is not. Hence we have $\fM^O\not\simeq \fM^I$.

In order to prove \eqref{main_i}, it is enough to construct a simple object
$\mu([x,y])$ which satisfies
\eq
&&\left\{\ba{l}
[x,y]\hconv\mu([x,y]) \simeq \fM^O,\\
\mu([x,y])\hconv[x,y]\simeq\fM^I,\\
\de(\mu([x,y]),[x,y])\le d_{i_x}.
\ea\right.
\label{eq:main}
\eneq

Indeed, we have $[x,y]\hconv\mu([x,y])\not\simeq \mu([x,y])\hconv[x,y]$
since  $\fM^O\not\simeq \fM^I$.  
 Since 
\eqn 
0< \de([x,y],\mu([x,y])) \le d_{i_x},
\eneqn
 we conclude that 
$\de([x,y],\mu([x,y])) = d_{i_x}$ by Lemma \cite[Lemma 3.11]{KKOP19A}. 
Thus we obtain the short exact sequence
$$0\to \mu([x,y])\hconv[x,y]\to [x,y]\tens \mu([x,y])\to[x,y]\hconv\mu([x,y])\to0$$
by \cite[Proposition 3.2.17]{KKKO18} 
and \cite[Proposition 2.11]{KKOP22}.
Thus we obtain Theorem~\ref{thm:main}~\eqref{main_i}.

\bigskip

We follow the notations in Section~\ref{sec:vertical_arrows}.
We shall divide the proof of the existence of $\mu([x,y])$
according to the configuration of adjacent horizontal arrows  as in Section \ref{sec:vertical_arrows}.

Set
\eqn
M^\Vi  &&\seteq \tens_{[x',y']\in \Vi} [x',y']^{\tens b_{[x',y'],[x,y]}} =\tens_{[x',y']\in \Vi} [x',y']^{\tens   -\sfc_{i_{x'},i} },  \\
M^\Vo &&\seteq \tens_{[x',y']\in \Vo} [x',y']^{\tens -b_{[x',y'],[x,y]}} =\tens_{[x',y']\in \Vo} [x',y']^{\tens    -\sfc_{i_{x'},i} } .
\eneqn

\subsection{Case: $[x_+,y] \to [x,y ] \to [x_-,y]$}\quad
Let $[x,y]\in \Fex$, and assume that $[x_+,y], [x_-,y]\in \F$.  
Then $x_-$  is the \efe of $[x_-,y]$,  and $x$ is the \efe of $[x,y]$.
Set  $i\seteq i_x$.
We have
$$\fM^O=[x_-,y]\tens M^\Vo,\qtq\fM^I=M^\Vi\tens [x_+,y].$$

Set 
\eqn
\mu([x,y])\seteq M^\Vo \hconv S_{x_-}.
\eneqn

Then we have
\eqn 
\de([x,y],\mu([x,y])) \le \de([x,y], M^\Vo )+\de([x,y], S_{x_-}) \le0+d_{i_x} =d_{i_x},
\eneqn
by Lemma \ref{lem:de_Sy+}.

Because $S_{x_-}$ is simple,  the composition 
\eqn
[x,y]\tens \mu( [x,y]) \monoto[]  [x,y] \tens (S_{x_-} \tens M^\Vo) \epito [x_-,y] \tens M^\Vo=\fM^O
\eneqn
does not vanish.  Hence $[x,y]\hconv\mu( [x,y])\simeq \fM^O$.

We shall show the following proposition:
\Prop\label{case1:main}
We have an epimorphism
\eq \label{eq:Vohxmhx}
&&\mu([x,y])\tens S_x \epito M^\Vi.
\eneq
\enprop
Admitting this proposition for a while, let
us prove \eqref{eq:main}.
We have a composition
\eqn
 \mu( [x,y]) \tens [x,y] \monoto\mu([x,y]) \tens (S_x \tens [x_+,y])
\epito M^\Vi\tens [x_+,y]=\fM^I,
\eneqn
which does not vanish since $S_x$ is simple. 
Thus we have
\eqn
[x,y]\hconv  \mu( [x,y])\simeq  \fM^O\qtq
\mu( [x,y]) \hconv [x,y]  \simeq \fM^I.
\eneqn

Thus we obtain \eqref{eq:main}.

Now Proposition~\ref{case1:main} is a consequence of
the following lemma and proposition.

\begin{lemma} \label{lem:Voxmx}
 The tensor product $M^\Vo \tens S_{x_-} \tens S_x$ has a simple head.
\end{lemma}
\begin{proof}
Since $x_-<x'$ for any $[x',y']\in\Vo$,
 $M^\Vo$ is a tensor product of a commuting family of
  \afr simple modules $Z$ such that $u\in[x_-+1,b]\setminus\st{x}$
  for any cuspidal component $S_u$ of $Z$. 
 
  Hence so is $M^\Vo$ by Lemma~\ref{lem:comad}, which implies that
there exist simples $X$ and $Y$ such that $M^\Vo \simeq X\hconv Y$,
every cuspidal component $S_u$ of $X$ satisfies  $u >x$ and
every cuspidal component $S_v$ of $Y$ satisfies  $x_- < v < x$.

Since every cuspidal component $S_v$ of $Y$ commutes with  $S_{x_-}$,  the tensor product $Y\tens S_{x_-}$ is simple and hence  
 $Y\tens S_{x_-}\tens S_x$ has a simple head.
Moreover,  
the pairs $(X,Y)$, $(X,S_{x_-})$ and $(X,S_x)$ are 
  unmixed,
since $(S_q, S_p)$ is unmixed whenever $q>p$.    
It follows that $X \tens (Y \tens S_{x_-}\tens S_x)$ has a simple head by Lemma \ref{lem:unmixed_simple_hd_qH} and  Lemma \ref{lem:unmixed_simple_hd_qa}. Hence $(X\hconv Y)\tens S_{x_-} \tens S_x$ has a simple head, as desired. 
\end{proof}

\begin{proposition} \label{prop:VoSSVi}
We have
\eqn 
M^\Vo \hconv (S_{x_-} \hconv S_x)  \simeq M^\Vi.
\eneqn
\end{proposition}

\begin{proof}
By the T-system,  we have
\eqn
S_{x_-} \hconv S_x \simeq \tens_{j \in I; \ \sfc_{i,j}<0,  \ x_- <x(j)^-} [x_{-}(j)^+, x(j)^-]^{ \tens -\sfc_{j,i}}. 
\eneqn

Hence in order to  prove  the proposition, it is enough to show that
for any $j\in I$  such that $\sfc_{i,j}<0$,
there exists an epimorphism
\eq \label{eq:Vohxmhx2}
 M^{\Vo^o_j}\tens M^{\Vo^e_j} \tens [x_-(j)^+,x(j)^-]^{ \tens -\sfc_{j,i}} \epito M^{\Vi^o_j}\tens M^{\Vi^e_j},
 \eneq 
where
\eqn M^{X}\seteq 
 \tens_{[x',y']\in X} [x',y']^{\tens b_{[x',y'],[x,y]}} =\tens_{[x',y']\in X} [x',y']^{\tens   -\sfc_{j,i} } \qtext{for} \  X= \Vo^o_j, \Vo^e_j,\Vi^o_j, \Vi^e_j.
\eneqn
 Indeed,  by tensoring \eqref{eq:Vohxmhx2} with respect to all $j\in I\setminus\st{i}$,
we obtain the epimorphism $M^\Vo \tens(S_{x_-} \hconv S_x)\epito M^\Vi$.
\medskip

Now, let us show the existence of an epimorphism in \eqref{eq:Vohxmhx2}.
If ${x_-(j)}^+ > x$, then $[x_{-}(j)^+, x(j)^-]\simeq \one$ and  $\Vo^o_j,=\Vo^e_j,=\Vi^o_j=\Vi^e_j=\emptyset$ by Lemma \ref{lem:>>Varrows}.
Hence we get \eqref{eq:Vohxmhx2}.  
Note that the tensor product with empty set of factors is understood as the tensor unit $\one$.

\medskip
Hence we may assume that ${x_-(j)}^+ < x$. 
Then the structure of $\Vi_j$ and $\Vo_j$ is described in
Proposition \ref{prop:>>Varrow_generic}.
There exist $w,z$ such that $ a\le w  \le z \le b$,
$[x_-(j)^+, z] \in \F$ with  \efe $x_-(j)^+$, and 
 $[x(j)^-, w] \in \F$ with  \efe $x(j)^-$.
We have
\eqn
\Vi^o_j &&=\{[x_-(j)^+, z]\},\\
Vo^o_j &&=
\begin{cases}
\{ [ x(j)^+,  w]\}  & \text{if} \ x(j)^-< w \\
 \emptyset & \text{if} \ x(j)^-=w.
\end{cases}
\eneqn
and
\eqn
\Vo^e_j&&=\set{[x^{(k)},y^{(k)}]}{1\le k\le t},\\
\Vi^e_j&&=\set{[x^{(k)},y^{(k+1)}]}{1\le k<t}\cup\st{[x^{(t)},w]}.
\eneqn 
with $y^{(1)}=z$.

Set $c\seteq -\sfc_{j,i}$ and $c'= -\sfc_{j,i} \delta( x(j)^-< w)$.
Then  we have

\scalebox{.91}{\parbox{\textwidth}{
    \eqn
&&\hs{-4ex}    \ba{l}
M^{\Vo^o_j} \tens M^{\Vo^e_j} \tens [x_-(j)^+,  x(j)^-]^{\tens c}  \\
\simeq  
[ x(j)^+,  w]^{\tens c'}   \tens [x^{(t)},y^{(t)}]^{\tens c} \tens [x^{(t-1)},y^{(t-1)}]^{\tens c} \tens \cdots \tens [x^{(2)},y^{(2)}]^{\tens c}  \tens [x^{(1)},y^{(1)}]^{\tens c} \tens [x_-(j)^+,  x(j)^-]^{\tens c} \\
\epito  
[ x(j)^+,  w]^{\tens c'}  \tens [x^{(t)},y^{(t)}]^{\tens c} \tens [x^{(t-1)},y^{(t-1)}]^{\tens c}\tens \cdots \tens [x^{(2)},y^{(2)}]^{\tens c} \tens [x^{(1)},x(j)^-]^{\tens c} \tens [x_-(j)^+,  y^{(1)}]^{\tens c}\\
\epito 
[ x(j)^+,  w]^{\tens c'}  \tens [x^{(t)},y^{(t)}]^{\tens c} \tens [x^{(t-1)},y^{(t-1)}]^{\tens c} \tens  \cdots \tens [x^{(2)},x(j)^-]^{\tens c} \tens [x^{(1)},y^{(2)}]^{\tens c} \tens [x_-(j)^+, \ y^{(1)}]^{\tens c} \\
\cdots\cdots \\
\epito [ x(j)^+,  w]^{\tens c'}  \tens [x^{(t)},x(j)^-]^{\tens c} \tens [x^{(t-1)},y^{(t)}]^{\tens c} \tens  \cdots \tens [x^{(2)},y^{(3)}] \tens [x^{(1)},y^{(2)}]^{\tens c} \tens [x_-(j)^+,  y^{(1)}]^{\tens c} \\
\epito [x^{(t)},w]^{\tens c} \tens [x^{(t-1)},y^{(t)}]^{\tens c} \tens \cdots \tens [x^{(2)},y^{(3)}]^{\tens c} \tens [x^{(1)},y^{(2)}]^{\tens c} \tens [x_-(j)^+, y^{(1)}]^{\tens c} \\
\simeq M^{\Vi^e_j} \tens M^{\Vi^o_j},\ea
\eneqn
}}

\noi
 where
the middle epimorphisms follow from  Lemma \ref{lem:T_long}. 
\end{proof}

\subsection{Case: $[x,y_-] \leftarrow [x,y ] \leftarrow [x,y_+]$}
Since the proof is similar to the preceding case, we are rather brief.

Let $[x,y]\in \Fex$ and assume that $[x,y_-], [x,y_+]\in \F$.  
Then $y$  is the \efe of $[x,y]$  and $y_+$ is the \efe of $[x,y_+]$.
Set  $i\seteq i_x$
.We have
$$\fM^O=[x,y_-]\tens M^\Vo\qtq\fM^I=M^\Vi\tens [x,y_+].$$
Set
\eqn
\mu([x,y]) \seteq S_{y_+} \hconv M^{V_i}.
\eneqn
Then we have
\eqn 
\de(\mu([x,y]),[x,y]) \le \de(S_{y_+}, [x,y])+ \de( M^\Vi,[x,y]) \le d_{i_x}.
\eneqn

We have the following non-zero composition of morphisms
\eqn
\mu([x,y]) \tens [x,y]  \monoto  M^{\Vi} \tens S_{y_+} \tens [x,y]  \to M^{\Vi} \tens  [x,y_+]=\fM^I.
\eneqn

In order to see
\eqn [x,y]\hconv \mu([x,y])\simeq [x,y_-]\tens M^{\Vo},\eneqn
it is enough to show that there exists an epimorphism
\eq
S_y\tens (S_{y_+} \hconv M^{\Vi}) \epito M^{\Vo}.
\label{eq:2main}
\eneq
Indeed, then we have
\eqn
[x,y] \tens (S_{y_+} \hconv M^{\Vi}) \monoto ([x,y_-] \tens S_y) \tens (S_{y_+} \hconv M^{\Vi}) 
\epito [x,y_-]\tens M^{\Vo}.
\eneqn

Now \eqref{eq:2main} follows from the following lemma and proposition, and we obtain
\eqref{eq:main}.
We omit the proof of the lemma below since it is similar to the one of Lemma \ref{lem:Voxmx}.
\begin{lemma}
 The tensor product  $S_y\tens S_{y_+} \tens M^{\Vi}$ has a simple head.
\end{lemma}

\begin{proposition}
$(S_y\hconv S_{y_+} ) \hconv M^{\Vi} \simeq M^{\Vo}$.
\end{proposition}
\Proof
The proof is similar to the proof of Proposition~\ref{prop:VoSSVi}
using Proposition~\ref{prop:>>Varrow_generic} instead of Proposition~\ref{propV:2main}.
By the T-system we have
\eqn
S_y\hconv S_{y_+}  
\simeq \tens_{y<y_+(j)^-} [y(j)^+,y_+(j)^-]^{\tens -\sfc_{j,i}}.
\eneqn

Set $c\seteq -\sfc_{j,i}$ and $c'= -\sfc_{j,i} \delta(  w< y(j)^+)$.
Then  we have

\scalebox{.92}{\parbox{\textwidth}{
    \eqn
&&\hs{-3ex}    \ba{l}
 [y(j)^+,  y_+(j)^-]^{\tens c} \tens  M^{\Vi^e_j} \tens M^{\Vi^o_j}  \\
\simeq  
[y(j)^+,  y_+(j)^-]^{\tens c} \tens  [x^{(1)},y^{(1)}]^{\tens c} \tens [x^{(2)},y^{(2)}]^{\tens c}\tens \cdots \tens [x^{(t-1)},y^{(t-1)}]^{\tens c}  \tens [x^{(t)},y^{(t)}]^{\tens c}   \tens [w,y(j)^-]^{\tens c'}\\
\epito 
[x^{(1)},  y_+(j)^-]^{\tens c} \tens  [y(j)^+,y^{(1)}]^{\tens c} \tens [x^{(2)},y^{(2)}]^{\tens c}\tens \cdots \tens [x^{(t-1)},y^{(t-1)}]^{\tens c}  \tens [x^{(t)},y^{(t)}]^{\tens c}   \tens [w,y(j)^-]^{\tens c'}\\
\epito 
[x^{(1)},  y_+(j)^-]^{\tens c} \tens  [x^{(2)},y^{(1)}]^{\tens c} \tens [y(j)^+,y^{(2)}]^{\tens c}\tens \cdots \tens [x^{(t-1)},y^{(t-1)}]^{\tens c}  \tens [x^{(t)},y^{(t)}]^{\tens c}   \tens [w,y(j)^-]^{\tens c'}\\
\cdots \\
\epito 
[x^{(1)},  y_+(j)^-]^{\tens c} \tens  [x^{(2)},y^{(1)}]^{\tens c} \tens [y^{(3)},y^{(2)}]^{\tens c}\tens \cdots \tens [x^{(t)},y^{(t-1)}]^{\tens c}  \tens [y(j)^+,y^{(t)}]^{\tens c}   \tens [w,y(j)^-]^{\tens c'}\\
\epito 
[x^{(1)},  y_+(j)^-]^{\tens c} \tens  [x^{(2)},y^{(1)}]^{\tens c} \tens [y^{(3)},y^{(2)}]^{\tens c}\tens \cdots \tens [x^{(t)},y^{(t-1)}]^{\tens c}  \tens [w,y^{(t)}]^{\tens c} \\
\simeq M^{\Vo^o_j} \tens M^{\Vo^e_j}.
\ea
\eneqn}}

\noi
By tensoring the epimorphisms
$[y(j)^+,  y_+(j)^-]^{\tens c} \tens  M^{\Vi^e_j} \tens M^{\Vi^o_j}
\epito M^{\Vo^o_j} \tens M^{\Vo^e_j}$ with respect to $j$ we obtain
$(S_y\hconv S_{y_+} )\tens M^{\Vi} \epito M^{\Vo}$.
\QED
Hence we obtain \eqref{eq:main}.

\subsection{Case: $[x,y_-] \leftarrow [x,y ] \to [x_-,y]$} \label{subsec:<>ses}

Let $[x,y]\in \Fex$ and assume that $[x,y_-], [x_-,y]\in \F$.  
Then $x_-$  is the \efe of $[x_-,y]$  and $y$ is the \efe of $[x,y]$.
Set  $i\seteq i_x$.

Set 
\eqn
\mu([x,y])\seteq (M^\Vo\tens  [x,y_-]) \hconv S_{x_-} \simeq  \hd(M^\Vo\tens ( [x,y_-] \tens S_{x_-}) )\simeq  M^\Vo\hconv  [x_-,y_-],
\eneqn
where the first isomorphism follows from the fact that ($M^{\Vo}$, $S_{x_-}$) is unmixed.

We have
\eqn 
\de(\mu([x,y]),[x,y]) \le \de( M^\Vo,[x,y])+ \de([x_-,y_-], [x,y]) \le d_{i_x}.
\eneqn
 
We have a composition of morphisms
\eqn
[x,y] \tens \mu([x,y]) \monoto  {[x,y]} \tens  (S_{x^-} \tens M^\Vo \tens [x,y_-] )\epito [x_-,y]\tens M^\Vo\tens [x,y_-].
\eneqn
which is non-zero since $S_{x_-}$ is simple.

\begin{lemma}
 The tensor product 
$M^\Vo \tens  [x_-,y_-] \tens [x,y]$ has a simple head.
\end{lemma}

\Proof
Let $M^\Vo = X\hconv Y$ such that every cuspidal component $S_u$ of $X$ satisfies that $y<u$ and every cuspidal component $S_v$ of $Y$ satisfies that
$ x_- <v<y$.
Then each $S_v$ commutes with $[x_-,y_-]$ 
and hence $Y\tens [x_-,y_-]$ is simple. It follows that $Y\tens [x_-,y_-] \tens [x,y]$ has a simple head since $[x,y]$ is \afr.  
Because $(X,Y)$, $(X, [x_-,y_-])$ and $(X,[x,y])$ are unmixed, 
we conclude that $X \tens Y \tens [x_-,y_-] \tens [x,y]$ has a simple head.
It follows that $(X \hconv Y) \tens [x_-,y_-] \tens [x,y]$ has a simple head, as desired.
\QED

Note that we have an epimorphism
\eqn
M^\Vo \tens  [x_-,y_-] \tens [x,y]  \epito
M^\Vo \tens ( [x_-,y_-] \hconv [x,y]  )
\simeq M^\Vo  \tens \tens_{j \in I; \ \sfc_{i,j}<0} [x_{-}(j)^+, y(j)^-]^{\tens -\sfc_{j,i}},  
\eneqn
by the $T$-system.
\begin{proposition}
There is an epimorphism
\eq \label{eq:<>epi}
M^\Vo  \tens \tens_{j \in I; \ \sfc_{i,j}<0} [x_{-}(j)^+, y(j)^-]^{\tens -\sfc_{j,i}} 
\epito M^\Vi.
\eneq
\end{proposition}
\begin{proof}
If either $[x_-(j)^+,y(j)^-] \in \F$ or  $x_-(j)^+>x$, then \eqref{eq:<>epi} holds by  Lemma \ref{lem:<>inF},
Corollary \ref{cor:<>x-j+>y} and Lemma \ref{lem:<> x<x-j+<y}.

Assume that $[x_-(j)^+,y(j)^-]\not\in \F$,  $x_-(j)^+<x$.

Set $c\seteq -\sfc_{j,i}$ and $c'= -\sfc_{j,i} \delta( u<y(j)^+)$.
Then by Proposition \ref{prop:<>generic},  we have
 \eqn
&& M^{\Vo_j}  \tens [x_-(j)^+,  y(j)^-]^{\tens c}  \\
\simeq  
&&  [x^{(t)},y^{(t)}]^{\tens c} \tens [x^{(t-1)},y^{(t-1)}]^{\tens c} \tens \cdots \tens [x^{(2)},y^{(2)}]^{\tens c}  \tens [x^{(1)},y^{(1)}]^{\tens c} \tens [x_-(j)^+,  y(j)^-]^{\tens c} \\
\epito  && 
 [x^{(t)},y^{(t)}]^{\tens c} \tens [x^{(t-1)},y^{(t-1)}]^{\tens c}\tens \cdots \tens [x^{(2)},y^{(2)}]^{\tens c} \tens [x^{(1)},y(j)^-]^{\tens c} \tens [x_-(j)^+,  y^{(1)}]^{\tens c}\\
\epito  && 
 [x^{(t)},y^{(t)}]^{\tens c} \tens [x^{(t-1)},y^{(t-1)}]^{\tens c} \tens  \cdots \tens [x^{(2)},y(j)^-]^{\tens c} \tens [x^{(1)},y^{(2)}]^{\tens c} \tens [x_-(j)^+, \ y^{(1)}]^{\tens c} \\
\cdots  && \\
\epito&&  [x^{(t)},y^{(t)}]^{\tens c} \tens [x^{(t-1)},y(j)^-]^{\tens c} \tens  \cdots \tens [x^{(2)},y^{(3)}]^{\tens c} \tens [x^{(1)},y^{(2)}]^{\tens c} \tens [x_-(j)^+, \ y^{(1)}]^{\tens c} \\
\epito &&  [x^{(t)},y(j)^-]^{\tens c'} \tens [x^{(t-1)},y^{(t)}]^{\tens c} \tens  \cdots \tens [x^{(2)},y^{(3)}] \tens [x^{(1)},y^{(2)}]^{\tens c} \tens [x_-(j)^+,  y^{(1)}]^{\tens c} \\
\simeq &&M^{\Vi(b)_j} \tens M^{\Vi(a)_j} \tens M^{\Vi(d)_j} \simeq M^{\Vi_j},
 \eneqn
as desired.
\QED

Then by the lemma above, we conclude that 
$$\mu([x,y]) \hconv [x,y] = (M^{\Vo}\hconv [x_-,y_-])\hconv [x,y] \simeq M^\Vi
\simeq\fM^I,$$
as desired.

Thus we obtain \eqref{eq:main}.

\subsection{Case: $[x_+,y] \rightarrow [x,y ] \leftarrow [x,y_+]$} 
\label{subsec:><ses}
Let $[x,y]\in \F$ and assume that $[x_+,y], [x,y_+]\in \F$.  
Then $y_+$  is the \efe of $[x,y_+]$  and $x$ is the \efe of $[x,y]$.
Set  $i\seteq i_x$.

Set 
\eqn
\mu[x,y]\seteq
S_{y_+} \hconv  ( [x_+,y] \tens M^\Vi)  \simeq \hd( S_{y_+} \tens  ( [x_+,y] \hconv M^\Vi) ) \simeq [x_+,y_+] \hconv M^{\Vi}
\eneqn
Then we have
\eqn 
\de(\mu([x,y]),[x,y]) \le \de( M^\Vi,[x,y])+ \de([x_+,y_+], [x,y])\le d_{i_x}.
\eneqn

We have a composition of morphisms
\eqn
\mu([x,y]) \tens [x,y]   \monoto    ( [x_+,y]\tens M^\Vi \tens S_{y_+} )  \tens [x,y]   \epito [x_+,y]\tens M^\Vi \tens [x,y_+]\simeq\fM^I,
\eneqn
which is non-zero since $S_{y_+}$ is simple.

\begin{lemma}
$[x,y] \tens  [x_+,y_+]\tens M^\Vi$ has a simple head.
\end{lemma}

\begin{proposition} There is an epimorphism
\eqn 
\tens_{j \in I; \ \sfc_{i,j}<0} [x(j)^+, y_+(j)^-]^{\tens -\sfc_{j,i}}  \tens M^\Vi
\epito M^\Vo.
\eneqn
\end{proposition}
\begin{proof}
We may assume that $[x(j)^+,y_+(j)^-]\not\in \F$,  $y< y_+(j)^-$. 

Set $c\seteq -\sfc_{j,i}$ and $c'= -\sfc_{j,i} \delta( x(j)^-<u)$.
Then by Proposition \ref{prop:><generic}, we have
 \eqn
&&   [x(j)^+,  y_+(j)^-]^{\tens c} \tens M^{\Vi_j} \\
\simeq  
&&  [x(j)^+,  y_+(j)^-]^{\tens c} \tens  [x^{(1)},y^{(1)}]^{\tens c} \tens [x^{(2)},y^{(2)}]^{\tens c} \tens \cdots \tens [x^{(t-1)},y^{(t-1)}]^{\tens c}  \tens [x^{(t)},y^{(t)}]^{\tens c}  \\
\epito &&  [x^{(1)},  y_+(j)^-]^{\tens c} \tens  [x(j)^+,y^{(1)}]^{\tens c} \tens [x^{(2)},y^{(2)}]^{\tens c} \tens \cdots \tens [x^{(t-1)},y^{(t-1)}]^{\tens c}  \tens [x^{(t)},y^{(t)}]^{\tens c}  \\
\epito &&  [x^{(1)},  y_+(j)^-]^{\tens c} \tens  [x^{(2)} ,y^{(1)}]^{\tens c} \tens [x(j)^+,y^{(2)}]^{\tens c} \tens \cdots \tens [x^{(t-1)},y^{(t-1)}]^{\tens c}  \tens [x^{(t)},y^{(t)}]^{\tens c}  \\
\cdots && \\
\epito &&  [x^{(1)},  y_+(j)^-]^{\tens c} \tens  [x^{(2)} ,y^{(1)}]^{\tens c} \tens \cdots \tens [x^{(t-1)},y^{(t-2)}]^{\tens c} \tens  [x(j)^+,y^{(t-1)}]^{\tens c}  \tens [x^{(t)},y^{(t)}]^{\tens c}  \\
\epito &&  [x^{(1)},  y_+(j)^-]^{\tens c} \tens  [x^{(2)} ,y^{(1)}]^{\tens c} \tens \cdots \tens [x^{(t-1)},y^{(t-2)}]^{\tens c} \tens  [x^{(t)},y^{(t-1)}]^{\tens c}  \tens [x(j)^+ ,y^{(t)}]^{\tens c'}   \\
\simeq &&M^{\Vo(d)_j} \tens M^{\Vo(a)_j}\tens M^{\Vo(b)_j} 
\simeq  M^{\Vo_j}, 
 \eneqn
as desired.
\QED

Since there is a morphism
\eqn
&&
[x,y] \tens [x_+,y_+] \tens M^\Vi   \epito  \tens_{j \in I; \ \sfc_{i,j}<0} [x(j)^+, y_+(j)^-]^{\tens -\sfc_{j,i}}  \tens M^\Vi,
\eneqn
we have 
$$[x,y]\hconv \mu([x,y])= [x,y]\hconv ([x_+,y_+]\hconv M^{\Vi}) \simeq
M^\Vo\simeq\fM^O,$$
as desired.

Thus we obtain \eqref{eq:main}.

\subsection{Case:  $[x,x] \rightarrow [x_-,x]$} 

Set
\eqn
\mu([x,x])=M^{\Vo}\hconv S_{x_-}.
\eneqn
Then we have
\eqn 
\de(\mu([x,y]),[x,y]) \le \de( M^\Vo,[x,y])+ \de(S_{x_-}, [x,y])\le d_{i_x},
\eneqn
There exists a non-zero composition of morphisms
\eqn
[x,x] \tens \mu([x,x]) \monoto   {[x,x]} \tens(S_{x_-}\tens M^{\Vo} ) \epito [x_-,x]\tens M^{\Vo}.
\eneqn

\begin{proposition}
There is  an epimorphism 
\eqn
M^{\Vo}\tens( \tens_{j\in I} [x_-(j)^+,x(j)^-]^{\tens -\sfc_{j,i}}) \epito M^{\Vi}.
\eneqn
\end{proposition}
\Proof
By Proposition \ref{prop:>_exceptional}, we may assume that $x_-(j)^+<x$ and $[x_-(j)^+,x(j)^-] \notin \F$.
Then the assertion follows from Proposition \ref{prop:>_generic}. 
Indeed, there are two cases: either $[x(j)^-,w] \in \F$ for some $w>x(j)^-$ or $[u,x(j)^+] \in \F$ for some $u<x(j)^+$.
Let  $c=\sfc_{j,i}$.

If $[x(j)^-,w] \in \F$ for some $w>x(j)^-$, then
we have

\scalebox{.91}{\parbox{\textwidth}{
    \eqn
\hs{-3ex}&&    \ba{l}
 M^{\Vo_j} \tens [x_-(j)^+,  x(j)^-]^{\tens c} = M^{\Vo^o_j} \tens M^{\Vo^e_j} \tens [x_-(j)^+,  x(j)^-]^{\tens c}  \\
\simeq  
 [ x(j)^+,  w]^{\tens c}   \tens [x^{(t)},y^{(t)}]^{\tens c} \tens [x^{(t-1)},y^{(t-1)}]^{\tens c} \tens \cdots \tens [x^{(2)},y^{(2)}]^{\tens c}  \tens [x^{(1)},y^{(1)}]^{\tens c} \tens [x_-(j)^+,  x(j)^-]^{\tens c} \\
\epito  
[ x(j)^+,  w]^{\tens c}  \tens [x^{(t)},y^{(t)}]^{\tens c} \tens [x^{(t-1)},y^{(t-1)}]^{\tens c}\tens \cdots \tens [x^{(2)},y^{(2)}]^{\tens c} \tens [x^{(1)},x(j)^-]^{\tens c} \tens [x_-(j)^+,  y^{(1)}]^{\tens c}\\
\epito   
[ x(j)^+,  w]^{\tens c}  \tens [x^{(t)},y^{(t)}]^{\tens c} \tens [x^{(t-1)},y^{(t-1)}]^{\tens c} \tens  \cdots \tens [x^{(2)},x(j)^-]^{\tens c} \tens [x^{(1)},y^{(2)}]^{\tens c} \tens [x_-(j)^+, \ y^{(1)}]^{\tens c} \\
\cdots \\
\epito [ x(j)^+,  w]^{\tens c}  \tens [x^{(t)},x(j)^-]^{\tens c} \tens [x^{(t-1)},y^{(t)}]^{\tens c} \tens  \cdots \tens [x^{(2)},y^{(3)}] \tens [x^{(1)},y^{(2)}]^{\tens c} \tens [x_-(j)^+,  y^{(1)}]^{\tens c} \\
\epito [x^{(t)},w]^{\tens c} \tens [x^{(t-1)},y^{(t)}]^{\tens c} \tens \cdots \tens [x^{(2)},y^{(3)}]^{\tens c} \tens [x^{(1)},y^{(2)}]^{\tens c} \tens [x_-(j)^+, y^{(1)}]^{\tens c} \\
\simeq  M^{\Vi_j(a)} \tens M^{\Vi_j(d)},
\ea
\eneqn }}

\noi
and
if $[u,x(j)^+] \in \F$ for some $u<x(j)^+$, then 
 \eqn
&& M^{\Vo_j} \tens [x_-(j)^+,  x(j)^-]^{\tens c}  =M^{\Vo^e_j} \tens [x_-(j)^+,  x(j)^-]^{\tens c}  \\
\simeq  
&&  [x^{(t)},y^{(t)}]^{\tens c} \tens [x^{(t-1)},y^{(t-1)}]^{\tens c} \tens \cdots \tens [x^{(2)},y^{(2)}]^{\tens c}  \tens [x^{(1)},y^{(1)}]^{\tens c} \tens [x_-(j)^+,  x(j)^-]^{\tens c} \\
\epito  && 
 [x^{(t)},y^{(t)}]^{\tens c} \tens [x^{(t-1)},y^{(t-1)}]^{\tens c}\tens \cdots \tens [x^{(2)},y^{(2)}]^{\tens c} \tens [x^{(1)},x(j)^-]^{\tens c} \tens [x_-(j)^+,  y^{(1)}]^{\tens c}\\
\epito  && 
 [x^{(t)},y^{(t)}]^{\tens c} \tens [x^{(t-1)},y^{(t-1)}]^{\tens c} \tens  \cdots \tens [x^{(2)},x(j)^-]^{\tens c} \tens [x^{(1)},y^{(2)}]^{\tens c} \tens [x_-(j)^+, y^{(1)}]^{\tens c} \\
\cdots  && \\
\epito&&  [x^{(t)},y^{(t)}]^{\tens c} \tens [x^{(t-1)},x(j)^-]^{\tens c} \tens  \cdots \tens [x^{(2)},y^{(3)}]^{\tens c} \tens [x^{(1)},y^{(2)}]^{\tens c} \tens [x_-(j)^+,  y^{(1)}]^{\tens c} \\
\epito &&  [x^{(t)},x(j)^-]^{\tens c} \tens [x^{(t-1)},y^{(t)}]^{\tens c} \tens  \cdots \tens [x^{(2)},y^{(3)}] \tens [x^{(1)},y^{(2)}]^{\tens c} \tens [x_-(j)^+,  y^{(1)}]^{\tens c} \\
\simeq &&M^{\Vi(b)_j} \tens M^{\Vi(a)_j} \tens M^{\Vi(d)_j} \simeq M^{\Vi_j},
 \eneqn
as desired.
\QED

The lemma below can be proved by the same argument in Lemma \ref{lem:Voxmx}.
\begin{lemma}
 The tensor product $M^{\Vo} \tens S_{x_-}  \tens S_x $ has a simple head. 
\end{lemma}

Since there is an epimorphism
\eqn
(M^{\Vo} \tens S_{x_-}) \tens [x,x]  \epito M^{\Vo}\tens( \tens_{j\in I} [x_-(j)^+,x(j)^-]^{\tens -\sfc_{j,i}}),
\eneqn
we conclude that
\eqn
\mu([x,x])\hconv [x,x]  \simeq M^{\Vi}, 
\eneqn
as desired.

Thus we obtain \eqref{eq:main}.
\subsection{Case:  $[x,x] \leftarrow [x,x_+]$} 

Set
\eqn
\mu([x,x])= S_{x_+} \hconv M^{\Vi}
\eneqn
Then we have
\eqn 
\de(\mu([x,y]),[x,y]) \le \de( M^\Vi,[x,y])+ \de(S_{x_+}, [x,y])\le d_{i_x}.
\eneqn

There exists a non-zero composition of morphisms
\eqn
 \mu([x,x]) \tens [x,x] \monoto   ( M^{\Vi} \tens S_{x_+}  ) \tens {[x,x]}   \epito  M^{\Vi} \tens [x,x_+].
\eneqn

We omit the proofs of the proposition and lemma below.
\begin{proposition}
There is  an epimorphism 
\eqn
( \tens_{j\in I} [x(j)^+,x_+(j)^-]^{\tens -\sfc_{j,i}}) \tens  M^{\Vi} \epito M^{\Vo}.
\eneqn
\end{proposition}

\begin{lemma}
The tensor product $S_x \tens S_{x_+}  \tens M^{\Vi}  $ has a simple head.
\end{lemma}

Since there is an epimorphism
\eqn
 [x,x] \tens  (S_{x_+}\tens M^{\Vi} )   \epito ( \tens_{j\in I} [x(j)^+,x_+(j)^-]^{\tens -\sfc_{j,i}}) \tens  M^{\Vi}
\eneqn
we conclude that
\eqn
 [x,x] \hconv \mu([x,x])  \simeq M^{\Vo}, 
\eneqn
as desired.

Thus we obtain \eqref{eq:main}.

\end{document}